\documentclass[11pt]{amsart}

\usepackage[utf8]{inputenc}
\usepackage[T1]{fontenc}

\usepackage{mathtools}

\usepackage{amsmath,amssymb, amscd, color}
\usepackage{graphicx}
\usepackage{hyperref}
\usepackage{comment}
\usepackage{subcaption}
\usepackage{float}
\usepackage{cleveref}
\usepackage{stmaryrd}
\usepackage{longtable}
\usepackage{colonequals}

\usepackage{todonotes}

\usepackage[theoremfont]{baskervillef}
\usepackage[baskerville, varg]{newtxmath}
\usepackage[scr=boondox,cal=boondoxo,bb=boondox,frak=euler]{mathalfa}

\usepackage[style=alphabetic,natbib=true,giveninits=true,url=false,maxbibnames=99]{biblatex}

\DeclareFieldFormat{title}{\mkbibquote{#1}}
\DeclareFieldFormat[online]{date}{Preprint (#1)}
\DeclareFieldFormat[misc]{date}{Personal communication (#1)}
\AtEveryBibitem{\clearfield{issn}}
\AtEveryCitekey{\clearfield{issn}}

\numberwithin{equation}{section}

\newtheorem{theorem}{Theorem}[section]

\newtheorem{lemma}[theorem]{Lemma}
\newtheorem{corollary}[theorem]{Corollary}

\newtheorem*{claim*}{Claim}

\newtheorem{definition}[theorem]{Definition}

\newtheorem{remark}[theorem]{Remark}

\newcommand{\calL}{\mathcal{L}}
\newcommand{\calM}{\mathcal{M}}

\newcommand{\NN}{\mathbb{N}}

\newcommand{\QQ}{\mathbb{Q}}
\newcommand{\RR}{\mathbb{R}}

\newcommand{\ZZ}{\mathbb{Z}}
\DeclareMathOperator{\aO}{\mathrm{O}}
\DeclareMathOperator{\ao}{\mathrm{o}}
\DeclareMathOperator{\dimH}{\mathrm{dim}_{\mathrm{H}}}
\DeclareMathOperator{\sizes}{\mathsf{s}}
\DeclareMathOperator{\sizer}{\mathsf{r}}

\newcommand{\bfirst}[1]{{#1}^{b}}
\newcommand{\alast}[1]{{#1}_{a}}

\DeclareMathOperator{\sign}{sign}

\makeatletter
\g@addto@macro\bfseries{\boldmath}
\makeatother

\begin{document}
	
	\title{Fractal dimensions of the Markov and Lagrange spectra near $3$}
	
	\author{Harold Erazo}
	\address[Harold Erazo]{IMPA, Estrada Dona Castorina 110, 22460-320, Rio de Janeiro, Brazil}
	\email{harolderaz@gmail.com}
	
	\author{Rodolfo Gutiérrez-Romo}
	\address[Rodolfo Gutiérrez-Romo]{Centro de Modelamiento Matemático, CNRS-IRL 2807, Universidad de Chile, Beauchef 851, Santiago, Chile.}
	\email{g-r@rodol.fo}
	\urladdr{http://rodol.fo}
	
	\author{Carlos Gustavo Moreira}
	\address[Carlos Gustavo Moreira]{SUSTech International Center for Mathematics, Shenzhen, Guangdong, People’s Republic of China; \hfill\break
		IMPA, Estrada Dona Castorina 110, 22460-320, Rio de Janeiro, Brazil}
	\email{gugu@impa.br}
	
	\author{Sergio Romaña}
	\address[Sergio Romaña]{Universidade Federal do Rio de Janeiro, Av. Athos da Silveira Ramos 149, Centro de Tecnologia - Bloco C - Cidade Universitária, Ilha do Fundão, cep 21941-909, Rio de Janeiro, Brasil}
	\email{sergiori@im.ufrj.br}
	
	\thanks{The first author is partially supported by CAPES. The second author was supported by Centro de Modelamiento Matemático (CMM), ACE210010 and FB210005, BASAL funds for centers of excellence from ANID-Chile, the MATH-AmSud 21-MATH-07 grant, and by ANID-Chile through the FONDECYT Iniciación 11190034 grant. The third author is partially supported by CNPq and FAPERJ. The fourth author is partially supported by FAPERJ-Bolsa Jovem Cientista do Nosso Estado No. E-26/201.432/2022.}
	
	\begin{abstract}
		The Lagrange spectrum $\calL$ and Markov spectrum $\calM$ are subsets of the real line with complicated fractal properties that appear naturally in the study of Diophantine approximations. It is known that the Hausdorff dimension of the intersection of these sets with any half-line coincide, that is, $\dimH(\calL \cap (-\infty, t)) = \dimH(\calM \cap (-\infty, t))\equalscolon d(t)$ for every $t \geq 0$. It is also known that $d(3)=0$ and $d(3+\varepsilon)>0$ for every $\varepsilon>0$.
		
		We show that, for sufficiently small values of $\varepsilon > 0$, one has the approximation $d(3+\varepsilon) = 2\cdot\frac{W(e^{c_0}|\log \varepsilon|)}{|\log \varepsilon|}+\aO\left(\frac{\log |\log \varepsilon|}{|\log \varepsilon|^2}\right)$, where $W$ denotes the Lambert function (the inverse of $f(x)=xe^x$) and $c_0=-\log\log((3+\sqrt{5})/2) \approx 0.0383$. We also show that this result is optimal for the approximation of $d(3+\varepsilon)$ by ``reasonable'' functions, in the sense that, if $F(t)$ is a $C^2$ function such that $d(3+\varepsilon) = F(\varepsilon) + \ao\left(\frac{\log |\log \varepsilon|}{|\log \varepsilon|^2}\right)$, then its second derivative $F''(t)$ changes sign infinitely many times as $t$ approaches $0$.
	\end{abstract}
	
	\maketitle
	\markboth{H. ERAZO, R. GUTIÉRREZ-ROMO, C.G. MOREIRA AND S. ROMAÑA}{DIMENSION OF THE LAGRANGE AND MARKOV SPECTRA NEAR $3$}
	
	\section{Introduction} \label{sec:introduction}
	
	\subsection{The Lagrange spectrum}
	
	The Lagrange spectrum is a subset of the real line which appears naturally in the study of Diophantine approximations of real numbers. 
	
	Consider an irrational real number $x \in \RR \setminus \QQ$. By Dirichlet's approximation theorem, there exist infinitely many pairs of integers $p, q$ with $q>0$ satisfying
	\[
	\left|x - \frac{p}{q}\right| < \frac{1}{q^2}.
	\]
	The previous result is not tight. Indeed, Hurwitz's theorem states that the following holds for infinitely many such pairs $p, q$:
	\[
	\left|x - \frac{p}{q}\right| < \frac{1}{\sqrt{5}q^2}.
	\]
	This is the best possible inequality of this type that holds for every irrational number $x$. Indeed, if $x = \frac{1 + \sqrt{5}}{2}$, the constant $\sqrt{5}$ \emph{cannot} be replaced by a larger constant while preserving the existence of infinitely many such pairs $p, q$ for which the corresponding inequality holds. However, for other irrational values of $x$ we may hope for better results. Following this idea, we define $L(x)$ as the supremum of the set of all $\ell > 0$ such that
	\[
	\left|x - \frac{p}{q}\right| < \frac{1}{\ell q^2}
	\]
	holds for infinitely many pairs of integers $p, q$ with $q>0$ (possibly with $L(x) = \infty$). The number $L(x)$ is known as the \emph{Lagrange value} of $x$, and the \emph{Lagrange spectrum} is defined as the set of all finite Lagrange values:
	\[
	\calL = \{ L(x) < \infty \ \mid\ x \in \RR \setminus \QQ \}.
	\]
	By means of the continued fraction expansion of $x$, it is possible to obtain a symbolic-dynamical characterization of the Lagrange spectrum. Indeed, consider the infinite sequence $(c_n)_{n \geq 0}$ such that
	\[
	x = [c_0; c_1, c_2, c_3, \ldots] = c_0 + \cfrac{1}{c_1 + \cfrac{1}{c_2 + \cfrac{1}{c_3 + \ddots}}},
	\]
	that is, $(c_n)_{n \geq 0}$ is the continued-fraction expansion of $x$. It is well-known that
	\[
	x - \frac{p_n}{q_n} = (-1)^n \frac{1}{(\alpha_{n+1} + \beta_{n+1})q_n^2}
	\]
	where we set $\alpha_{n+1} = [c_{n+1}; c_{n+2}, c_{n+3}, \ldots]$, $\beta_{n+1} = [0; c_n, c_{n-1}, \dotsc, c_1]$, and where $p_n/q_n = [c_0; c_1, c_2, \dotsc, c_n]$. It is also known that these {\it convergents} $p_n/q_n$ of the continued-fraction expansion of $x$ are the best rational approximations of $x$ for instance in the following sense: if $p, q$ are integers with $q>0$ and $|x-p/q|<\frac1{2q^2}$, then $p/q=p_n/q_n$ for some $n\in\NN$. From these facts, we obtain the following expression for the Lagrange value of $x$:
	\[
	L(x) = \limsup_{n \to \infty} (\alpha_{n+1} + \beta_{n+1}).
	\]
	If we define $\beta_{n+1}' = [0; c_n, c_{n-1}, \dotsc, c_1, 1, \dotsc, 1, \dotsc]$, we also have that
	\[
	L(x) = \limsup_{n \to \infty} (\alpha_{n+1} + \beta_{n+1}'),
	\]
	since the trailing sequence of $1$'s does not change the value in the limit.
	
	It follows that 
	\[
	\calL = \left\{ \limsup_{n \to \infty} \lambda(\sigma^n(\omega)) \ \mathbin{\Big|}\ \omega \in (\NN^*)^\ZZ \right\},
	\]
	where, for $\omega=(\omega_n)_{n\in \ZZ} \in (\NN^*)^\ZZ$, $\lambda(\omega) = [\omega^+] + [0; \omega^-]$, with $\omega^+ = (\omega_n)_{n \geq 0}$ and $\omega^- = (\omega_{-n})_{n \ge 1}$.
	
	We refer the reader to the expository article by Bombieri \cite{B:markoff_numbers} and to the books by Cusick--Flahive \cite{Cusick-Flahive}, and by Lima--Matheus--Moreira--Romaña \cite{LMMR} for a more detailed account on these constructions.
	
	\subsection{The Markov spectrum}
	
	The Markov spectrum is another fractal subset of the real line which is very closely related to the Lagrange spectrum. Using the symbolic-dynamical definition of the Lagrange spectrum a starting point, it can be defined similarly as
	\[
	\calM = \left\{ \sup_{n \in \ZZ} \lambda(\sigma^n(\omega)) \ \mathbin{\Big|}\ \omega \in (\NN^*)^\ZZ \right\}.
	\]
	We denote by $\mathsf{m}(\omega)=\sup_{n \in \ZZ} \lambda(\sigma^n(\omega))$ the \emph{Markov value} of $\omega \in (\NN^*)^\ZZ$.
	
	This set is also related to some Diophantine approximation problems. Indeed, it encodes the (inverses of) minimal possible values of real indefinite quadratic forms with normalized discriminants (equal to $1$). Nevertheless, throughout this article we will only use the symbolic-dynamical definitions of $\calL$ and $\calM$.
	
	\subsection{Structure of the Lagrange and Markov spectra}
	
	Both the Lagrange and Markov spectra have been intensively studied since the seminal work of Markov \cite{M:formes_quadratiques}. In particular, it is well-known that
	\[
	\calL \cap [0, 3) = \calM \cap [0, 3) = \left\{ \sqrt{5} < \sqrt{8} < \frac{\sqrt{221}}{5} < \dotsb \right\},
	\]
	that is, $\calL$ and $\calM$ coincide below $3$ and consist of a sequence of \emph{explicit} quadratic surds accumulating only at $3$. Moreover, it is also possible to explicitly characterize the sequences $\omega \in (\NN^*)^\ZZ$ associated with Markov values less than or equal to $3$ \cite[Theorem 15]{B:markoff_numbers}.
	
	On the other hand, the behavior of these sets after $3$ remains somewhat mysterious. Indeed, it is known that $\calL \subseteq \calM$ and some authors conjectured that these sets are \emph{equal}; Freĭman disproved this conjecture only in 1968 \cite{F:non-coincidence}. Much more is now known in this regard: the Hausdorff dimension of the complement $\calM \setminus \calL$ lies strictly between $0$ and $1$ \cite{MM:markov_lagrange}.
	
	Even if the previous paragraph suggests that these sets are somewhat different, they are known to coincide before $3$ and after large enough values. Indeed, Hall showed in 1947 that $\calL$ (and thus also $\calM$) contains a half-line $[c, \infty)$ \cite{H:ray}; any such ray is hence known as a \emph{Hall ray}. After several years, Freĭman found the largest Hall ray to be $[c_{\text{F}}, \infty)$, where $c_{\text{F}} \approx 4.5278\ldots$ is an explicit quadratic surd known as Freĭman's constant \cite{F:constant}. These results in turn imply that $\calL$ and $\calM$ coincide starting at $c_{\text{F}}$, so they both contain the half-line $[c_{\text{F}}, \infty)$.
	
	There are more striking similarities between these two sets. In particular, their Hausdorff dimensions coincide when truncated: the third author showed that
	\[
	\dimH(\calL \cap (-\infty, t)) = \dimH(\calM \cap (-\infty, t))
	\]
	for every $t > 0$ \cite{M:geometric_properties_Markov_Lagrange}. Clearly, this result shows that, when studying the Hausdorff dimension of such truncated versions, one can choose to use either $\calL$ or $\calM$.
	
	Let
	\[
	d(t) \colonequals \dimH(\calL \cap (-\infty, t)) = \dimH(\calM \cap (-\infty, t)).
	\]
	
	Moreira also proved \cite{M:geometric_properties_Markov_Lagrange} the following nice formula:
	\[
	d(t)=\min\{1,2\cdot D(t)\},
	\]
	where $D(t)=\dimH(K_t)$, and 
	\begin{multline*}
		K_t=\{[0;c_1,\dots,c_n,\dots]\ \mid\ \text{there exists $(c_{-n})_{n\geq 0}\in(\NN^*)^{\NN}$ such that} \\ 
		[c_k;c_{k+1},\dots,]+[0;c_{k-1},c_{k-2},\dots]\leq t,\forall k\in\ZZ\}.
	\end{multline*}
	
	In fact, he showed \cite[Lemma 2]{M:geometric_properties_Markov_Lagrange} that $D(t)=\overline{\dim}_{\mathrm{B}}(K_t)=\dimH(K_t)$ and $\overline{\dim}_{\mathrm{B}}(K_t+K_t)=\dimH(K_t+K_t)=\min\{1,2\cdot\dimH(K_t)\}$, where $\overline{\dim}_{\mathrm{B}}$ denotes the upper box dimension. Indeed, that lemma states that, given any $\eta>0$, there is a Gauss--Cantor set $K(B)\subseteq K_t$ such that
	\[
	\dimH(K(B))>(1-\eta)\overline{\dim}_{\mathrm{B}}(K_t),
	\]
	so
	\[
	(1-\eta)\overline{\dim}_{\mathrm{B}}(K_t)\leq\dimH(K(B))\leq\dimH(K_t)\leq \overline{\dim}_{\mathrm{B}}(K_t).
	\]
	Letting $\eta\to 0$ shows the first equality. The second equality follows from the fact that
	\[
	\calM\cap(-\infty,t)\subseteq(\NN^*\cap[1,t])+K_t+K_t
	\]
	and the inequalities
	\begin{align*}
		d(t)&=\dimH(\calM\cap(-\infty,t)) \\
		&\leq\dimH(K_t+K_t)\leq\overline{\dim}_{\mathrm{B}}(K_t+K_t)\leq2\cdot\overline{\dim}_{\mathrm{B}}(K_t)=2\cdot\dimH(K_t).
	\end{align*}
	
	\subsection{The Hausdorff dimension near 3}

	The goal of this article is to determine the behavior of $d(t)$ near $t = 3$. By work of the third author \cite{M:geometric_properties_Markov_Lagrange}, we have that $d(t) > 0$ for every $t > 3$. On the contrary, $d(t) = 0$ for every $t \leq 3$, as $\calL \cap (-\infty,3] = \calM \cap (-\infty,3]$ is countable.
	
	Our main objective is to determine the modulus of continuity of $d(t)$ near $3$. The first result we obtained in this direction was the following:
	
	\smallbreak
	
	There exist constants $C_1, C_2 > 0$ such that, for any sufficiently small $\varepsilon > 0$, one has
	\begin{equation}\label{eq:main}
		C_1 \frac{\log |\log \varepsilon|}{|\log \varepsilon|} \leq d(3 + \varepsilon) \leq C_2 \frac{\log |\log \varepsilon|}{|\log \varepsilon|}.
	\end{equation}

	\smallbreak
	
	Let us explain how this partial result is obtained. Our methods are mainly combinatorial and the proofs of the upper and lower bounds on $d(t)$ are done in separate sections. 
	
	To establish the upper bound, we extend some results in Bombieri's article \cite{B:markoff_numbers} to (factors of) sequences with Markov value slightly larger than $3$. In this way, we can analyze the sequences $\omega \in \{1,2\}^\ZZ\subseteq (\NN^*)^\ZZ$ that produce such Markov values; we show that they are not that different from those with Markov value less than or equal to $3$.
	
	To make this more precise, let $\Sigma(t) = \{\omega \in (\NN^*)^\ZZ \ \mid\ \sup_{n \in \ZZ} \lambda(\sigma^n(\omega)) \leq t\}$. We define $\Sigma(t, n)$ to be the set of length-$n$ subwords of sequences in $\Sigma(t)$. We have the following:
	
	\smallbreak
	
	\begin{theorem} \label{thm:equalities}
		There exists a constant $B > 1$ such that
		\[
		\Sigma(3 + B^{-n}, n) = \Sigma(3, n) = \Sigma(3 - B^{-n}, n)
		\]
		for every sufficiently large integer $n$. In fact we can take $B=6^3=216$ and $n \geq 68$.
	\end{theorem}
	
	\smallbreak
	
	The previous theorem can be interpreted as follows: given a bi-infinite word, whose Markov value is exponentially close to $3$ (smaller than $3 + B^{-n} = 3+6^{-3n}$), then its length-$n$ subwords are indistinguishable from those in $\Sigma(3, n)$. That is to say, a length-$n$ window cannot detect the patterns of symbols that make their Markov values different from $3$; they are only present when considering windows of larger lengths.
	
	Since the words before 3 are well understood, we will construct alphabets that allow us to write words in $\Sigma(3+B^{-n},n)$ as \emph{weakly renormalizable words} (see \Cref{def:weaklyrenormalizable}). This construction is inspired by the ``exponent-reducing'' construction by Bombieri, which is detailed in \Cref{sec:bombieri}. Indeed, the inductive procedure of reducing exponents can also be regarded as replacing the alphabet in which a word is written with a more complicated alphabet (so some exponents are ``captured'' by the letters of the new alphabet). The construction is inductive, so we will develop it as a \emph{renormalization algorithm} (\Cref{lem:renormalizingr}). This algorithm is used to obtain a proof \Cref{thm:equalities}.
	
	\Cref{thm:equalities} allows us to reduce the proof of the upper bound to a simple counting. Indeed, we show in \Cref{lem:O(n^3)} that $|\Sigma(3, n)| = \mathrm{O}(n^3)$, which implies that $|\Sigma(3 +B^n, n)| = \mathrm{O}(n^3)$. This is enough to establish the upper bound by covering $K_t$ with small intervals in the standard way and using this counting.
	
	To show that the lower bound holds, we prove that $d(3+e^{-r})$ (where $r\in\NN^\ast$) is larger than the Hausdorff dimension of a suitable Gauss--Cantor set; recall that a Gauss--Cantor set is a subset of the real line defined by numbers with continued-fraction expansions that obey certain patterns. Finally, the Hausdorff dimension of a Gauss--Cantor set can be estimated by the (relatively elementary) methods in the book by Palis--Takens \cite[Chapter 4]{PalisTakens}, and, hence the proof of \eqref{eq:main} is complete.
	
	While these methods are enough to prove inequalities \eqref{eq:main}, they are actually sufficient to obtain an asymptotic approximation of $d(t)$. In fact, to prove \eqref{eq:main}, only the results in \Cref{sec:weakly_renormalizable} and (a simplification of the results) in \Cref{sec:lower_bound} are needed.
	
	We will now state our main results, which give more precise estimates of $d(t)$ for $t$ close to $3$. Let $f\colon[-1, +\infty)\to [-e^{-1},+\infty)$ be given by $f(x)=x e^x$ and recall that the Lambert $W$ function is the function $W\colon [-e^{-1},+\infty)\to [-1,+\infty)$ given by $W = f^{-1}$. Our main result is the following:
	\medbreak
	\begin{theorem}\label{thm:mainstrong}
		Let $d(t) = \dimH(\calL \cap [0, t)) = \dimH(\calM \cap [0, t))$. Then, for all sufficiently small $\varepsilon$, we have
		\[
		d(3+\varepsilon)=2\cdot\frac{W(e^{c_0}|\log \varepsilon|)}{|\log \varepsilon|}+\aO\left(\frac{\log |\log \varepsilon|}{|\log \varepsilon|^2}\right),
		\]
		where $c_0=-\log\log((3+\sqrt{5})/2) \approx 0.0383$.
	\end{theorem}
	
	The main idea behind the upper bound of \Cref{thm:mainstrong} is again the construction of alphabets that allow us to write finite subwords of $\Sigma(3+e^{-r})$ as \emph{weakly renormalizable words}. Then, using the fact that windows of sizes comparable to $r$ must have a very similar structure with those before 3 (that are well understood because of the work of Bombieri \cite{B:markoff_numbers}), we can find long forced continuations of finite subwords of size comparable to $r$ of words of $\Sigma(3+e^{-r})$. Here, by \emph{size} we no longer mean the length of a word, but rather the size of the interval it induces by continued fraction expansions. Using the covering of $K_t$ constructed with finite subwords of $\Sigma(3+e^{-r})$, we can control the size of a subcovering by smaller intervals (associated with longer words), depending on the structure of each word, so intervals with few continuations contribute less to the dimension. It turns out that there are some configurations which contribute more than others to the dimension of these sets, namely configurations obtained by alternate concatenations of large blocks of $1$'s with blocks $22$.
	
	To be more precise about this last statement, define the Gauss--Cantor set
	\begin{equation*}
		C_n\colonequals K(\{221^n,1\})=\{[0;\gamma_1,\gamma_2,\dots]\ \mid\ \gamma_i\in\{221^n,1\},\forall i\geq 1\},
	\end{equation*}
	and let $\varepsilon_n\colonequals \max L(C_n)$, so $\varepsilon_n$ is of the order of $((3+\sqrt{5})/2)^{-n}$. We have, from \Cref{thm:mainstrong} and from the proof of its lower bound (\Cref{sec:lower_bound}), that
	\begin{align*}
		d(3+\varepsilon_n)&=2\cdot\frac{W(e^{c_0}|\log \varepsilon_n|)}{|\log \varepsilon_n|}+\aO\left(\frac{\log |\log \varepsilon_n|}{|\log \varepsilon_n|^2}\right) \\
		&=2\cdot\dimH(C_n)+\aO\left(\frac{\log |\log \varepsilon_n|}{|\log \varepsilon_n|^2}\right) \\
		&=\dimH(L(C_n))+\aO\left(\frac{\log |\log \varepsilon_n|}{|\log \varepsilon_n|^2}\right) \\
		&=d(3+\varepsilon_{n-1})+\aO\left(\frac{\log |\log \varepsilon_n|}{|\log \varepsilon_n|^2}\right).
	\end{align*}
	
	One natural follow-up question is if it is possible to find a better approximation of $d(t)$ near $3$. The next theorem shows that this is not possible for ``reasonable'' (or explicit) approximations: for such reasonable approximations, the error term is optimal. We prove the following:
	\medbreak
	\begin{theorem}\label{thm:optimal}
		Let $d(t) = \dimH(\calL \cap [0, t)) = \dimH(\calM \cap [0, t))$. There exists sequences $(x_k), (y_k)$ and constants $0<C_1<C_2$, with $0<C_1\varphi^{-4k}=x_k<\frac{3}{2}x_k<y_k=C_2\varphi^{-4k}$, where $\varphi=(1+\sqrt{5})/2$ is the golden mean, such that
		\[
		d(3+y_k)-d(3+x_k)=\aO\left(\frac{1}{k^2}\right).
		\]
		In particular, if $F$ is a twice continuously-differentiable function satisfying
		\[
		d(3+\varepsilon)=F(\varepsilon)+\ao\left(\frac{\log|\log \varepsilon|}{|\log \varepsilon|^2}\right), 
		\]
		then its second derivative $F''(\varepsilon)$ changes sign infinitely many times as $\varepsilon$ approaches $0$.
	\end{theorem}
	
	In fact, we will prove that $d(3+y_k)-d(3+x_k) = \aO\left(\frac1{k^2}\right) = \ao\left(\frac{\log k}{k^2}\right)$, while
	\[
	\frac{W(e^{c_0}|\log y_k|)}{|\log y_k|}-\frac{W(e^{c_0}|\log x_k|)}{|\log x_k|}>\tilde c\frac{\log k}{k^2},
	\]
	for a positive constant $\tilde c$, which implies that the error term in the approximation of $d(3+\varepsilon)$ by any reasonable function of $\varepsilon$ is at least of the order of $\frac{\log |\log \varepsilon|}{|\log \varepsilon|^2}$. In this sense, $(3+x_k,3+y_k)$ is an ``almost plateau'' for the dimension function $d(t)$ (the variation of $d(t)$ in these intervals is much smaller than the variation of its reasonable approximations). Indeed, we have proven that $d(3+\varepsilon)$ is very well approximated by
	\[
	g_1(\varepsilon)=2\cdot\frac{W(e^{c_0}|\log \varepsilon|)}{|\log \varepsilon|},
	\]
	and that it is also asymptotic to the simpler function $g_2(\varepsilon)=2\cdot \frac{\log|\log \varepsilon|}{|\log \varepsilon|}$. Moreover, given constants $0<c_1<c_2$, we have, that
	\[
	g_j(c_2 \varepsilon)-g_j(c_1\varepsilon)=(2\log(c_2/c_1)+\ao(1))\frac{\log|\log \varepsilon|}{|\log \varepsilon|^2},
	\]
	for $j \in \{1,2\}$, so reasonable functions $\tilde g(\varepsilon)$ which are asymptotic to $g_1$ and $g_2$ should satisfy $\tilde g(c_2 \varepsilon)-\tilde g(c_1 \varepsilon)\ge \log(c_2/c_1)\frac{\log|\log \varepsilon|}{|\log \varepsilon|^2}$ for $\varepsilon>0$ small enough.

	While the estimates in the third author's work \cite{M:geometric_properties_Markov_Lagrange} in principle would allow us to obtain some information regarding the modulus of continuity, those estimates are very far from being optimal (this is particularly true for the upper estimates). Thus, we rely on the methods described above instead of on the general methods in the third author's previous work.
	
	This article is organized as follows: Section 2 contains some preliminary notations and facts that we will use later on. By analyzing the combinatorics of finite words, we develop a renormalization algorithm which we use to prove \Cref{thm:equalities} in Section 3. Using the understanding of finite subwords, we will find large forced extensions which by a delicate analysis of the sizes and counting of them, will give us the upper bound of \Cref{thm:mainstrong} in Section 4. In Section 5 we present the construction and analysis of a suitable Gauss--Cantor set, which allows us to establish the lower bound in \Cref{thm:mainstrong} and, thus, to finish the proof of the main theorem. Finally, we study how the bad cuts produce gaps in their respective Markov values in Section 6, which allow us to prove the optimality of our approximation in \Cref{thm:optimal}.
	
	\noindent {\bf Acknowledgements:} We would like to thank Carlos Matheus and Jamerson Bezerra for helpful conversations about the subject of this paper. We also would like to thank Moubariz Garaev, Harald Helfgott and Lola Thompson for organizing and inviting us to the meeting Number Theory in the Americas/Teoría de Números en América held in Casa Matemática Oaxaca, in 2019, where this work started. 
	
	We also thank the anonymous referees for their helpful and insightful comments that greatly improved the exposition of this article.
	
	\section{Preliminaries} \label{sec:preliminaries}
	
	Our goal is to study the function
	\[
	d(t) \colonequals \dimH(\calL \cap (-\infty, t)) = \dimH(\calM \cap (-\infty, t))
	\]
	near $t = 3$.
	If a sequence $\omega \in (\NN^*)^\ZZ$ contains $3$, then $\lambda(\omega) > 3.52$, which is ``much larger'' than $3$, so we can ignore such sequences. Thus, throughout the entire article, a \emph{word} is made up of letters of the alphabet $\{1, 2\}$. Words can be either finite, infinite or bi-infinite. If $w$ is a finite word, we denote its \emph{length} by $|w|$, that is, the amount of letters $1$ or $2$ that are needed to write $w$.
	
	We will also consider \emph{sections} of words, which consist of a word together with a choice of a \emph{splitting point} marked with a vertical bar. A section of a bi-infinite word can be interpreted as a \emph{shift} of the original word. We usually write sections as $\omega = P^*|Q$, where $P \in (\NN^*)^{\NN^*}$ and $Q \in (\NN^*)^{\NN}$ are an infinite words, and $P^* \in (\NN^*)^{-\NN^*}$ denotes the \emph{transpose} of $P$, that is, $P^*_{-k} = P_k$ for every $k \in \NN^*$.
	
	\subsection{Words in \texorpdfstring{$\Sigma(3)$}{Sigma(3)}} \label{sec:bombieri} Bombieri \cite{B:markoff_numbers} showed that bi-infinite words in $\Sigma(3)$ have to follow very special patterns (which is essentially a restatement of much older results by Markov \cite{Markoff}, as stated in the book by Cusick--Flahive \cite{Cusick-Flahive}). Indeed, he showed \cite[Lemma 9]{B:markoff_numbers} that $\omega$ is a word in the letters $a = 22$ and $b = 11$ (that is, the number of consecutive ones or twos is always even or infinite), and he also showed \cite[Lemma 11]{B:markoff_numbers} that, if $\omega \in \Sigma(3)$, then $\omega$ has to be of four possible forms:
	\begin{itemize}
		\item Constant, that is, $\omega = a^\infty$ or $\omega = b^\infty$;
		\item Degenerate, that is, $\omega = b^\infty ab^\infty$ or $\omega = a^\infty ba^\infty$;
		\item Type I, that is, $\omega = \ldots ab^{e_i}ab^{e_{i+1}}a \ldots$ with every $e_i \geq 1$; or
		\item Type II, that is, $\omega = \ldots ba^{e_i}ba^{e_{i+1}}b \ldots$ with every $e_i \geq 1$.
	\end{itemize}
	The \emph{exponents} $(e_i)_{i \in \ZZ}$ that appear in Type I and Type II elements of $\Sigma(3)$ also have to be of some special forms, but we will not use them explicitly.
	
	Now, let $U$ and $V$ be the Nielsen substitutions given by
	\[
	U \colon \begin{matrix}
		a &\mapsto &ab \\
		b &\mapsto &b
	\end{matrix}, \qquad
	V \colon \begin{matrix}
		a &\mapsto &a \\
		b &\mapsto &ab.
	\end{matrix}
	\]
	This substitutions have inverses defined in the free group $\mathsf{F}\langle a, b\rangle$ given by
	\[
	U^{-1} \colon \begin{matrix}
		a &\mapsto &ab^{-1} \\
		b &\mapsto &b
	\end{matrix}, \qquad
	V^{-1} \colon \begin{matrix}
		a &\mapsto &a \\
		b &\mapsto &a^{-1}b.
	\end{matrix}
	\]
	Bombieri also proved \cite[Lemma 14]{B:markoff_numbers} that if $\omega \in \Sigma(3)$, then both $U(\omega)$ and $V(\omega)$ belong to $\Sigma(3)$. These words can be described explicitly. Indeed, if we write $\omega = \ldots ab^{e_i}ab^{e_{i+1}}a \ldots$ where each $e_i \geq 0$, then
	\[
	U(\omega) = \ldots ab^{e_i+1}ab^{e_{i+1}+1}a \ldots
	\]
	Similarly, if we write $\omega = \ldots ba^{e_i}ba^{e_{i+1}}b \ldots$ with each $e_i \geq 0$, then
	\[
	V(\omega) = \ldots ba^{e_i+1}ba^{e_{i+1}+1}b \ldots
	\]
	
	Furthermore, we have that if $\omega$ is of Type I, then $U^{-1}(\omega)$ is well-defined and belongs to $\Sigma(3)$. Similarly, if $\omega$ is of Type II, then $V^{-1}(\omega)$ is well-defined and belongs to $\Sigma(3)$. These can be described by
	\[
	U^{-1}(\omega) = \ldots ab^{e_i-1}ab^{e_{i+1}-1}a \ldots
	\]
	and
	\[
	V^{-1}(\omega) = \ldots ba^{e_i-1}ba^{e_{i+1}-1}b \ldots,
	\]
	where $e_i \geq 1$ for each $i \in \ZZ$ by definition.
	
	We will now include an useful lemma which is implicit in Bombieri's work.
	
	\begin{lemma}
		\label{lem:bombieri_old}
		A nonempty finite word $w$ belongs to $\Sigma(3,|w|)$ if and only if there exists $W \in \langle U, V\rangle$ such that $w$ is a factor of $W(ab)$.
	\end{lemma}
	\begin{proof} 
		We will first show that if $w$ is a factor of $W(ab)$ for some $W \in \langle U, V\rangle$, then it belongs to $\Sigma(3,|w|)$. This is shown by Bombieri \cite[Theorem 15]{B:markoff_numbers}, as the word $\omega = \ldots W(ab)W(ab)W(ab) \ldots$ belongs to $\Sigma(3)$.
		
		We will now show that if $w \in \Sigma(3, |w|)$, then it is a factor of $W(ab)$ for some $W \in \langle U, V\rangle$.
		
		Let $\omega$ be a bi-infinite word in $\Sigma(3)$ containing $w$ as a factor. We know that $\omega$ can only be constant, degenerate, of Type I or of Type II.
		
		Assume first that $\omega$ is constant. Then, $w$ is a factor of $a^k$ or $b^k$ for some $k \geq 1$. Observe that $U^{k-1}(ab) = ab^k$ and $V^{k-1}(ab) = a^kb$, so the result follows in this case. Assume now that $\omega$ is degenerate. If $w$ is constant, we reduce to the previous case. Otherwise, $w$ is factor of $b^k a b^k$ or $a^k b a^k$ for some $k \geq 0$. Since $U^kV(ab) = ab^kab^{k+1}$ and $V^kU(ab) = a^{k+1}ba^kb$, we also obtain the result in this case.
		
		We will then suppose that $\omega$ is of Type I or of Type II. Hence, $U^{-1}(\omega)$ or $V^{-1}(\omega)$ is well-defined and belongs to $\Sigma(3)$. Recall that these automorphisms act by reducing all exponents by $1$.
		
		By iteratively applying the appropriate automorphism, $U^{-1}$ or $V^{-1}$, we obtain a (possibly finite) sequence of bi-infinite words $\omega = \omega^{(1)}, \omega^{(2)}, \ldots$ This process only stops if $\omega^{(n)}$ is constant or degenerate for some $n \in \mathbb{N}^*$. In this latter case, take $W' \in \langle U, V \rangle$ such that $\omega = W'(\omega^{(n)})$. By definition, there exists a factor $\theta$ of $\omega^{(n)}$ such that $W'(\theta)$ contains $w$. Since $\omega^{(n)}$ is constant or degenerate, we know that its factors satisfy the statement of the lemma. Thus, $\theta$ is contained in a word of the form $W''(ab)$ for some $W'' \in \langle U, V\rangle$. We obtain that $w$ is then a factor of $W'W''(ab)$.
		
		Finally, assume that the process never stops, so we obtain an infinite sequence $(\omega^{(k)})_{k \in \mathbb{N}^*}$ of bi-infinite words. Possibly by first making $w$ longer so it can be written in the alphabet $\{a, b\}$, we can apply the same sequence of operations to the finite word $w$, that is, reduce its exponents by $1$ in the same way that the exponents of the bi-infinite words in the sequence $(\omega^{(k)})_{k \in \mathbb{N}^*}$ are being reduced by $1$. In this way, we obtain a sequence $w = w^{(1)}, w^{(2)}, \ldots$ of (possibly empty) finite words. We claim that $w^{(n)}$ is constant and nonempty for some $n \geq 1$. Indeed, if $w^{(k)}$ is not constant for some $k \geq 1$, then $w^{(k+1)}$ is nonempty, since only the exponents of exactly one of the letters, $a$ or $b$, are reduced by the operation taking $w^{(k)}$ to $w^{(k+1)}$. Moreover, we have that $|w^{(k+1)}| < |w^{(k)}|$, as some exponents are reduced. If $w^{(k+1)}$ is again not constant, we can continue the process inductively. Since $w$ is finite, this process has to stop, so some word in the sequence must be constant. This completes the proof as we already know that constant words satisfy the statement of the lemma.
	\end{proof}
	
	\subsection{Constraints for words} For a finite word $u = u_1 u_2 \ldots u_n \in (\mathbb N^*)^n$, we define $u^*$ as the transpose of $u$, that is, $u^* = u_n u_{n-1} \ldots u_1 \in (\mathbb{N}^*)^n$. Moreover, we set
	\[
	M_u = \begin{cases}
		12 & |u| \text{ is even} \\
		21 & |u| \text{ is odd}
	\end{cases}, \qquad
	m_u = \begin{cases}
		21 & |u| \text{ is even} \\
		12 & |u| \text{ is odd}.
	\end{cases}
	\]
	
	Now, given a section $w = u^*|v$ of a finite word $w$, we define
	\[
	\lambda^+(w) = [v m_u^\infty] + [0u M_v^\infty], \qquad \lambda^-(w) = [v M_u^\infty] + [0u m_v^\infty].
	\]
	These quantities are the largest and smallest values of $\lambda$ that a section of a bi-infinite word containing $w$ can attain, respectively. Thus, they induce restrictions on which finite words can be factors of bi-infinite words whose Markov values are known to be bounded in some way.
	
	\subsection{Useful notation} We will set some notation that will be used through all the article; some of it was borrowed from the third author's work \cite{M:geometric_properties_Markov_Lagrange}.
	
	For a finite word $\alpha \in (\mathbb N^*)^n$ written as $\alpha = c_1 c_2 \ldots c_n$, we define its \emph{size} by $\sizes(\alpha)\colonequals |I(\alpha)|$, where $I(\alpha)$ is the interval
	\[
	I(\alpha) \colonequals \{x\in[0,1]\ \mid\ x=[0; c_1,c_2,\dots,c_n,t], t\ge 1\} \cup \{[0, c_1,c_2, \dotsc, c_n]\}
	\]
	consisting of the numbers in $[0,1]$ whose continued fractions start with $\alpha$. The set $I(\alpha)$ is a closed interval in $[0,1]$.
	
	If we take $p_0=0$, $q_0=1$, $p_1=1$, $q_1=c_1$ and, for each integer $k\ge 0$, we take $p_{k+2}=c_{k+2}p_{k+1}+p_k$ and $q_{k+2}=c_{k+2} q_{k+1}+q_k$, then the endpoints of $I(\alpha)$ are $[0;c_1,c_2,\dots,c_n]=p_n/q_n$ and $[0;c_1,c_2,\dots,c_{n-1},c_n+1]=\frac{p_n+p_{n-1}}{q_n+q_{n-1}}$. Thus,
	\[
	\sizes(\alpha)=\left| \frac{p_n}{q_n}-\frac{p_n+p_{n-1}}{q_n+q_{n-1}} \right| = \frac1{q_n(q_n+q_{n-1})},
	\]
	since $p_nq_{n-1}-p_{n-1}q_n=(-1)^{n-1}$. We define $\sizer(\alpha)=\lfloor\log (\sizes(\alpha)^{-1}) \rfloor$, which controls the order of magnitude of the size of $I(\alpha)$. Observe that $\sizer(\alpha) \leq r$ if and only if $\sizes(\alpha) > e^{-r-1}$.
	
	We also define, for $r \in \mathbb N$, the set
	\[
	Q_r= \{\alpha=c_1 c_2 \ldots c_n \ \mid\ \sizer(\alpha) \ge r, \sizer(c_1 c_2 \ldots c_{n-1})<r\}.
	\]
	Observe that $\alpha \in Q_r$ if and only if $\sizes(\alpha) < e^{-r-1}$ and $\sizes(\alpha') \geq e^{-r-1}$, where $\alpha'$ is the word obtained by removing the last letter from $\alpha$. Informally, this means that the interval $I(\alpha)$ is ``small'', while the interval $I(\alpha')$ is ``not as small'', so the last letter of cannot be removed from $\alpha$ without changing the order of magnitude of $|I(\alpha)|$.
	
	Let us recall some estimates from the third author's work \cite{M:geometric_properties_Markov_Lagrange} that will be useful for us. Indeed, for any finite words $\alpha$, $\beta$, we have that
	\[
	\frac{1}{2} \sizes(\alpha)\sizes(\beta)<\sizes(\alpha\beta)<2\sizes(\alpha)\sizes(\beta);
	\]
	it follows that $\sizer(\alpha)+\sizer(\beta)-1\le \sizer(\alpha\beta)\le \sizer(\alpha)+\sizer(\beta)+2$ \cite[Lemma A.2]{M:geometric_properties_Markov_Lagrange}. By Euler's property of continuants (\Cref{lem:eulersproperty}), if $\alpha=c_1c_2\cdots c_m$ and $\beta=d_1 d_2\cdots d_n$ are finite words, then we have
	\[
	q_{m+n}(\alpha\beta)=q_m(\alpha)q_n(\beta)+q_{m-1}(c_1c_2\cdots c_{m-1})q_{n-1}(d_2 d_3\cdots d_n),
	\]
	and, thus,
	\[
	q_m(\alpha)q_n(\beta)<q_{m+n}(\alpha\beta)<2q_m(\alpha)q_n(\beta).
	\]
	
	Finally, recall that $\Sigma(t) = \{\omega \in (\NN^*)^\ZZ \ \mid\ \sup_{n \in \ZZ} \lambda(\sigma^n(\omega)) \leq t\}$ and that $\Sigma(t, n)$ is the set of length-$n$ subwords of sequences in $\Sigma(t)$. In this context, we define $\Sigma^{(r)}(3+\delta)$ as the set of the words $w\in Q_r$ belonging to $\Sigma(3+\delta,|w|)$.
	
	\section{Weakly renormalizable words} \label{sec:weakly_renormalizable}
	
	The main goal of this section is to prove \Cref{thm:equalities}. For this task, we will prove several lemmas that allow us to understand the structure of  $\Sigma(3,n)$.
	
	\subsection{Basic facts about \texorpdfstring{$\lambda$}{lambda}} We start by showing some basic facts about the function $\lambda$ that will be useful throughout the article.
	
	\begin{lemma}\label{lem:121_212_forbidden}
		Let $\omega \in \Sigma(3.06)$. Then, $\omega$ does not contain $121$ or $212$ as subwords.
	\end{lemma}
	\begin{proof}
		Assume that $\omega \in \Sigma(3.06)$. Observe that $\lambda^-(1|21) > 3.15$, so the word $121$ does not appear in $\omega$. Now, if $212$ is a subword of $\omega$, so is $2212$. This is not possible since $\lambda^-(2|212) > 3.06$.
	\end{proof}

	\begin{lemma} \label{lem:calc_s}
		Let $\omega$ be a bi-infinite word in $1$ and $2$ not containing $121$ and $212$ and such that $\omega = R^* w^*b|awS$, where $w$ is a finite word, $R = R_1 R_2 \ldots$, $S = S_1 S_2 \ldots$ and $R_1 \neq S_1$, with $R_i, S_i \in \{1,2\}$ for each $i$. Then
		\[
		\sizes(bwb)< \sign([w,S]-[w,R])(\lambda(\omega)-3)<\sizes(bw1).
		\]
		In particular if $w$ has even length, $R_1=1$ and $S_1=2$, then 
		\[
		\sizes(bwb)< \lambda(\omega)-3<\sizes(bw1).
		\]
	\end{lemma}
	
	\begin{proof}
		First observe that $[2; 2, w, R] + [0; 1, 1, w, R] = 3$. Thus, we have that
		\begin{align*}
			\lambda(R^*w^* 11| 22w S)&=[2;2,w,S]+[0;1,1,w,R]\\
			&=3+[0;1,1,w,R]-[0;1,1,w,S].
		\end{align*}
		We obtain that
		\begin{align*}
			&\lambda(R^*w^* 11| 22w S) - 3 \\
			{}={} &[0;1,1,w,R]-[0;1,1,w,S] \\
			{}={} &\sign([w,S]-[w,R])\cdot |[0;1,1,w,R]-[0;1,1,w,S]|.
		\end{align*}
		
		Let $x = [0;1,1,w,R]$ and $y = [0;1,1,w,S]$. We will write the continued-fraction expansion of these numbers as
		\begin{align*}
			x &= [0;u_1, u_2,\dotsc, u_{\ell}, u_{\ell+1}, u_{\ell+2}, \ldots], \\
			y &= [0;u_1, u_2,\dotsc, u_{\ell}, v_{\ell+1}, v_{\ell+2}, \ldots],
		\end{align*}
		where $u_1 = u_2 = 1$ and $u_{\ell+1} \neq v_{\ell+1}$. With this notation, we have
		\[
		\sign([w,S]-[w,R])(\lambda(\omega) - 3) = |x - y|.
		\]

		Let $(p_n/q_n)_{n \in \mathbb N}$ be the sequence of convergents of $x$. More explicitly, we have that $p_n/q_n = [0;u_1,u_2,\dots,u_n]$.
		
		If we put $\alpha_{\ell+1} = [u_{\ell+1}; u_{\ell+2}, u_{\ell+3}, \ldots]$, then
		\[
		x=[0; u_1, u_2, \dots, u_{\ell}, \alpha_{\ell+1}] = \frac{\alpha_{\ell+1}p_{\ell}+p_{\ell-1}}{\alpha_{\ell+1}q_{\ell}+q_{\ell-1}}.
		\]
		Similarly, let $\beta_{\ell+1} = [v_{\ell+1}; v_{\ell+2}, v_{\ell+3}, \ldots]$. We then have that
		\[
		y=\frac{\beta_{\ell+1}p_{\ell}+p_{\ell-1}}{\beta_{\ell+1}q_{\ell}+q_{\ell-1}},
		\]
		since the sequence of convergents of $y$ coincides with $(p_n/q_n)_{n \in \NN}$ up to $n = \ell$. Thus,
		\begin{align}
			|x-y| &= \left|\frac{\alpha_{\ell+1} p_{\ell}+p_{\ell-1}}{\alpha_{\ell+1} q_{\ell}+q_{\ell-1}}-\frac{\beta_{\ell+1} p_{\ell}+p_{\ell-1}}{\beta_{\ell+1} q_{\ell}+q_{\ell-1}}\right| \notag \\
			&= \left|\frac{(\alpha_{\ell+1}-\beta_{\ell+1})(p_{\ell}q_{\ell-1}-p_{\ell-1}q_{\ell})}{(\alpha_{\ell+1}q_{\ell}+q_{\ell-1})(\beta_{\ell+1} q_{\ell}+q_{\ell-1})}\right| \notag \\
			&= \left|\frac{(\alpha_{\ell+1}-\beta_{\ell+1})(-1)^{\ell-1}}{(\alpha_{\ell+1}q_{\ell}+q_{\ell-1})(\beta_{\ell+1}q_{\ell}+q_{\ell-1})}\right| \notag \\
			&= \frac{|\alpha_{\ell+1}-\beta_{\ell+1}|}{(\alpha_{\ell+1}q_{\ell}+q_{\ell-1})(\beta_{\ell+1} q_{\ell}+q_{\ell-1})}, \label{eq:x-y_1}
		\end{align}
		where we used that $p_{\ell}q_{\ell-1}-p_{\ell-1}q_{\ell} = (-1)^{\ell + 1}$.
		
		Since we are only interested in continued fractions whose partial quotients are $1$ or $2$, we can assume, without loss of generality, that $u_{\ell+1}=2$ and $v_{\ell+1}=1$. We denote $\alpha=\alpha_{\ell+1}$, $\beta=\beta_{\ell+1}$ and $\lambda=q_{\ell-1}/q_{\ell} \in (0,1)$. Thus,
		\begin{equation} \label{eq:x-y_2}
			|x-y|=\frac{\alpha-\beta}{q_{\ell}^2(\alpha+\lambda)(\beta+\lambda)}=\frac1{q_{\ell}^2}\left(\frac1{\beta+\lambda}-\frac1{\alpha+\lambda}\right).
		\end{equation}
		We obtain that $|x - y|$ is (for fixed $q_{\ell-1}$ and $q_{\ell}$) an increasing function of $\alpha$, and a decreasing function of $\beta$. By analyzing \Cref{eq:x-y_1,eq:x-y_2} we deduce that:
		\begin{itemize}
			\item the quantity $|x - y|$ is minimized when $\alpha$ is minimized, and $\beta$ is maximized. This happens when
			\begin{align*}
				\alpha=\alpha_0 &\colonequals [2;\overline{2,1,1,1,2,2}] = \frac{21 + 2\sqrt{210}}{21} \approx 2.3801, \\
				\beta=\beta_0 &\colonequals [1;\overline{1,2,2,2,1,1}] = \frac{6 + \sqrt{210}}{12} \approx 1.7076.
			\end{align*}
			
			\item the quantity $|x - y|$ is maximized when $\alpha$ is maximized, and $\beta$ is minimized. This happens when
			\begin{align*}
				\alpha=\alpha_1 &\colonequals [2;\overline{1,1,1,2,2,2}] = \frac{21 + 2\sqrt{210}}{19} \approx 2.6306, \\
				\beta=\beta_1 &\colonequals [1;\overline{2,2,2,1,1,1}] = \frac{12 + 2\sqrt{210}}{29} \approx 1.4132,
			\end{align*}
		\end{itemize}

		On the other hand, $bwb= u_1 u_2 \ldots u_{\ell} 1 1$, so 
		\begin{align*}
			\sizes(bwb) &=|[0; u_1, u_2, \dots, u_{\ell},1,1]-[0; u_1, u_2, \dots, u_{\ell},1,1,1]| \\
			&=\left|\frac{2p_{\ell}+p_{\ell-1}}{2q_{\ell}+q_{\ell-1}}-\frac{3p_{\ell}+2p_{\ell-1}}{3q_{\ell}+2q_{\ell-1}}\right| \\
			&=\frac1{(2q_{\ell}+q_{\ell-1})(3q_{\ell}+2q_{\ell-1})}=\frac1{q_{\ell}^2}\frac1{(2+\lambda)(3+2\lambda)}.
		\end{align*}
		
		Similarly
		\begin{align*}
			\sizes(bw1)&=|[0; u_1, u_2, \dots, u_{\ell},1]-[0; u_1, u_2, \dots, u_{\ell},1,1]| \\ &=\left|\frac{p_{\ell}+p_{\ell-1}}{q_{\ell}+q_{\ell-1}}-\frac{2p_{\ell}+p_{\ell-1}}{2q_{\ell}+q_{\ell-1}}\right|
			\\
			&=\frac1{(q_{\ell}+q_{\ell-1})(2q_{\ell}+q_{\ell-1})}=\frac1{q_{\ell}^2}\frac1{(1+\lambda)(2+\lambda)}.
		\end{align*}

		We then have 
		\begin{align*}\frac{|x-y|}{\sizes(bwb)}&\ge (\alpha_0-\beta_0)\frac{(2+\lambda)(3+2\lambda)}{(\alpha_0+\lambda)(\beta_0+\lambda)}\\
			&\ge(\alpha_0-\beta_0)\frac{(2+1/3)(3+2/3)}{(\alpha_0+1/3)(\beta_0+1/3)}\approx 1.03895>1,
		\end{align*}
		since the maps $f_1(\lambda)=\frac{2+\lambda}{\alpha_0+\lambda}$ and $f_2(\lambda)=\frac{3+2\lambda}{\beta_0+\lambda}$ are increasing and
		\[
		\lambda=q_{\ell-1}/q_{\ell}=q_{\ell-1}/(u_{\ell}q_{\ell-1}+q_{\ell-2})\ge q_{\ell-1}/(2q_{\ell-1}+q_{\ell-2})\ge 1/3.
		\]
		Analogously,
		\begin{align*}
			\frac{|x-y|}{\sizes(bw1)}&\le (\alpha_1-\beta_1)\frac{(1+\lambda)(2+\lambda)}{(\alpha_1+\lambda)(\beta_1+\lambda)}\\
			&\le(\alpha_1-\beta_1)\frac{(1+1)(2+1)}{(\alpha_1+1)(\beta_1+1)}\approx 0.83374<1,
		\end{align*}
		since the maps $g_1(\lambda)=\frac{1+\lambda}{\beta_1+\lambda}$ and $g_2(\lambda)=\frac{2+\lambda}{\alpha_1+\lambda}$ are increasing and $\lambda\le 1$.
	\end{proof}
	
	\begin{remark}
		The Markov value of $\omega = R^*11|22S$ coincides with the Markov value of $\sigma(\omega)^* = S^*2|211R$ \cite[Lemma 5]{B:markoff_numbers}.
	\end{remark}
	
	It is not difficult to adapt the proof above to obtain a more explicit (but weaker) version of this lemma which depends only on the length of $w$:
	
	\begin{lemma} \label{lem:R_S_bound}
		Let $\omega$ be a bi-infinite word in $1$ and $2$ not containing $121$ and $212$ and such that $\omega = R^* 11|22 S$ with $R = R_1 R_2 \ldots$ and $S = S_1 S_2 \ldots$ and $R \neq S$. Let $\ell$ be the smallest nonnegative integer such that $R_\ell \neq S_\ell$. Then,
		\[
		\frac{1}{7} (3 - 2\sqrt{2})^{\ell} < \sign([S]-[R])(\lambda(R^* 11| 22 S) - 3 ) < \frac{1}{7} \left(\frac{3-\sqrt{5}}2\right)^{\ell}.
		\]
		In particular, if $w = w^*$ and $\ell = |w|$ is even, then
		\[
		3-\frac{1}{7}\left(\frac{3-\sqrt{5}}2\right)^{\ell+1}
		< \lambda( (wba)^\infty w b| a w (baw)^\infty) < 3 -\frac{1}{7}(3-2\sqrt{2})^{\ell+1} \].
	\end{lemma}

	We will usually use the previous lemma in the following way. Consider a finite word $w$ in the alphabet $\{a, b\}$. Assume that $ba$ is a factor of $w$. Then, we write $w = u^*b|av$, where the vertical bar indicates a \emph{cut}, that is, the position at which we compute the Markov value. Now, let $\ell$ be the smallest nonnegative integer such that $u_\ell \neq v_\ell$ and assume that $u_\ell = b$ and $v_\ell = a$. In other words, $w$ contains the factor $b\theta^*b|a\theta a$, where the vertical bar marks the same position as the cut in $w$. By the previous lemma, the Markov value of \emph{any} infinite word in the alphabet $\{a, b\}$ containing $w$ is at least $3+\frac{1}{7}(3 - 2\sqrt{2})^{2\ell-1}$. Similarly, if $w$ contains $ab$ as a factor, then we can also write $w = u^*a|bv$. Assume now that the smallest nonnegative integer $\ell$ such that $u_\ell \neq v_\ell$ satisfies $u_\ell = a$ and $v_\ell = b$. Then, the Markov value of \emph{any} infinite word in the alphabet $\{a, b\}$ containing $w$ is at least $3+\frac{1}{7}(3 - 2\sqrt{2})^{2\ell-1}$.
	
	In particular, if we assume that $\omega$ is an infinite word in the alphabet $\{a, b\}$ and that its Markov value is sufficiently small, then no finite factor $w$ of $\omega$ can contain patterns as above. This ultimately allows us deduce that some letters are \emph{forced} inside an infinite word containing a finite word.
	
	For the sake of concreteness, we will demonstrate an usage of the previous lemma by showing that no bi-infinite word in $\Sigma(3.0007)$ contains the factor $w = bbab|aa$. Let $\omega$ be a bi-infinite word containing $w$. We start by considering the cut $bb|abaa$. By the previous lemma, if $aa$ does not appear at the left of $w$ in $\omega$, then $\lambda(\omega) > 3+\frac{1}{7}(3 - 2\sqrt{2})^3 > 3.0007$. Thus, we assume that $\omega$ contains $aabbabaa$ as a factor. We can now consider consider a second cut, $aa|bbabaa$. This cut shows that $\lambda(\omega) > 3+\frac{1}{7}(3 - 2\sqrt{2})^1 > 3.0007$, which completes the example.

	We now show that sequences of $1$'s or $2$'s of odd length are forbidden if we assume that the Markov value of a word is sufficiently close to $3$ (relative to the size of the interval it defines).
	\begin{lemma}\label{lem:1}
		Let $r\in\NN$ with $r\geq 5$. Let $c,c'\in\{1,2\}$ with $c\neq c'$. Let $w = c' c^n c'$, for some integer $n \geq 1$, and suppose that $w\in\Sigma(3+e^{-r}, |w|)$. If $\sizer(c^n)\leq r-4$ then $n$ is even. 
	\end{lemma}
	
	\begin{proof}
		Note that $w \neq 121$ and $w \neq 212$ by \Cref{lem:121_212_forbidden}, so $n > 1$. Without loss of generality, we can assume that $w$ is the shortest word of this form satisfying $w \in \Sigma(3+e^{-r}, |w|)$. Let $\omega \in  \Sigma(3+e^{-r})$ be a bi-infinite word such that $w$ is a factor of $\omega$. Assume by contradiction that $n=2k+1$. We will show that $\lambda(\omega) > 3 + e^{-r}$.
		
		Suppose $c=1$. We have a section $\omega = R^*11|22S$ with
		\begin{align*}
			R &= R_1R_2R_3 \ldots = 1^{2k-1}2\ldots \\
			S &= S_1S_2S_3 \ldots = 2^p 1^q 2^r\ldots
		\end{align*}
		By \Cref{lem:calc_s}, $p > 0$ implies that $\lambda(\omega) > 3 + \sizes(bb) = 3 + \frac{1}{40}$, which contradicts the assumption on $w$. Thus, we have that $p = 0$. Let $\ell$ be the smallest positive integer such that $R_\ell \neq S_\ell$. We have two cases:
		\begin{itemize}
			\item If $q > 2k - 1$, then $\ell = 2k$. Since we are assuming that $n = 2k + 1$, we have that $\ell < n$. Moreover, we have that $[S] > [R]$ since $S_\ell < R_\ell$ and $\ell$ is even.
			\item If $q \leq 2k - 1$, then it is even as, otherwise, it would contradict the assumption on $k$. Thus, $q \leq 2k - 2$ and $\ell = q + 1 < n$. Hence, we have that $[S] > [R]$ as $S_\ell > R_\ell$ and $\ell$ is odd.
		\end{itemize}
		In any case, by the assumption on $n$ we obtain from \Cref{lem:calc_s} that
		\[
		\lambda(\omega) > 3 + \sizes(111^{\ell-1}11) \geq 3 + \sizes(1^{n+3}) \geq 3+\sizes(1^{n})e^{-3}>3+e^{-r},
		\]
		where the last inequality holds as $\sizer(1^n) \leq r - 4$.
		
		Now suppose $c=2$, so we have a section $\omega = R^* 11|22 S$ with
		\begin{align*}
			R &= R_1R_2R_3\ldots = 1^p 2^q 1^r\ldots \\
			S &= S_1S_2S_3\ldots = 2^{2k - 1} 1 \ldots
		\end{align*}
		If $p > 0$, \Cref{lem:calc_s} shows that $\lambda(\omega) > 3 + \sizes(bb) = 3 + \frac{1}{40}$, so we have that $p = 0$. Let $\ell$ be the smallest positive integer such that $R_\ell \neq S_\ell$. We have two cases:
		\begin{itemize}
			\item If $q > 2k - 1$, then $\ell = 2k$. Since we are assuming that $n = 2k + 1$, we have that $\ell < n$. Moreover, we have that $[S] > [R]$ since $S_\ell < R_\ell$ and $\ell$ is even.
			\item If $q \leq 2k - 1$, then it is even as, otherwise, it would contradict the assumption on $k$. Thus, $q \leq 2k - 2$ and $\ell = q + 1 < n$. Hence, we have that $[S] > [R]$ as $S_\ell > R_\ell$ and $\ell$ is odd.
		\end{itemize}
		In any case, by the assumption on $n$ we obtain from \Cref{lem:calc_s} that
		\[
		\lambda(\omega) > 3 + \sizes(112^{\ell-1}11) \geq 3 + \sizes(2^{n})e^{-2}>3+e^{-r}
		\]
		where the last inequality holds as $\sizer(2^n) \leq r - 3$.
	\end{proof}	
	
	Whenever we want a version of some lemma that depends only on the length of a word instead of on the size of the interval that it defines (since we want to prove \Cref{thm:equalities} which is stated in terms of lengths of words), we can either repeat the proof using \Cref{lem:R_S_bound} instead of \Cref{lem:calc_s}, or directly compare $\sizer$ with the length using \Cref{lem:rcomparedwithsize}. For example, we can show that sequences of $1$'s or $2$'s of odd length are forbidden:
	
	\begin{lemma} \label{lem:no_odd_ones_or_twos}
		Let $n$ be sufficiently large so that
		\[
		\frac{1}{6^n} < \frac{1}{7}(3-2\sqrt{2})^n;
		\]
		for the sake of concreteness, we can take $n \geq 68$. 
		Let $\omega \in \Sigma(3 + 6^{-n})$. Then, $\omega$ does not contain $12^{2k+1}1$ or $21^{2k+1}2$ as subwords if $2k + 1 < n$.
	\end{lemma}
	
	\subsection{Nielsen substitutions and sequences with Markov value close to 3}
	
	Recall the Nielsen substitutions
	\[
	U \colon \begin{matrix}
		a &\mapsto &ab \\
		b &\mapsto &b
	\end{matrix}, \qquad
	V \colon \begin{matrix}
		a &\mapsto &a \\
		b &\mapsto &ab.
	\end{matrix}
	\]
	Let $T$ be the tree obtained by successive applications of the substitutions $U$ and $V$, starting at the root $ab$. Let $P$ be the set of vertices of $T$ and let $P_n$, for $n \geq 0$, be the set of elements of $P$ that whose distance to the root $ab$ is exactly $n$. Recall from \Cref{lem:bombieri_old} that a finite word $w$ belongs to $\Sigma(3, |w|)$ if and only if it is a factor of a word in $P$.
	
	Given a pair of words $(u, v)$, we also define the operations $\overline{U}(u, v) = (uv, v)$ and $\overline{V}(u, v) = (u, uv)$. Let $\overline{T}$ be the tree obtained by successive applications of the operations $\overline{U}$ and $\overline{V}$, starting at the root $(a, b)$. Let $\overline{P}$ be the set of vertices of $\overline{T}$ and let $\overline{P}_n$, for $n \geq 0$, be the set of elements of $\overline{P}$ that whose distance to the root $(a, b)$ is exactly $n$.
	
	Let $\mathsf{c}$ be the concatenation operator, that is, $\mathsf{c}(u, v) = uv$.
	
	\begin{lemma}\label{lem:alphabets}
		Let $(\alpha, \beta) \in \overline{P}$. Then, there exists $W \in \langle U, V\rangle$ such that $\alpha = W(a)$ and $\beta = W(b)$. In particular, the sets $\mathsf{c}(\overline{P})$ and $P$ are equal.
	\end{lemma}
	
	\begin{proof}
		We will prove a stronger equality: $\mathsf{c}(\overline{P}_n) = P_n$ for each $n \geq 0$. It is enough to show one inclusion as both sets have cardinality $2^n$.
		
		We proceed by induction. We claim that, for every $n \geq 0$ and $(u, v) \in \overline{P}_n$, there exists $W \in \langle U, V\rangle$ such that $u = W(a)$ and $v = W(b)$. The base case, for $n = 0$, is clear.
		
		Now, let $(u, v) \in \overline{P}_{n-1}$ for $n \geq 1$. We will prove the claim for $(uv, v) \in \overline{P}_{n}$. Indeed, we have that there exists $W \in \langle U, V\rangle$ such that $u = W(a)$ and $v = W(b)$. Observe that $WU(a) = W(ab) = W(a)W(b) = uv$ and $WU(b) = W(b) = v$. The proof for $(u, uv) \in \overline{P}_n$ is analogous.
	\end{proof}
	
	To state the following lemmas, we need to fix some useful notation. Let $\alpha$ and $\beta$ be finite words and assume that $\alpha$ starts with $a$, and that $\beta$ ends with $b$. We write $\alpha = a\alpha^+$ and $\beta = \beta^- b$. Then, we define $\bfirst{\alpha} = b\alpha^+$ and $\alast{\beta} = \beta^-a$. That is, $\bfirst\alpha$ is obtained by replacing the first letter of $\alpha$ (which is $a$ by assumption) with $b$, and, similarly, $\alast\beta$ is obtained by replacing the last letter of $\beta$ (which is $b$ by assumption) with $a$.
	
	\begin{lemma}\label{lem:structure}
		For every $(\alpha,\beta)\in \overline P$, $\alpha$ starts with $a$, $\beta$ ends with $b$. Moreover, every word $\alpha^k\beta$, with $k\ge 1$, starts with $\alast\beta$, and every word $\alpha\beta^k$, with $k\ge 1$, ends with $\bfirst\alpha$. In particular, every sufficiently large word in $(\alpha,\beta)$ starts with $\beta^-$ and ends with $\alpha^+$, and we always have the equality $\alpha\beta=(\alast\beta)(\bfirst\alpha)$.
	\end{lemma}
	
	\begin{proof}
		For $(\alpha,\beta)=(a,b)$, we clearly have that $\alpha$ starts with $a$, $\beta$ ends with $b$, $\alpha^+=\beta^-=\emptyset$, $\bfirst\alpha=b$, $\alast\beta=a$, $\alpha^k\beta=a^k b$ starts with $a=\alast\beta$ for every $k\ge 1$, and $\alpha\beta^k=ab^k$ ends with $b=\bfirst\alpha$ for every $k\ge 1$.
		
		By induction, if $(A,B)=(\alpha,\alpha\beta)$ then $A = \alpha$ starts with $a$, and $B = \alpha\beta$ ends with $b$. Since $B=\alpha\beta$ ends with $\bfirst\alpha=\bfirst A$, then, for every $k\ge 1$, $AB^k$ also ends with $\bfirst A$. Now, fix $k \geq 1$. By induction, have that $\alpha^k\beta$ starts with $\alast\beta$, so $A^kB=\alpha^{k+1}\beta=\alpha\alpha^k\beta$ starts with $\alpha \alast\beta = \alast{(\alpha\beta)}=\alast B$.
		
		On the other hand, if $(A,B)=(\alpha\beta,\beta)$, then clearly $A$ starts with $a$, and $B$ ends with $b$. Since $A=\alpha\beta$ starts with $\alast\beta=\alast B$, then, for every $k\ge 1$, $A^kB$ starts with $\alast B$. Furthermore, since $\alpha\beta^k$ ends with $\bfirst\alpha$, $AB^k=\alpha\beta^{k+1}=\alpha\beta^k\beta$ ends with $\bfirst\alpha \beta=\bfirst{(\alpha \beta)}=\bfirst A$ for every $k\ge 1$. The inductive argument is therefore complete.
		
		Finally, the remaining equality $\alpha\beta=(\alast\beta)(\bfirst\alpha)$ follows immediately since $|\alpha\beta|=|(\alast\beta)(\bfirst\alpha)|$ (and, as we have just proved, $\alpha\beta$ starts with $\alast\beta$ and ends with $\bfirst\alpha$).
	\end{proof}
	
	\begin{remark} \label{rem:palindromic}
		Every word in $P$ is of the form $a\theta b$, with $\theta$ palindromic, i.e., $\theta$ coincides with its transpose $\theta^*$, as stated in Bombieri's article \cite[Proof of Theorem 15]{B:markoff_numbers}. Since $\alpha$ starts with $a$ and $\beta$ ends with $b$, this is equivalent to $(\alpha\beta)^* = \alast{(\bfirst{(\alpha\beta)})}$. In other words, both $\bfirst\alpha$ and $\alast\beta$ are palindromic for every pair $(\alpha, \beta) \in \overline{P}$. We will now present an alternative proof of this fact.
		
		As in the previous lemma, we will proceed by induction; the base case is clear. Suppose that $\bfirst\alpha$ and $\alast\beta$ are palindromic. Then $\bfirst{(\alpha\beta)}$ is palindromic, since both the word $\bfirst{(\alpha\beta)}=\bfirst\alpha \beta$ and the word $(\bfirst{(\alpha\beta)})^*=(\bfirst\alpha \beta)^*=\beta^* \bfirst\alpha$ are obtained from $\alpha\beta=(\alast\beta)(\bfirst\alpha)$ by replacing the first letter (which is $a$) with $b$, and therefore coincide. Similarly, $\alast{(\alpha\beta)}$ is also palindromic. Thus, the result holds for both $(\alpha\beta, \beta)$ and $(\alpha, \alpha\beta)$, which completes the inductive proof.
	\end{remark}
	
	\begin{lemma}\label{lem:determine_alphabet}
		Suppose that a word $w$ can be written as a concatenation $\tau\alpha\beta\tau'$ for some words $\tau$, $\tau'$, $\alpha$ and $\beta$, with $(\alpha,\beta)\in \overline P_n$ and $n \in \NN$. If there exist $(A,B)\in \overline P_n$, $k\ge 1$ and $w_1,\dotsc,w_k\in \{A,B\}$ such that $w=w_1\ldots w_k$, then $(A,B)=(\alpha,\beta)$ and there exists $1\le j<k$ such that $w_1\ldots w_{j-1}=\tau$, $w_j=\alpha$, $w_{j+1}=\beta$ and $w_{j+2}\ldots w_k=\tau'$.  
	\end{lemma}
	
	\begin{proof}
		As usual, we proceed by induction. The result is trivial for the base case $(\alpha,\beta)=(a,b)\in \overline P_0$. Assume now that $(\alpha,\beta)=(uv,v)$ for $(u,v)\in\overline P_{n-1}$, where $n \geq 1$. Let $w$ be a word such that $w = \tau \alpha \beta \tau'$ for some words $\tau, \tau'$ and assume that there exist $(A, B) \in \overline{P}_n$, $k \geq 1$ and $w_1, \dotsc, w_k \in \{A, B\}$ such that $w = w_1 \ldots w_k$.
		
		Since $w = \tau \alpha \beta \tau'$ and $\alpha \beta = uvv$, there exist $\sigma = \tau$ and $\sigma' = v\tau'$ such that $w = \sigma u v \sigma'$. Thus, by induction, if $w$ can be written as concatenation of words from a pair in $\overline{P}_{n-1}$, then the pair is necessarily $(u,v)$ and the words $u$, $v$ and $v$ appear consecutively in this decomposition. This is indeed the case as each $w_j$ for $1 \leq j \leq k$ is a concatenation of words from a pair in $\overline{P}_{n-1}$, so $w$ can be written in this way as well.
		
		We conclude that $(A, B) = \overline{U}(u,v) = (uv, v)$ or $(A, B) = \overline{V}(u,v) = (u, uv)$. Indeed, if this did not hold, then we would be able to find a different pair in $\overline{P}_{n-1}$ whose words can be concatenated to obtain $w$. Finally, if $(A, B) = \overline{V}(u, uv)$, then it would not be possible for the words $u$, $v$, $v$ to appear consecutively. We conclude that $(A, B) = (\alpha, \beta)$.
		
		The case where $(\alpha,\beta)=(u,uv)$ for $(u, v) \in \overline{P}_{n-1}$ is analogous.  
	\end{proof}
	
	We can now relate the length of a factor of a word in $P$ with the length of the smallest word in $P$ containing it:
	
	\begin{lemma} \label{lem:fundamental_period_bounded}
		Let $w$ be a factor of a word in $P$. Then, the length of the shortest word in $P$ containing $w$ is strictly smaller than $3|w|$.
	\end{lemma}
	
	\begin{proof}
		Let $(\alpha, \beta) \in \overline{P}$ such that $\alpha\beta$ contains $w$ and such that $|\alpha\beta|$ is minimal for this property. We will assume $|\alpha|>|\beta|$ (the case $|\alpha|<|\beta|$ is analogous, and the case $|\alpha|=|\beta|$ only occurs in the trivial case $\alpha=a, \beta=b$, in which we may replace the constant $3$ with $2$). Hence, we may write $\alpha=\tilde\alpha\beta^r$ for some $r\ge 1$, where $(\tilde\alpha, \beta) \in \overline P$ and $|\tilde\alpha|\le|\beta|$. We then have the bounds $(r+1)|\beta| < |\alpha\beta| \leq (r+2)|\beta|$.
		
		Observe that $w$ must intersect both $\alpha$ and $\beta$ by minimality of $|\alpha\beta|$. Indeed, if $w$ only intersects $\alpha = \tilde\alpha\beta^r$ or $\beta$, then the shorter word $\tilde\alpha\beta^{r-1}\beta$ corresponding to the pair $(\tilde\alpha\beta^{r-1}, \beta) \in \overline{P}$ contradicts the minimality of $|\alpha\beta|$. Now, if $w$ intersects the prefix $\tilde\alpha$ of $\alpha$, then it contains $\beta^r$ strictly, and so the ratio $|w|/|\alpha\beta|$ is larger than $r/(r+2)\ge 1/3$. Thus, from now on we may assume that $w$ is contained in $\beta^{r+1}$. Moreover, $r$ is minimal for this property as, otherwise, the pair $(\tilde\alpha\beta^{r-1}, \beta) \in \overline{P}$ again contradicts the minimality of $|\alpha\beta|$. Thus, $w = u \beta^{r-1} v$, where $u$ is a nonempty suffix of $\beta$ and $v$ is a nonempty prefix of $\beta$.
		
		Assume that $r \geq 2$. By \Cref{lem:structure}, we have that $\tilde{\alpha}\beta = (\alast{\beta})(\bfirst{\tilde\alpha})$. We now claim that $|u| \geq |\bfirst{\tilde\alpha}|$. Indeed, assume that this is not the case. Then, $w$ is contained in $\alpha = \tilde\alpha\beta^r$, since any proper suffix of $\bfirst{\tilde\alpha}$ is also a proper suffix of $\tilde\alpha$. This contradicts the minimality of $|\alpha\beta|$ as before. Thus, since $|\bfirst{\tilde\alpha}|=|\tilde\alpha|$, the ratio $|w|/|\alpha\beta|$ is at least $((r-1)|\beta|+|\tilde\alpha|)/((r+1)|\beta|+|\tilde\alpha|)$, which is larger than $(r-1)/(r+1) \geq 1/3$ since $r \geq 2$.
		
		We will now address the remaining case where $r=1$. First observe that if $|\tilde\alpha| = |\beta|$, then $\tilde\alpha = a$ and $\beta = b$. Hence, $\alpha = ab$ and $w = b^2$. We then have that $|w|/|\alpha\beta| = 2/3$. We can therefore assume from now on that $|\beta| > |\tilde\alpha|$ and we may write $\beta = \tilde\alpha^j\tilde\beta$ for some $j \geq 1$, where $(\tilde\alpha, \tilde\beta) \in \overline P$ and $|\tilde\beta|\le|\tilde\alpha|$.
		
		We have that $w$ is a factor of $\beta^2=\tilde\alpha^j\tilde\beta\tilde\alpha^j\tilde\beta$ and that it intersects both copies of $\beta$. Hence, $w = u v$, where $u$ is a nonempty suffix of $\beta = \tilde\alpha^j \tilde\beta$ and $v$ is a nonempty prefix of $\beta = \tilde\alpha^j \tilde\beta = \tilde\alpha\tilde\alpha^{j-1}\tilde\beta$.
		
		By \Cref{lem:structure}, $\tilde\alpha^j \tilde\beta$ ends with $\bfirst{\tilde\alpha}$. We claim that $\bfirst{\tilde\alpha}$ is a suffix of $u$. Indeed, if this were not the case, then $u$ would be a suffix of $\tilde\alpha$ and, hence, $w$ would be contained in the shorter word $\tilde\alpha\tilde\alpha^j \tilde\beta$ corresponding to the pair $(\tilde\alpha, \tilde\alpha^j\tilde\beta) \in \overline P$, which is not possible by the minimality of $|\alpha\beta|$. Similarly, \Cref{lem:structure} implies that $\tilde\alpha^j\tilde{\beta} = \tilde\alpha\tilde\alpha^{j-1}\tilde\beta$ starts with $\tilde\alpha^{j-1}\alast{\tilde{\beta}}$. We claim that $\tilde{\alpha}^{j-1}\alast{\tilde\beta}$ is a prefix of $v$. Indeed, if this were not the case, then $v$ would be a prefix of $\tilde\alpha^{j-1}\tilde\beta$ and, hence, $w$ would be contained in the shorter word $\tilde\alpha^j\tilde\beta\tilde\alpha^{j-1}\tilde\beta$ corresponding to the pair $(\tilde\alpha^j\tilde\beta, \tilde\alpha^{j-1}\tilde\beta) \in \overline P$, which is not possible by the minimality of $|\alpha\beta|$. Finally, we conclude from $|\bfirst{\tilde\alpha}|=|\tilde\alpha|$ and $|\alast{\tilde\beta}|=|\tilde\beta|$, that the ratio $|w|/|\alpha\beta|$ is at least $(j|\tilde\alpha|+|\tilde\beta|)/((2j+1)|\tilde\alpha|+2|\tilde\beta|)$, which is larger than $j/(2j+1)\ge 1/3$. 	
	\end{proof}

	\begin{remark}
		The general bound in the previous lemma cannot be improved. Indeed, for the word $w = bab^{k+1}a$ for $k \geq 1$, we have that $|\alpha\beta|$ is minimal for the pair $(\alpha, \beta) = (ab^kab^{k+1},ab^{k+1}) \in \overline{P}$. Since $|w| = 2(k+4)$ and $|\alpha\beta| = 2(3k+5)$, the ratio $|w|/|\alpha\beta|$ is arbitrarily close to $1/3$ when $k$ is sufficiently large.
		
		The previous example corresponds to the first case of the proof of the previous lemma, namely when $w$ intersects the prefix $\tilde\alpha$ of $\alpha = \tilde\alpha\beta^r$. In the two remaining cases of the proof, nevertheless, the bound \emph{can} be improved as we do below.
		
		Assume then that $w$ does not intersect $\tilde\alpha$. As in the previous proof, we first consider the case where $r \geq 2$. Then, we may replace the constant $3$ with $2 + \varepsilon$ for any $\varepsilon > 0$. Indeed, observe first that if $|\tilde\alpha| = |\beta|$, then $\tilde\alpha = a$ and $\beta = b$, so $\alpha = ab^r$ and $w = b^{r+1}$. Thus, the ratio $|w|/|\alpha\beta|$ is $(r+1)/(r+2) \geq 3/4 \geq 1/2$. Otherwise, if $|\beta| > |\tilde\alpha|$, we write $\beta = \tilde\alpha^j\tilde\beta$ for $j \geq 1$ and $(\tilde\alpha, \tilde\beta) \in \overline{P}$. We have that $w = uv$, where $u$ is a suffix of $\beta^r = \tilde\alpha^j\tilde\beta\beta^{r-1}$ and $v$ is prefix of $\beta = \tilde\alpha^j\tilde\beta$. By \Cref{lem:structure}, $\tilde\alpha^j\tilde\beta$ ends with $\bfirst{\tilde\alpha}$ and we claim that $\bfirst{\tilde\alpha}\beta^{r-1}$ is a suffix of $u$. Indeed, if this were not the case, then $u$ would be a suffix of $\tilde\alpha\beta^{r-1}$, so $w$ would be contained in the shorter word $\alpha = \tilde\alpha\beta^r$ corresponding to the pair $(\tilde\alpha\beta^{r-1},\beta) \in \overline{P}$.
		Similarly, if we put $\hat\beta = \tilde\alpha^{j-1}\tilde\beta$ we have that $(\tilde\alpha, \hat\beta) = (\tilde\alpha, \tilde\alpha^{j-1}\tilde\beta) \in \overline{P}$, so \Cref{lem:structure} implies that $\tilde\alpha\hat{\beta}$ starts with $\alast{\hat{\beta}}$. We claim that $\alast{\hat\beta}$ is a prefix of $v$. Indeed, if we assume otherwise, then $v$ is a prefix of $\hat\beta$ and, thus, $w$ is contained in the shorter word $(\tilde\alpha\hat\beta)^r\hat\beta$ corresponding to the pair $(\tilde\alpha\hat\beta, (\tilde\alpha\hat\beta)^{r-1}\hat\beta) \in \overline{P}$, a contradiction. Therefore, the ratio $|w|/|\alpha\beta|$ is at least $((r-1)|\beta|+|\tilde\alpha|+|\hat\beta|)/((r+1)|\beta|+|\tilde\alpha|)=r|\beta|/((r+1)|\beta|+|\tilde\alpha|)$, which is larger than $r/(r+2)\ge 1/2$. 
		
		Finally, we analyze the case where $r = 1$ and show that we can replace the constant $3$ with $5/2 + \varepsilon$ for any $\varepsilon > 0$. Recall that $\beta = \tilde\alpha^j\tilde\beta$, so the result is clear when $j \geq 2$ as $j/(2j+1) \geq 2/5$. Thus, we will assume that $j = 1$, so $\alpha = \tilde\alpha\tilde\alpha\tilde\beta$ and that $\beta = \tilde\alpha\tilde\beta$. If $j = 1$ and $|\tilde\alpha| = |\tilde\beta|$, then $\tilde\alpha = a$, $\tilde\beta = b$, $\alpha = aab$ and $\beta = ab$. Since $w$ intersects both $\alpha$ and $\beta$, we have that $|w| \geq 4$, so we obtain $|w|/|\alpha\beta| \geq 2/5$ once again. We will then assume that $|\tilde\beta| < |\tilde\alpha|$. We have that $w$ is a factor of $\tilde\alpha\tilde\beta\tilde\alpha\tilde\beta$, which is in turn a factor of the shorter word $\tilde\alpha\tilde\beta\tilde\alpha\tilde\beta\tilde\beta$ corresponding to the pair $(\tilde\alpha\tilde\beta, \tilde\alpha\tilde\beta\tilde\beta) \in \overline{P}$, a contradiction.
	\end{remark}
	
	The previous lemma allows us to control the size of the set $\Sigma(3, n)$:
	
	\begin{corollary} \label{lem:O(n^3)}
		For all $n \geq 1$, we have $|\Sigma(3, n)| \le 9n^3$.
	\end{corollary}
	\begin{proof}
		If $w\in \Sigma(3,n)$, then there is $(\alpha,\beta)\in \overline{P}$ with $|\alpha\beta|<3n$ such that $w$ is a factor of $\alpha\beta$ by \Cref{lem:bombieri_old,lem:fundamental_period_bounded}. Notice now that the pair $(\alpha,\beta)\in \overline{P}$ is determined by the irreducible fraction $|\alpha|/|\beta|$: indeed, $|\alpha|=|\beta|$ if and only if $\alpha=a$ and $\beta=b$; if $|\alpha|>|\beta|$, then $\alpha=\tilde\alpha\beta^k$ for some positive integer $k$ with $(\tilde\alpha,\beta)\in \overline{P}$ and $|\tilde\alpha|\le|\beta|$, and, thus, $|\alpha|/|\beta|=k+|\tilde\alpha|/|\beta|$; and, if $|\alpha|<|\beta|$, then $\beta=\alpha^k\tilde\beta$ and $|\alpha|/|\beta|=1/(k+|\tilde\beta|/|\alpha|)$. Hence, our claim follows by induction on the number of elements of the continued fraction of $|\alpha|/|\beta|$. 
		
		The number of such fractions $|\alpha|/|\beta|$ is bounded by the number of pairs $(i,j)$ of positive numbers with $i+j\le 3n$, which is $3n(3n-1)/2<9n^2/2$. Since a word of size smaller than $3n$ has at most $2n$ factors of size $n$, there are at most $2n \cdot 9n^2/2=9n^3$ elements in $\Sigma(3,n)$.
	\end{proof}
	
	Recall that $U$ and $V$ are the Nielsen operators given by $U(a) = ab$, $U(b) = b$, $V(a) = a$ and $V(b) = ab$.
	
	\begin{lemma} \label{lem:reverse}
		For any finite word $w$ in the alphabet $\{a, b\}$, we have the identities $bU(w^*) = U(w)^*b$ and $V(w^*)a = aV(w)^*$. In particular, if $w$ is a palindrome, then $b U(w)$ and $V(w) a$ are palindromes as well.
	\end{lemma}
	
	\begin{proof}
		This was already done by Bombieri \cite[Proof of Theorem 15]{B:markoff_numbers}, but for the sake of completeness we include a short proof by induction.
		
		These identities are trivial if $|w| = 0$. We then assume that they hold for words of length $n-1$ for $n \geq 1$; let $w$ be a word of such length. If $\tilde{w} = aw$, then
		\begin{align*}
			b U(\tilde{w}^*) &= bU(w^* a) = bU(w^*)ab = U(w)^*bab = U(\tilde{w})^*b \\
			V(\tilde{w}^*) a &= V(w^* a)a = V(w^*)aa = a V(w)^* a = a V(\tilde{w})^*.
		\end{align*}
		On the other hand, if $\tilde{w} = bw$, then
		\begin{align*}
			b U(\tilde{w}^*) &= b U(w^* b) = bU(w^*)b = U(w)^*bb = U(\tilde{w})^*b \\
			V(\tilde{w}^*) a &= V(w^* b) a = V(w^*)aba = aV(w)^*ba = aV(\tilde{w})^*.
		\end{align*}
		
		Assume now that $w$ is a palindrome. Then,
		\begin{align*}
			(b U(w))^* &= U(w)^* b = b U(w^*) = b U(w) \\
			(V(w) a)^* &= a V(w)^* = V(w^*) a = V(w) a.
		\end{align*}
	\end{proof}
	
	The following lemma shows that bi-infinite words with Markov value exponentially close to $3$ (relative to the size of the interval they induce) cannot contain both $\alpha\alpha$ and $\beta\beta$ if $(\alpha, \beta) \in \overline{P}$. Recall that
	\[
	\sizer(w) = \lfloor\log (\sizes(\alpha)^{-1}) \rfloor = \lfloor\log (|I(\alpha)|^{-1}) \rfloor.
	\]
	
	\begin{lemma}\label{lem:aa_bbr}
		Let $(\alpha,\beta)\in \overline P$. If $w$ is a finite word in the alphabet $\{\alpha,\beta\}$ starting with $\alpha\alpha$ and ending by $\beta\beta$ such that $\sizer(w)\leq r$, then the Markov value of any bi-infinite word containing $w$ as a factor is larger than $3+e^{-r}$. Moreover, if $w$ contains $\alpha\alpha\beta\beta$ as a factor and $\sizer(w)\leq 2r$, we have that the Markov value of any bi-infinite word containing $w$ as a factor is larger than $3+e^{-r}$.
	\end{lemma}
	
	\begin{proof}
		We will present this proof in several steps. The first step, that we label as ``Step 0'', is not strictly necessary; it is contained in the other more general steps. In this step we make an estimate depending on $|w|$, and it is weaker than the estimate of the statement, which depends on $\sizer(w)$. However, we included it since it contributes to the understanding of the overall strategy. 
		
		\medbreak
		
		\noindent \textbf{Step 0:} Assume that $\alpha = a$ and $\beta = b$. Without loss of generality, we can assume that $w = aa(ba)^k bb$, since, otherwise there is a factor of $w$ of this form and we may replace $w$ by this factor. We will consider two cuts of this word. One cut, to which we will refer as the ``first cut'' is $aa(ba)^k|bb$, while the ''second cut'' is $aa(ba)^kbb|$. We start by applying \Cref{lem:R_S_bound} to the first cut. This immediately shows that $k \geq 1$, as, otherwise, any bi-infinite word containing $w$ has a Markov value of at least $3 + \frac{1}{7}$ (in the general case this is not immediate; it is treated in Step 2). Hence, we assume that $k \geq 1$.
		
		Let $\omega$ be a bi-infinite word containing $w$ and assume by contradiction that its Markov value is smaller than $3 + \frac{1}{7}(3 - 2\sqrt{2})^{|w|}$. We continue drawing conclusions from \Cref{lem:R_S_bound}: the first cut shows that $\omega$ must contain an $a$ to the right of $w$. Thus, $\omega$ contains $w' = aa(ba)^k|bb|a$, where we again marked both cuts. We now use these cuts to conclude inductively that $w'$ must be followed with $(ba)^{k-1}$ in $\omega$: each $b$ is forced by the second cut (since there is a $b$ at the symmetric position with respect to the second cut), and it is followed with an $a$ by the first cut (since there is an $a$ at the symmetric position with respect to the first cut).
		
		Set $\gamma = (ba)^{k-1}$. Between both cuts, we have the word $bb$ which we will write as $b\theta b$ with $\theta = \emptyset$ (in the general case, $\theta$ can be more complicated). At the left of the first cut, we have a word of the form $(\theta b a \gamma a)^*a$, while the second cut is followed with $a\gamma$. Thus, $\omega$ contains the word $w'' = (\theta b a \gamma a)^*a|b\theta b|a\gamma$, where we again marked the first and second cuts.
		
		The structure above is precisely the configuration that we will try to replicate the general case, as it already leads to a large Markov value. Indeed, using the first cut again, we obtain that $w''$ is followed with an $a$ in $\omega$. Finally, $w''$ can also be written as $w'' = (\theta^* b a \theta b a \gamma a)^* b|a\gamma$ (where only the second cut is marked). Since $\theta^*ba\theta ba\gamma$ starts with $\gamma b$, we obtain that $w''$ is followed with a $b$ inside $\omega$, which contradicts that it is followed with an $a$ as we obtained before. In other words, we have shown that any bi-infinite word containing $w$ has a Markov value of at least $3 + \frac{1}{7}(3 - 2\sqrt{2})^{|\gamma|+1} > 3 + \frac{1}{7}(3 - 2\sqrt{2})^{|w|}$.
		
		\medbreak
		
		\noindent \textbf{Step 1:} We now start treating the general case, so assume that $w$ starts with $\alpha\alpha$ and ends with $\beta\beta$. Since $(\alpha, \beta) \in \overline P$, \Cref{lem:alphabets} shows that there exists some $W \in \langle U, V\rangle$ such that $\alpha = W(a)$ and $\beta = W(b)$. Thus, $w$ is the image by $W$ of a word in the alphabet $\{a, b\}$ starting with $aa$ and ending with $bb$. Without loss of generality, we assume that $w = \alpha\alpha (\beta\alpha)^k \beta\beta$ with $k \geq 0$, as, otherwise, $w$ contains a factor of this form and we may replace $w$ with this factor.
		
		\medbreak
		
		\noindent \textbf{Step 2:} In this step, we assume that $k = 0$, so $w = \alpha\alpha\beta\beta$. We claim that $w$ contains a cut of the form $\tau a|b\theta b$, where $\tau^*$ starts with $\theta a$ and $\theta$ is a palindromic word. This leads to a contradiction by \Cref{lem:calc_s}, as then the Markov value of any bi-infinite word containing $w$ is at least $3 + \sizes(b\theta b)$ which is larger than $3+e^{-r}$ by the following computation. 
		
		By hypothesis, we have that $\sizes(w) \geq e^{-2r-1}$, so we obtain that $\sizes(\alpha\alpha\beta\beta) \geq e^{-2r-1}$. Write $\theta=\theta_1 \ldots \theta_n$. By \eqref{eq:sisesgap} we have the inequalities $\sizes(a\theta a)^{-1}\geq 961q_n(\theta)^2$, and $\sizes(b\theta b)^{-1}\leq 162q_n(\theta)^2$. Therefore,
		
		\[
		e^{-2r-1} \leq \sizes(\alpha\alpha\beta\beta)\leq \sizes(a \theta a b \theta b) \leq 2\sizes(a\theta a)\sizes(b\theta b)\leq  \frac{324}{961}\sizes(b\theta b)^2,
		\]
		hence $\sizes(b\theta b)\geq e^{-r}$.

		We proceed by induction: in the base case, we have $\theta = \emptyset$ and $\tau = a$. Now, observe that
		\[
		U(\tau ab\theta b) = U(\tau) a|b b U(\theta) b = \tilde\tau a|b \tilde\theta b,
		\]
		where $\tilde\tau = U(\tau)$ and $\tilde\theta = b U(\theta)$, and we have adjusted the position of the cut. We claim that $\tilde{\tau}^* = U(\tau)^*$ starts with $b U(\theta) a = \tilde{\theta} a$. Indeed, since $\tau^*$ starts with $\theta a$, we have that $\tau$ ends with $a \theta^*$. Thus, $U(\tau)$ ends with $ab U(\theta^*)$. Therefore, $U(\tau)^*$ starts with $U(\theta^*)^*b a$, which is equal, by \Cref{lem:reverse}, to $b U(\theta) a = \tilde\theta a$ (since $\theta$ is a palindrome).
		
		On the other hand, observe that
		\[
		V(\tau ab\theta b) = V(\tau) a a|b V(\theta) ab = \tilde\tau a|b \tilde\theta b,
		\]
		where $\tilde\tau = V(\tau)a$ and $\tilde\theta = V(\theta)a$, and we have adjusted the position of the cut. We claim that $\tilde{\tau}^* = a V(\tau)^*$ starts with $V(\theta)aa = \tilde\theta a$. Indeed, first observe that, by \Cref{lem:reverse}, $\tilde{\tau}^* = V(\tau^*)a$. Now, we consider two cases. If $\tau^* = \theta a$, then $\tilde{\tau}^* = V(\theta a)a = V(\theta)aa = \tilde\theta a$. Otherwise, $\tau^*$ starts with $\theta a c$ where $c \in \{a, b\}$, so $\tilde{\tau}^*$ starts with $V(\theta a c) = V(\theta) a V(c)$. Since $V(c)$ starts with $a$ whether $c = a$ or $c = b$, we obtain that $\tilde{\tau}^*$ starts with $V(\theta) a a = \tilde{\theta} a$.
		
		Since, by Step 1, there exists $W \in \langle U, V \rangle$ such that $W(a) = \alpha$ and $W(b) = \beta$, this concludes the proof when $k = 0$.
		
		\medbreak
		
		\noindent \textbf{Step 3:} In this step we leverage the structure found in Step 0 when $k \geq 1$ and shows that it also leads to a large Markov value in a more general context. Assume now that we have a word $w$ with two cuts of the form $w = \tau a|b\theta b|$ such that:
		\begin{enumerate}
			\item there exists a word $\gamma$ such that $\tau$ ends with $(\theta b a \gamma a)^*$; and \label{it:1}
			\item $\theta^* b a \theta b a \gamma$ starts with $\gamma b$. \label{it:2}
		\end{enumerate}
		We have shown that \eqref{it:1} and \eqref{it:2} hold for the base case $w = aa(ba)^kbb$ with $\tau = aa(ba)^{k-1}b$ and $\theta = \emptyset$.
		
		Then, as before, the Markov value of any bi-infinite word $\omega$ containing $w$ is at least $3 + e^{-r}$. To see this, we will again use that, by \Cref{lem:calc_s}, some of the letters surrounding $w$ are forced in $\omega$ for the Markov value to remain below this value; eventually this will not be possible anymore. Indeed, an $a$ is forced after $w$ by the first cut, since $\tau$ ends with $(\theta b a)^*$. Moreover, \Cref{lem:1} shows that the configuration $\tau a|b \theta b|a$ is followed by $\gamma$: each $a$ of $\gamma$ is forced by the first cut (since $\tau$ ends with $(\theta b a \gamma)^*$), while each $b$ of $\gamma$ is forced by the second cut (since $\theta^* b a \theta b a \gamma$ starts with $\gamma b$). Finally, the first cut forces an $a$ after $\tau a|b \theta b|a\gamma$ (since $\tau$ ends with $(\theta b a \gamma a)^*$), while, on the contrary, the second cut forces a $b$ after $\tau a|b \theta b|a\gamma$ (since $\theta^* b a \theta b a \gamma$ starts with $\gamma b$). Thus, we obtain that the Markov value of any bi-infinite word containing $w$ is at least $3 + e^{-r}$.
		
		To be more precise about this last part, observe that $\gamma$ cannot be followed by $12$ or $21$, because otherwise we will find a sequence of the form $c' c^sc'$ where $c,c'\in\{1,2\}$ with $c\neq c'$ and $s$ odd, but using \Cref{prop:rtranspose}, the fact that $\sizes$ is monotone and the fact that $w$ ends with $(\theta b a \gamma a)^*(ab\theta b)$
		\[
		\sizes(c^s)\geq\sizes(a\gamma a)\geq2^{-1}\sizes(a\gamma^*ab\theta^*)\geq2^{-2}\sizes(abb)^{-1}\sizes(w)= (775/2)\sizes(w)
		\]
		
		whence $\sizer(c^s)\leq r-4$, a contradiction with \Cref{lem:1}. If $\gamma$ is followed by $b$, then writing the first cut as $\omega^\ast=R^\ast b\eta^\ast b|a\eta aS$ with $\eta=\theta ba\gamma $ we have, by \Cref{lem:calc_s}, that
		\[
		\lambda(\omega)=\lambda(\omega^\ast)\geq 3+\sizes(b\eta b).
		\]
		Since $a\eta^\ast ab\theta b$ is a subword of $w$, we have that $\sizer(a\eta^\ast ab\theta b)\leq r$ by \Cref{lem:bounds_subwords}. In particular, $\sizes(a\eta^\ast ab\theta b)\geq e^{-r-1}$. On the other hand, by \Cref{prop:rtranspose} one has that $\sizes(a\eta^\ast ab\theta b)\leq 4\sizes(a\eta a)\sizes(b\theta b)\leq \sizes(b\eta b)/3$, whence $\sizes(b\eta b)\geq e^{-r}$.

		Similarly, if the word $\gamma$ is followed by $a$, then, by writing the second cut as $\omega=R^\ast b\gamma^\ast b|a\gamma aS$, we have
		\[
		\lambda(\omega)\geq 3+\sizes(b\gamma b). 
		\]
		Finally, since $\gamma$ is a subword of $\eta=\theta ba\gamma$, by \Cref{lem:bounds_subwords} again we get that $\sizes(b\gamma b)\geq \sizes(b\eta b)\geq e^{-r}$.
		
		\medbreak
		
		\noindent \textbf{Step 4:} We now show inductively that the previous structure (namely properties \eqref{it:1} and \eqref{it:2}) persists when we apply $U$ or $V$ to $w = \tau a|b \theta b|$. First, observe that, after adjusting the position of the cuts, we have that
		
		\begin{equation} \label{eq:UV}
			U(w) = U(\tau) a|b b U(\theta)b| \quad \text{ and } \quad V(w) = V(\tau) a a|b V(\theta) ab|.
		\end{equation}
		
		Thus, we have that $U(w)=\tilde\tau a|b\tilde\theta b|$, with $\tilde\tau=U(\tau)$ and $\tilde\theta=bU(\theta)$. Let $\tilde\gamma=bU(\gamma)$. Then, since $w$ satisfies \eqref{it:1}, $\tilde\tau=U(\tau)$ ends with
		\begin{align*}
			U((\theta ba\gamma a)^*)&=U(a\gamma^*ab\theta^*)=abU(\gamma^*)abbU(\theta^*)\\&=aU(\gamma)^*babU(\theta)^*b
			=(\tilde\theta ba\tilde\gamma a)^*
		\end{align*}
		where we used \Cref{lem:reverse}. This shows that \eqref{it:1} holds for $U(w)$. Similarly, this lemma shows that
		\begin{align*}
			\tilde\theta^*ba\tilde\theta ba\tilde\gamma&=U(\theta)^*bbabU(\theta)babU(\gamma)\\
			&=bU(\theta^*)babU(\theta)babU(\gamma)\\
			&=b U(\theta^* b a \theta b a \gamma).
		\end{align*}
		This word starts with $bU(\gamma b)=b U(\gamma) b=\tilde\gamma b$, since $\theta^*ba\theta ba\gamma$ starts with $\gamma b$, as $w$ satisfies \eqref{it:2}. Hence, we obtain that \eqref{it:2} also holds for $U(w)$.
		
		Now, from \eqref{eq:UV}, we have $V(w)=\tilde\tau a|b\tilde\theta b|$, with $\tilde\tau=V(\tau) a$ and $\tilde\theta=V(\theta)a$. Let $\tilde\gamma=V(\gamma)a$. Then, since $w$ satisfies \eqref{it:1}, $\tilde\tau=V(\tau) a$ ends with
		\begin{align*}
			V((\theta ba\gamma a)^*)a=aV(\theta ba\gamma a)^*=a(V(\theta)abaV(\gamma)a)^*=(\tilde\theta ba\tilde\gamma a)^*
		\end{align*}
		where we used \Cref{lem:reverse}. This shows that \eqref{it:1} holds for $V(w)$. Similarly, this lemma shows that
		\begin{align*}
			\tilde\theta^*ba\tilde\theta ba\tilde\gamma&=(V(\theta)a)^*baV(\theta)abaV(\gamma)a\\
			&=aV(\theta)^*baV(\theta)abaV(\gamma)a\\
			&=V(\theta^*)abaV(\theta)abaV(\gamma)a \\
			&=V(\theta^* b a \theta b a \gamma)a.
		\end{align*}
		This word starts with $V(\gamma b) = V(\gamma)ab=\tilde\gamma b$, since $\theta^*ba\theta ba\gamma$ starts with $\gamma b$, as \eqref{it:2} holds for $w$. Hence, we obtain that \eqref{it:2} also holds for $V(w)$.
		
		Since, by Step 1, there exists $W \in \langle U, V \rangle$ such that $W(a) = \alpha$ and $W(b) = \beta$, this concludes the proof when $k \geq 1$.
	\end{proof}
	
	In order to consider other possible cases, such as words starting with $\beta\beta$ and ending with $\alpha\alpha$, we will show some symmetry properties of the pairs in $\overline P$.
	\begin{lemma}\label{lem:symmetry}
		Let $(u,v) \in \overline P$. If $(\alpha, \beta) = (u, uv)$, then $\alpha^k \beta = (\bfirst{u} \alpha^k \alast{v})^*$. Similarly, if $(\alpha, \beta) = (uv, v)$, then $\alpha \beta^k = (\bfirst{u} \beta^k \alast{v})^*$.
	\end{lemma}
	
	\begin{proof}
		Assume first that $(\alpha, \beta) = (u, uv)$. We have that $(u, u^k v) \in \overline P$ for any $k \geq 1$. 
		Now, recall that, by \Cref{lem:structure}, $u u^k v = \alast{(u^k v)} \bfirst{u} = u^k \alast{v} \bfirst{u}$. Moreover, both $\bfirst{u}$ and $u^k \alast{v}$ are palindromic by \Cref{rem:palindromic}. Thus,
		\[
		\alpha^k \beta = u u^k v = u^k \alast{v} \bfirst{u} = (u^k \alast{v})^* (\bfirst{u})^* = (\bfirst{u} u^k \alast{v})^* = (\bfirst{u} \alpha^k \alast{v})^*.
		\]
		Similarly, if $(\alpha, \beta) = (uv, v)$, we have that $(uv^k, v) \in \overline P$ for any $k \geq 1$. Now, using \Cref{lem:structure} again, we obtain that $uv^k v = \alast{v} \bfirst{(uv^k)} = \alast{v} \bfirst{u} v^k$, where both $\alast{v}$ and $\bfirst{u} v^k$ are palindromic by \Cref{rem:palindromic}. Hence,
		\[
		\alpha \beta^k = uv^k v = \alast{v} \bfirst{u} v^k = (\alast{v})^* (\bfirst{u} v^k)^* = (\bfirst{u} v^k \alast{v})^* = (\bfirst{u} \beta^k \alast{v})^*.
		\]
	\end{proof}
	
	\begin{lemma}\label{lem:symmetry2}
		Let $(u, v) \in \overline{P}$ and let $e_1, \dotsc, e_k \geq 1$. If $(\alpha, \beta) = (u, uv)$, then
		\[
		\bfirst{u}\beta\alpha^{e_1} \beta \alpha^{e_2} \beta \ldots \alpha^{e_k} \alast{v}=(\alpha^{e_k}\beta\alpha^{e_{k-1}}\beta\ldots\beta\alpha^{e_1}\beta\beta)^*,
		\]
		while if $(\alpha, \beta) = (uv, v)$, then
		\[
		\bfirst{u}\beta^{e_1}\alpha\beta^{e_2}\alpha \ldots \beta^{e_k}\alpha \alast{v} = (\alpha\alpha \beta^{e_k} \alpha \beta^{e_{k-1}} \alpha \ldots \alpha \beta^{e_1})^*.
		\]
	\end{lemma}
	
	\begin{proof}
		Assume first that $(\alpha, \beta) = (u, uv)$. Then, by \Cref{lem:symmetry,lem:structure},
		\begin{align*}
			\bfirst{u}\beta\alpha^{e_1} \beta \alpha^{e_2} \beta \ldots \alpha^{e_k} \alast{v}&=\bfirst{u}(uv)\alpha^{e_1}(uv)\ldots (uv)\alpha^{e_k}\alast{v} \\
			&=(\bfirst{u})(\alast{v}\bfirst{u})\alpha^{e_1}(\alast{v}\bfirst{u})\ldots (\alast{v}\bfirst{u})\alpha^{e_k}\alast{v} \\
			&=(\bfirst{u}\alast{v})(\bfirst{u}\alpha^{e_1}\alast{v})\bfirst{u}\ldots \alast{v}(\bfirst{u}\alpha^{e_k}\alast{v}) \\
			&=(uv)^\ast(\alpha^{e_1}\beta)^\ast\ldots (\alpha^{e_k}\beta)^\ast \\
			&=(\alpha^{e_k}\beta\alpha^{e_k-1}\beta\ldots\beta\alpha^{e_1}\beta\beta)^*.
		\end{align*}
		Now, take $(\alpha, \beta) = (uv, v)$. Then, by \Cref{lem:symmetry,lem:structure},
		\begin{align*}
			\bfirst{u}\beta^{e_1}\alpha \ldots \beta^{e_k}\alpha \alast{v} &= (\bfirst{u}) \beta^{e_1} (uv) \ldots \beta^{e_k} (uv) \alast{v} \\
			&= (\bfirst{u}) \beta^{e_1} (\alast{v}\bfirst{u}) \ldots \beta^{e_k} (\alast{v}\bfirst{u}) \alast{v} \\
			&= (\bfirst{u} \beta^{e_1} \alast{v})\bfirst{u} \ldots \beta^{e_k} \alast{v}(\bfirst{u} \alast{v}) \\
			&= (\alpha\beta^{e_1})^* \ldots (\alpha \beta^{e_k})^* (uv)^* \\
			&= (\alpha \alpha \beta^{e_k}\alpha \beta^{e_k-1} \ldots \alpha\beta^{e_1})^*.
		\end{align*}
	\end{proof}
	
	The three previous lemmas imply that we obtain a large Markov value in the case where
	$(\alpha,\beta)=(u,uv)$ for any word of the form $\bfirst{u} \beta\ldots \alpha\alpha \alast{v}$, and in the case where $(\alpha,\beta)=(uv,v)$ for any word of the form $\bfirst{u} \beta\beta\ldots \alpha \alast{v}$.
	
	We now define the notion of a weakly renormalizable word, which is central to our methods as it is used to find suitable alphabets in which words can be written.
	
	\begin{definition}\label{def:weaklyrenormalizable}
		Let $(\alpha,\beta)\in \overline P$ and $w\in\langle a,b\rangle$ be a finite word. We say that $w$ is \emph{$(\alpha,\beta)$-weakly renormalizable} if we can write $w=w_1 \gamma w_2$ where $\gamma$ is a word (called the \emph{renormalization kernel}) in the alphabet $\{\alpha,\beta\}$ and $w_1, w_2$ are (possibly empty) finite words with $|w_1|, |w_2|<\max\{|\alpha|,|\beta|\}$ such that $w_2$ is a prefix of $\alpha\beta$ and $w_1$ is a suffix of $\alpha\beta$, with the following restrictions:
		
		If $(\alpha,\beta)=(u,uv)$ for some $(u,v)\in\overline P$ and $\gamma$ ends with $\alpha$, then $|v|\le |w_2|$. If $(\alpha,\beta)=(uv,v)$ for some $(u,v)\in\overline P$ and $\gamma$ starts with $\beta$, then $|u|\le |w_1|$.
	\end{definition}
	\begin{definition}\label{def:semirenormalizable}
		Let $(\alpha,\beta)\in \overline P$ and $w\in \langle 1,2\rangle$ be a finite word. We say that $w$ is \emph{$(\alpha,\beta)$-semi renormalizable} if there is an extension $\tilde{w}$ of at most two digits, one to the left and one to the right such that $\tilde{w}$ is \emph{$(\alpha,\beta)$-weakly renormalizable}.
	\end{definition}
	
	\medbreak
	
	The previous definition is motivated by the following ideas. Given an alphabet $\{\alpha, \beta\}$ with $(\alpha, \beta) \in \overline{P}$, it may not be possible to write a word $w$ in terms of $\alpha$ and $\beta$. Nevertheless, it may very well be possible to write ``most'' of $w$ in terms of $\alpha$ and $\beta$, preceded by and followed by some short trailing words. These words are $w_1$ and $w_2$ in the previous definition, and the condition ensuring that they are short is that $|w_1|, |w_2| < \max\{|\alpha|, |\beta|\}$. Indeed, if, for example, $|w_1| \geq \max\{|\alpha|, |\beta|\}$, then either $w_1$ ends with $\alpha$ or $\beta$ in $\{\alpha, \beta\}$ (so our choice of renormalization kernel was spurious; it should be longer), or it does not (so $w$ is actually not well described by the alphabet $\{\alpha, \beta\}$). To further ensure that $w_1$ and $w_2$ are well-adjusted to the chosen alphabet, we also require them to be a prefix or suffix of $\alpha\beta$; then $w$ is contained in $\alpha\beta \gamma \alpha\beta$, where the renormalization kernel $\gamma$ can be written in the alphabet $\{\alpha, \beta\}$.
	
	Finally, we need to ensure that the first and last letters of the renormalization kernel are chosen appropriately. This follows from the following lemma (which is essentially already contained in \Cref{def:weaklyrenormalizable}).
	
	\begin{lemma}
		Let $(\alpha, \beta) \in \overline{P}$ and $w \in \langle a, b \rangle$ be an $(\alpha, \beta)$-weakly renormalizable word. Write $w = w_1 \gamma w_2$ as in \Cref{def:weaklyrenormalizable}.
		
		If $(\alpha, \beta) = (u, uv)$ for some $(u, v) \in \overline{P}$ and $\gamma$ ends with $\alpha = u$, then $w_2$ starts with $\alast{v} \neq v$. Moreover, the word $\theta$ consisting of the last $|u|$ letters of $\gamma$ followed by the first $|v|$ letters of $w_2$ is different from $\beta$.
		
		Similarly, if $(\alpha, \beta) = (uv, v)$ for some $(u, v) \in \overline{P}$ and $\gamma$ starts with $\beta = v$, then $w_1$ ends with $\bfirst{u} \neq u$. Moreover, the word $\theta$ consisting of the last $|u|$ letters of $w_1$ followed by the first $|v|$ letters of $\gamma$ is different from $\alpha$.

	\end{lemma}
	
	\begin{proof}
		Assume first that $(\alpha, \beta) = (u, uv)$ and that $\gamma$ ends with $\alpha = u$. Then, \Cref{def:weaklyrenormalizable} ensures that $|v| \leq |w_2|$. Since $w_2$ is a prefix of $\alpha\beta$ of length at least $|v|$, \Cref{lem:symmetry} implies that $w_2$ starts with $\alast{v} \neq v$ (since $v$ ends with $b$, and $\alast{v}$ is palindromic by \Cref{rem:palindromic}). Now, $\theta$ ends with the first $|v|$ letters of $w_2$, so it ends with $\alast{v} \neq v$. Therefore, it cannot be equal to $\beta = uv$.
		
		Similarly, if $(\alpha, \beta) = (uv, v)$ and $\gamma$ starts with $\beta = v$, then \Cref{def:weaklyrenormalizable} ensures that $|u| \leq |w_1|$. Since $w_1$ is a suffix of $\alpha\beta$ of length at least $u$, \Cref{lem:symmetry} implies that $w_1$ ends with $\bfirst{u} \neq u$ (since $u$ starts with $a$, and $\bfirst{u}$ is palindromic by \Cref{rem:palindromic}). Now, $\theta$ starts with the last $u$ letters of $w_1$, so it starts with $\bfirst{u} \neq u$. Therefore, it cannot be equal to $\alpha = uv$.
	\end{proof}
	
	The previous lemma can be understood as follows. Since the renormalization kernel $\gamma$ is the part of $w = w_1 \gamma w_2$ that can be written in the alphabet $(\alpha, \beta)$, it should be as long as possible (in the sense that $w_1$ and $w_2$ are just ``short trailing words''). Hence, if $(\alpha, \beta) = (u, uv)$ and $\gamma$ ends with $\alpha = u$, then the word $w_2$ should not start with $v$ since, otherwise, $\gamma$ should instead end with $\beta = uv$ (and $w_2$ should be shorter). Similarly, if $(\alpha, \beta) = (uv, v)$ and $\gamma$ starts with $\beta = v$, then the word $w_1$ should not start with $u$ since, otherwise, $\gamma$ should instead start with $\alpha = uv$ (and $w_1$ should be shorter). All of these undesirable cases are ruled out by the previous lemma.
	
	Exhibiting a word as being $(\alpha, \beta)$-weakly renormalizable is nontrivial in general and, to complicate matters even further, the choice of alphabet $(\alpha, \beta) \in \overline P$ is not clear to begin with. Nevertheless, any word in the alphabet $\{a, b\}$ is trivially $(a, b)$-weakly renormalizable (by setting the renormalization kernel equal to the entire word).
	
	On the other hand, there are subwords of words in $\langle a, b\rangle$ that can fail to be weakly renormalizable (for any alphabet) with nontrivial kernel, because they are missing one digit at one (or both) of their ends. For example, the word of even length $w=21\ldots 1$ is a subword of $b^\infty ab^\infty$, and hence it belongs to $\Sigma(3,n)$. However it can only be exhibited as an $(\alpha,\beta)$-weakly renormalizable word by $w=w_1w_2$. This is why we introduce the notion of $(\alpha, \beta)$-semi renormalizable in \Cref{def:semirenormalizable}. Indeed, the previous example $w$ is $(a, b)$-semi renormalizable (with nontrivial kernel), since $2w1 \in \langle a, b\rangle$.
	
	With these considerations, we will now present a renormalization algorithm: if we have a $(u, v)$-weakly renormalizable word with a nonempty renormalization kernel, we can exhibit this word as being $(\alpha, \beta)$-weakly renormalizable for $(\alpha, \beta) \in \{(uv, v), (u, uv)\}$ chosen appropriately.
	
	\begin{lemma}[Renormalization algorithm]\label{lem:renormalizingr}
		Let $w \in \Sigma(3 + e^{-r}, |w|)$ satisfying $\sizer(w)\leq r$. If $w$ is $(u, v)$-weakly renormalizable as $w = w_1 \gamma w_2$ with $\gamma \neq \emptyset$, then $w$ is $(\alpha, \beta)$-weakly renormalizable for some $(\alpha, \beta) \in \{(uv, v), (u, uv)\}$. Moreover, if $\gamma$ starts with $u$ or ends with $v$, then $w_1$ or $w_2$, respectively, does not change for the renormalization with alphabet $(\alpha,\beta)$.
	\end{lemma}
	
	Before proving the previous lemma, we will discuss the intuition behind this algorithm. The main inspiration is the ``exponent-reducing'' procedure discussed in \Cref{sec:bombieri}. Indeed, if a word $w$ is $(u, v)$-weakly renormalizable, then it is of the form $w = w_1 \gamma w_2$, where $\gamma$ is written in terms of $u$ and $v$. The word $\gamma$ cannot contain factors of the form $uu \ldots vv$ or $vv \ldots uu$ (as discussed in the proof below), so it is written as powers of $u$ (respectively, $v$) followed by single instances of $v$ (respectively, $u$). Hence, we can choose a new alphabet $(\alpha, \beta) = (u, uv)$ (respectively, $(\alpha, \beta) = (uv, v)$) so that all exponents are now reduced by $1$ when $\gamma$ is written in the new alphabet $(\alpha, \beta)$. This simplifies the structure of the renormalization kernel at the cost of making the alphabet more complex. The renormalization algorithm should be, hence, applied inductively a certain number of times to ensure that the complexity of both the renormalization kernel and the alphabet remain reasonable (see for example \Cref{cor:renormalizing_3n} and the proof of \Cref{thm:equalities} to see how this is used).

	\begin{proof}
		
		We will explicitly exhibit $w$ as being $(\alpha, \beta)$-renormalizable as $w = \tilde{w}_1 \tilde\gamma \tilde{w}_2$ for some $(\alpha, \beta) \in \{(uv, v), (u, uv)\}$.
		
		By \Cref{lem:aa_bbr} and the comments after \Cref{lem:symmetry2}, some patterns on a weakly renormalizable word imply that $w \notin \Sigma(3 + e^{-r}, |w|)$, and so are forbidden: this holds if $\gamma$ contains both the factors $uu$ and $vv$ (in any order), and also in the following situations:
		
		\begin{enumerate}
			\item If $(u,v)=(\eta,\eta\theta)$ for some $(\eta, \theta) \in \overline{P}$, $\gamma$ starts with $v$ and contains the factor $uu$, and $|w_1|\ge |u|$.
			\item If $(u,v)=(\eta\theta,\theta)$ for some $(\eta, \theta) \in \overline{P}$, $\gamma$ ends with $u$ and contains the factor $vv$, and $|w_2|\ge |v|$.
		\end{enumerate}

		We first assume that $w$ does not contain the factor $vv$ and we analyze the following subcases (where $s$ and $e_j$ are positive integers for $1 \leq j \leq k$):
		
		\noindent\textbf{Case 1:} If $\gamma = u^{e_1} v u^{e_2} v \ldots u^{e_k} v$, we take $\alpha = u$, $\beta = uv$ and
		\[
		\tilde\gamma = \alpha^{e_1-1}\beta\alpha^{e_2-1}\beta\ldots\alpha^{e_k-1}\beta, \quad \tilde{w}_1 = w_1, \quad \tilde{w}_2 = w_2.
		\]
		Indeed, $\tilde{w}_1 = w_1$ is a suffix of $uv$ by hypothesis, so it is also a suffix of $\alpha\beta = u^2v$. Moreover, $\tilde{w}_2 = w_2$ is a prefix of $uv$ by hypothesis and to show that it is also a prefix of $\alpha\beta$ we consider two cases. If $|w_2| < |v|$, then $w_2$ is a prefix of $\alast{v}$, since $uv$ starts with $\alast{v}$ by \Cref{lem:structure}. The same lemma also shows that $u^2v$ starts with $\alast{v}$, so $w_2$ is a prefix of $\alpha\beta = u^2v$. Otherwise, we must have $|w_2| < |u|$, since $|w_2| < \max\{|u|, |v|\}$. Thus, $\tilde{w}_2 = w_2$ is a proper prefix of $u$ and, hence, of $\alpha\beta = u^2v$.
		
		\smallbreak
		
		\noindent\textbf{Case 2:} If $\gamma = v u^{e_1} v u^{e_2} v \ldots u^{e_k} v$ we consider two cases. If $|w_1| < |u|$, we take $\alpha = u$, $\beta = uv$ and
		\[
		\tilde\gamma = \alpha^{e_1 - 1} \beta \alpha^{e_2 - 1}\beta \ldots \alpha^{e_k-1}\beta, \quad \tilde{w}_1 = w_1 v, \quad \tilde{w}_2 = w_2.
		\]
		Indeed, recall that $uv$ ends with $\bfirst{u}$ and that $w_1$ is a suffix of $uv$. Since $|w_1| < |u|$, $w_1$ is also a suffix of $u$ (as $\bfirst{u}$ and $u$ are equal up to the first letter). We obtain that $w_1$ is a suffix of $u$, so $\tilde{w}_1 = w_1v$ is a suffix of $\alpha\beta = u^2v$. Moreover, $\tilde{w}_2$ is a prefix of $\alpha\beta = u^2v$ by the exact same proof of the previous case: it is either shorter than $v$ (in which case it is a proper prefix of $\alast{v}$ and, hence, of $u^2v$ by \Cref{lem:structure}), or shorter than $u$ (in which case it is a prefix of $u$ and, hence, of $u^2v$).
		
		Otherwise, we have $|u| \leq |w_1| < |v|$, so $(u, v) = (\eta, \eta\theta)$ for some pair $(\eta, \theta) \in \overline{P}$. Since $w_1$ is a suffix of $uv = \eta^2\theta$ and $|w_1| \geq |u| = |\eta|$, we have that $w_1$ ends with $\bfirst{\eta}$ by \Cref{lem:structure}. If $e_j > 1$ for some $1 \leq j \leq k$, then $w$ contains a factor of the form $\bfirst{\eta}v \ldots uu\alast{\theta}$. In fact, since $w_1$ ends with $\bfirst{\eta}$, we have that $w$ contains a word of the form $w' = \bfirst{\eta}v \ldots u^{e_j-2}uuv$, where $1 \leq j \leq k$ is chosen so $e_j > 1$. Moreover, $v = \eta\theta$ starts with $\alast{\theta}$ by \Cref{lem:structure}, so $w'$ contains, in turn, a word of the form $\bfirst{\eta}v \ldots uu\alast{\theta}$. This contradicts that $w \in \Sigma(3 + e^{-r}, |w|)$ by \Cref{lem:aa_bbr,lem:symmetry,lem:symmetry2}.
		
		We assume then that $e_j = 1$ for every $1 \leq j \leq k$ and take $\alpha = uv$, $\beta = v$ and
		\[
		\tilde\gamma = \beta\alpha^k, \quad \tilde{w}_1 = w_1, \quad \tilde{w}_2 = w_2.
		\]
		Indeed, we have that $\tilde{w}_2 = w_2$ is a prefix of $uv^2$ since it is a prefix of $uv$. Moreover, if $|w_1| < |v|$ then $\tilde{w}_1 = w_1$ is a suffix of $uv^2$ as it is a suffix of $uv$, and if $|w_1| < |u|$ then $w_1$ is a proper suffix of $\bfirst{u}$ (by \Cref{lem:structure}), so it is also a suffix of $uv^2$ (by \Cref{lem:structure} again). Finally, since $\tilde{\gamma}$ starts with $\beta$ and $(\alpha, \beta) = (uv,v)$, we have to check that $|u| \leq |\tilde{w}_1| = |w_1|$, but this holds by hypothesis.
		
		\smallbreak
		
		\noindent\textbf{Case 3:} If $\gamma = u^{e_1} v u^{e_2} v \ldots u^{e_k} v u^s$, we must have that $|v| \leq |u w_2|$. Indeed, if $|v| > |u w_2|$, then $(u, v) = (\eta, \eta\theta)$ for some alphabet $(\eta, \theta)$. Since $\gamma$ ends with $u$, by definition of $(u, v)$-weakly renormalizability we have that $|\theta| \leq |w_2|$. Hence, $|v| = |\eta\theta| = |u\theta| \leq |uw_2|$, a contradiction.
		
		Let $r \in \{0,1\}$ then be such that $|v| \leq |u^r w_2| < |uv| \leq |u^{r+1}w_2|$. Then, we choose $\alpha = u$, $\beta = uv$ and
		\[
		\tilde{\gamma} = \alpha^{e_1-1} \beta \alpha^{e_2-1}\beta \ldots \alpha^{e_k-1}\beta \alpha^{s-r}, \quad \tilde{w}_1 = w_1, \quad \tilde{w}_2 = u^r w_2.
		\]
		
		Indeed, $\tilde{w}_1 = w_1$ is a suffix of $\alpha\beta = u^2v$ since it is a suffix of $uv$. 
		Now, if $r = 0$, then $|w_2| < |u|$ since $|w_2| < \max\{|u|, |v|\}$ and $|v| \leq |w_2|$ by hypothesis. Since $w_2$ is a prefix of $uv$, it is actually a prefix of $u$ and, hence, of $\alpha\beta = u^2v$. If $r = 1$, we have that $w_2$ is a prefix of $uv$ and, thus, $\tilde{w}_2 = u w_2$ is a prefix of $\alpha\beta = u^2v$.
		
		Since $\tilde{\gamma}$ ends with $\alpha$ if $r = 0$, we have to check that $|v| \leq |\tilde{w}_2|$. This holds since $|v| \leq |\tilde{w}_2| = |w_2|$ in this case.
		
		\smallbreak
		
		\noindent\textbf{Case 4:} Finally, if $\gamma = vu^{e_1}vu^{e_2}v\ldots u^{e_k}vu^s$, we combine the discussions of the previous two cases. More precisely, we assume first that $|u| \leq |w_1| < |v|$. If $e_j > 1$ for some $1 \leq j \leq k$ or $s > 1$, then we obtain a contradiction with the hypothesis that $w \in \Sigma(3 + e^{-r}, |w|)$ by \Cref{lem:aa_bbr,lem:symmetry,lem:symmetry2}. Indeed, in this case we have that $(u, v) = (\eta, \eta\theta)$ for some $(\eta, \theta) \in \overline{P}$, so $w$ contains a factor of the form $\bfirst{\eta} v \ldots uu \alast{\theta}$ as in the second case if $e_j > 1$ for some $1 \leq j \leq k$. On the other hand, if $s > 1$, then $w$ contains a word of the form $w' = \bfirst{\eta}v \ldots u^{s-2}uu w_2$. Now, observe that the fact that $\gamma$ ends with $u$ and the definition of $(u, v)$-renormalizability imply that $w_2$ starts with $\alast{\theta}$. Hence, $w'$ contains a word of the form $\bfirst{\eta}v \ldots uu \alast{\theta}$. This leads to the same contradiction with \Cref{lem:aa_bbr,lem:symmetry,lem:symmetry2}.
		
		In the case where $|u| \leq |w_1| < |v|$, $e_j = 1$ for every $1 \leq j \leq k$ and $s = 1$, we take $\alpha = uv$, $\beta = v$ and
		\[
		\tilde{\gamma} = \beta \alpha^k, \quad \tilde{w}_1 = w_1, \quad \tilde{w}_2 = u w_2.
		\]
		Since $w_1$ is a suffix of $uv$, then $\tilde{w}_1 = w_1$ is a suffix of $\alpha\beta = u^2v$. Now, observe that $|u| \leq |w_1| < |v|$ implies that $|w_2| < |v|$, since by hypothesis we have that $|w_1|, |w_2| < \max\{|u|, |v|\}$. Thus, by \Cref{lem:structure}, $w_2$ is a proper prefix of $\alast{v}$, so it is also a prefix of $v$. We then obtain that $\tilde{w}_2 = uw_2$ is a prefix of $\alpha\beta = uv^2$.
		
		Otherwise, if $|w_1| < |u|$ we take $\alpha = u$, $\beta = uv$ and argue as in the third case. More precisely, let $r \in \{0,1\}$ be such that $|v| \leq |u^r w_2| < |uv| \leq |u^{r+1}w_2|$ and take
		\[
		\tilde{\gamma} = \alpha^{e_1-1}\beta\alpha^{e_2-1}\beta \ldots \alpha^{e_k-1}\beta\alpha^{s-r}, \quad \tilde{w}_1 = w_1 v, \quad \tilde{w}_2 = u^rw_2.
		\]
		We have that $\tilde{w}_2 = u^rw_2$ is a prefix of $\alpha\beta$ by the same arguments of the third case, and $r$ is chosen so $|v| \leq |\tilde{w}_2|$. Moreover, $\tilde{w}_1 = w_1v$ is a suffix of $\alpha\beta = u^2v$ since \Cref{lem:structure} and the fact that $|w_1| < |u|$ imply that $w_1$ is a proper suffix of $\bfirst{u}$, so that it is also a suffix of $u$. This finishes the last subcase.
		
		\smallbreak
		
		We now assume that $w$ contains the factor $vv$, so, in particular, it does not contain the factor $uu$. We analyze the following subcases (where $s$ and $e_j$ are a positive integers for $1 \leq j \leq k$):
		
		\noindent\textbf{Case 1:} If $\gamma = u v^{e_1} u v^{e_2} \ldots u v^{e_k}$, we take $\alpha = uv$, $\beta = v$ and
		\[
		\tilde{\gamma} = \alpha \beta^{e_1-1} \alpha \beta^{e_2-1} \ldots \alpha \beta^{e_k-1}, \quad \tilde{w}_1 = w_1, \quad \tilde{w}_2 = w_2.
		\]
		
		\smallbreak
		
		\noindent\textbf{Case 2:} If $\gamma = u v^{e_1} u v^{e_2} \ldots u v^{e_k} u$, we take $\alpha = uv$, $\beta = v$ and
		\[
		\tilde{\gamma} = \alpha \beta^{e_1 - 1} \alpha \beta^{e_2 - 1} \ldots \alpha \beta^{e_k - 1}, \quad \tilde{w}_1 = w_1, \quad \tilde{w}_2 = u w_2.
		\]
		
		\smallbreak
		
		\noindent\textbf{Case 3:} If $\gamma = v^s u v^{e_1} u e^{e_2} \ldots u v^{e_k}$, we take $r \in \{0,1\}$ such that $|u| \leq |w_1 v^r| < |uv| \leq |w_1 v^{r+1}|$, and define $\alpha = uv$, $\beta = v$ and
		\[
		\tilde{\gamma} = \beta^{s-r} \alpha \beta^{e_1-1} \alpha \ldots \beta^{e_{k-1}-1} \alpha \beta^{e_k-1}, \quad \tilde{w}_1 = w_1 v^r, \quad \tilde{w}_2 = w_2.
		\]
		
		\smallbreak
		
		\noindent\textbf{Case 4:} If $\gamma = v^s u v^{e_1} u e^{e_2} \ldots u v^{e_k} u$, we take $r \in \{0,1\}$ such that $|u| \leq |w_1 v^r| < |uv| \leq |w_1 v^{r+1}|$, and define $\alpha = uv$, $\beta = v$ and
		\[
		\tilde{\gamma} = \beta^{s-r} \alpha \beta^{e_1-1} \alpha \ldots \beta^{e_{k-1}-1} \alpha \beta^{e_k-1}, \quad \tilde{w}_1 = w_1 v^r, \quad \tilde{w}_2 = uw_2.
		\]
		
		\smallbreak
		
		Observe that the cases where $e_j = 1$ for all $1 \leq j \leq k$ cannot arise in the previous subcases, since we are explicitly assuming that $w$ contains the factor $vv$. The arguments showing that these choices satisfy the definition of $(\alpha, \beta)$-renormalizabilty are analogous to those of the previous cases (where the factor $vv$ was not present). Thus, this concludes the proof.
	\end{proof}
	
	Once again, this lemma could be stated in terms of the length of $w$, as in the following corollary:
	\begin{corollary} \label{cor:renormalizing_length}
		Let $w \in \Sigma(3 + (3 + 2\sqrt{2})^{-(|w|+1)}, |w|)$ be a finite word. If $w$ is $(u, v)$-weakly renormalizable as $w = w_1 \gamma w_2$ with $\gamma \neq \emptyset$, then $w$ is $(\alpha, \beta)$-weakly renormalizable for some $(\alpha, \beta) \in \{(uv, v), (u, uv)\}$.
	\end{corollary} 
	\begin{proof}
		We have that $\sizer(w) \leq (n+1) \log(3 + 2\sqrt{2})$ by \Cref{lem:rcomparedwithsize}. Taking $r = (n+1) \log(3 + 2\sqrt{2})$ and applying the previous lemma we obtain this result.
	\end{proof}
	
	We will now present a series of corollaries of the renormalization algorithm. We start with the version that is needed for the proof of \Cref{thm:equalities}.
	
	\begin{corollary} \label{cor:renormalizing_3n}
		Let $n \geq 68$ and let $w \in \Sigma(3 + 6^{-3n}, 3n)$. Then, there exists  an alphabet $(\alpha, \beta) \in \overline{P}$ satisfying $|\alpha|, |\beta| < n$ and $|\alpha\beta| \geq n$ such that $w$ is $(\alpha, \beta)$-semi renormalizable.
	\end{corollary}
	
	\begin{proof}

		Since $n \geq 68$, \Cref{lem:no_odd_ones_or_twos} holds, so, possibly up to adding one letter to the left and one to the right, $w$ is a word in the alphabet $\{a,b\}$. As previously discussed, $w$ is trivially $(a,b)$-weakly renormalizable with $w_1 = w_2 = \varnothing$ and $\gamma = w$. Observe that $w$ satisfies the first hypothesis of \Cref{cor:renormalizing_length}. Indeed, this follows from the fact that $6^{-3n} \leq (3 + 2\sqrt{2})^{-(3n+3)}$ for every $n \geq 61$. By \Cref{cor:renormalizing_length}, we can apply the renormalization algorithm inductively as long as the renormalization kernel is nonempty; this produces a a finite sequence of alphabets. We will show that the sought-after alphabet $(\alpha, \beta)$ is the first alphabet in the sequence that satisfies $|\alpha\beta| \geq n$.
		
		We will first show that such an alphabet exists. Assume that $(u, v) \in \overline{P}$ is one of the alphabets of the sequence. If $|uv| < n$, then $|w_1|, |w_2| < n$, since
		\[
		|w_1|, |w_2| \leq \max\{|u|, |v|\} < |uv| < n,
		\]
		where $w_1$ and $w_2$ are the words obtained in this step of the algorithm by the decomposition $w = w_1 \gamma w_2$. Hence, $\gamma \neq \varnothing$, since
		\[
		|\gamma| = |w| - |w_1| - |w_2| > 3n - n - n = n.
		\]
		Thus, we can apply the algorithm again if $|uv| < n$. Since the length of an alphabet increases with each inductive application of the algorithm, we will eventually find an alphabet $(\alpha, \beta) \in \overline{P}$ satisfying $|\alpha\beta| \geq n$. Assume that $(\alpha, \beta) \in \overline{P}$ is the first alphabet in the sequence satisfying this condition.
		
		It remains to show that $|\alpha|, |\beta| < n$. Assume by contradiction that this is false. Assume further that $(\alpha, \beta) = (uv, v)$ for some alphabet $(u, v) \in \overline{P}$; the case where $(\alpha, \beta) = (u, uv)$ is similar. We then have that $|\alpha| \geq n$.
		
		Observe that the alphabet $(u, v)$ satisfies that $|uv| = |\alpha|$, which contradicts that $(\alpha, \beta)$ is the first alphabet in the sequence of inductive applications of the algorithm satisfying this inequality. Thus, the proof of the corollary is complete.
		
	\end{proof}
	
	\begin{remark}
		Clearly, the previous corollary holds for any $w \in \Sigma(3 + B^{-3n}, 3n)$, where $B > 3 + 2\sqrt{2}$ and $n \in \NN^*$ is large enough (depending on $B$).
	\end{remark}
	
	The following corollaries are straightforward consequences of the renormalization algorithm and are thus presented here. Nevertheless, they are not used in the proof of \Cref{thm:equalities} and will be only used in the next section. Recall that a word belongs to $\Sigma^{(r-2)}(3+e^{-r})$ if it belongs to both $\Sigma(3 + e^{-r}, |w|)$ and $Q_{r-2}$.
	
	\begin{corollary}
		Let $r\in\NN$ and let $w \in \Sigma^{(r-2)}(3 + e^{-r})$ be a finite word. If $w$ is $(u, v)$-weakly renormalizable as $w_1 \gamma w_2$ with $\gamma \neq \varnothing$, then $w$ is $(\alpha, \beta)$-weakly renormalizable for some $(\alpha, \beta) \in \{(uv, v), (u, uv)\}$.
	\end{corollary}
	
	\begin{proof}
		Observe that if $w=c_1 \ldots c_n \in Q_r$, then $\sizer(c_1 \ldots c_{n-1})\leq r-3$, so
		\[
		\sizes(w)^{-1}\leq 2\sizes(c_n)^{-1}\sizes(c_1 \ldots c_{n-1})^{-1}\leq 12e^{r-2},
		\]
		which implies that $\sizer(w)\leq r$. We then use \Cref{lem:renormalizingr}.
	\end{proof}
	
	\begin{corollary} \label{cor:r/6}
		Let $(\alpha,\beta)\in  \overline{P}$ and let $w \in \Sigma^{(r-2)}(3 + e^{-r})$ be a word starting with $\alpha$ or $\beta$. Then, by extending $w$ by at most one digit to the right, $w$ is $(\alpha, \beta)$-weakly renormalizable for some alphabet $(\alpha, \beta)$ satisfying $|\alpha\beta| \geq r/6$. 
	\end{corollary}
	
	\begin{proof}
		<        By \Cref{lem:1} we know that $w$ does not contain ``internal'' blocks of $1$'s or $2$'s with odd length, that is, words of the form $c' c^n c'$ for $c, c' \in \{1, 2\}$ with $c \neq c'$ for some odd $n \in \NN$. Since $w$ starts with $\alpha$ or $\beta$, it starts with an even block as well. On the other hand, $w$ can possibly end with an odd block of $1$'s or $2$'s. If $w$ ends with an odd block of $2$'s, then $w=\gamma w_2$ is $(a,b)$-weakly renormalizable where $\gamma\in\langle a,b\rangle$ and $w_2=2$. In case that it ends with an odd block of $1$'s, we just need to extend $w= c_1 \ldots c_n$ to $\tilde{w}=w1=c_1 \ldots c_n 1$. In this case 
		\begin{equation*}
			\sizes(c_1 \ldots c_{n} 1)^{-1}\leq2\sizes(c_1 \ldots c_{n-1})^{-1}\sizes(1,1)^{-1}\leq 12e^{r-2}
		\end{equation*}
		which gives $\sizer(c_1 \ldots c_{n} 1)\leq r$. 
		
		We claim that $w$ or $\tilde{w}$ is $(\alpha,\beta)$-weakly renormalizable for an alphabet $(\alpha, \beta)$ satisfying $|\alpha|,|\beta|<|w|$ with $|\alpha\beta|\geq |w|/2$. Indeed, if $|\alpha\beta|<|w|/2$, then writing $w=w_1\gamma w_2$ gives $|w_1|+|w_2|<2|\alpha\beta|\leq |w|$. We obtain that $\gamma\neq\emptyset$, so we can continue applying the algorithm. Here, we skipped most details as this is very similar to the proof of \Cref{cor:renormalizing_3n}.
		
		We remark that if, for some iteration of the algorithm, we obtain $\gamma=\alpha^r$ (respectively, $\gamma=\beta^r$), then the algorithm increases the size of the alphabet, but does not change the renormalization kernel $\gamma$. In these cases, we have that $w$ is a subword of $\alpha\beta\alpha^r\alpha\beta$ and of $\alpha^{r+1}\beta\alpha^{r+1}\beta$ (respectively, of $\alpha\beta\beta^r\alpha\beta$ and of $\alpha\beta^{r+1}\alpha\beta^{r+1}$), and so it belongs to $\Sigma(3,|w|)$.
		
		Now, let $(\alpha, \beta)$ be such an alphabet. If $r \leq 24$, then $|\alpha\beta| \geq 4 \geq r/6$ since $|\alpha|, |\beta| \geq 2$. If $r > 24$, we have that
		\[
		|\alpha\beta|\geq |w|/2\geq (r-2)/(2\log(3+2\sqrt{2}))-1/2\geq r/6,
		\]
		where we are using \Cref{lem:rcomparedwithsize}.
	\end{proof}
	
	\begin{corollary}\label{lem:local_extension}
		Let $(\alpha, \beta) \in \overline{P}$ with $|\alpha\beta|< r/6$ and $w\in\Sigma^{(r-2)}(3+e^{-r})$. If $w$ contains $\alpha\beta$, then $w$ is $(\alpha,\beta)$-semi renormalizable, say $\tilde{w}=w_1\gamma w_2$. Moreover, if $w$ starts (ends) with $\alpha\beta$, then $w_1=\emptyset$ ($w_2=\emptyset$).
	\end{corollary}
	
	\begin{proof}
		First note that $w$ is trivially $(a,b)$-semi renormalizable, say $\tilde{w}=\gamma_0$ where $\gamma_0\in\langle a,b\rangle$.  Now we apply inductively the renormalization algorithm (\Cref{lem:renormalizingr}) to obtain a sequence of alphabets $(A_j,B_j)\in\overline{P}_j$ such that for all $0\leq j\leq m$, the word $\tilde{w}$ is $(A_j,B_j)$-weakly renormalizable for each $j$ and $|A_mB_m|\geq r/6$. 
		
		On the other hand, since $(\alpha,\beta)\in\overline{P}$ there exists a sequence of alphabets $(\alpha_i,\beta_i)\in\overline{P}_i$ such that $\alpha\beta\in\langle\alpha_i,\beta_i\rangle$ for all $0\leq i\leq n$ and $(\alpha_n,\beta_n)=(\alpha,\beta)$. Since $\alpha\beta$ starts with $a=\alpha_0$ and ends with $b=\beta_0$ (\Cref{lem:structure}), inductively we obtain that $\alpha\beta$ starts with $\alpha_i$ and ends with $\beta_i$. In particular $\alpha\beta$ contains $\alpha_i\beta_i$.  
		
		Write $\tilde{w}=w_1\gamma_jw_2$ as in the definition of $(A_j,B_j)$-weakly renormalizable. Using the fact that $\alpha\beta$ contains $\alpha_j\beta_j$, gluing some words $\tau$ and $\tau'$ we get
		\[
		\tau\alpha_j\beta_j\tau'=A_jB_j\gamma_jA_jB_j\in\langle A_j,B_j\rangle,
		\]
		hence by \Cref{lem:determine_alphabet} we obtain that $(A_j,B_j)=(\alpha_j,\beta_j)$ for all $0\leq j\leq n$. In particular $m>n$, because otherwise  $r/6\leq|A_mB_m|=|\alpha_m\beta_m|<r/6$. This shows that $\tilde{w}$ is $(\alpha,\beta)$-weakly renormalizable.

		Now assume that $w$ starts with $\alpha\beta$ (the other case is analogous). Observe that there is no need to complete the word to the left. We will show that $w_1=\emptyset$ for all $0\leq j\leq n$. Note that we already showed that $w_1$ is empty for $(\alpha_0,\beta_0)=(a,b)$. If $w_1$ becomes nonempty for $k+1$ for some $0\leq k\leq n$, it must happen that $\tilde{w}=\gamma_k w_2$ starts with $\beta_k$ (because of the renormalization algorithm). But $w$ starts with $\alpha\beta$, which in turn starts with $\alpha_k^s\beta_k$, which leads to a contradiction because it starts with $(\beta_k)_a$ by \Cref{lem:structure}. Since $(\alpha_n,\beta_n)=(\alpha,\beta)$ this finishes the proof.
	\end{proof}
	
	Finally, to end this section we prove \Cref{thm:equalities}.
	\begin{proof}[\emph{\textbf{Proof of \Cref{thm:equalities}}}]
		Consider $n \geq 68$, then we claim that $\Sigma(3 + 6^{-3n}, n)=\Sigma(3, n)$. Indeed, let $\theta$ an element of $\Sigma(3 + 6^{-3n}, n)$. By definition, $\theta$ can be continued indefinitely to the left and right so, in particular, there exists a word $\tau \in \Sigma(3 + 6^{-3n}, 3n)$ obtained by gluing words of size $n$ at each side of $\theta$. By \Cref{cor:renormalizing_3n}, there exists $(\alpha,\beta)\in \overline P$ with $|\alpha|, |\beta|<n$ and $|\alpha\beta|\ge n$ such that $\tau$ is $(\alpha,\beta)$-semi renormalizable. Writing $\tilde{\tau}=w_1\gamma w_2$ as in the definition of weak renormalization, we have $|w_1|, |w_2|<\max\{|\alpha|, |\beta|\}<n$, so $\theta$ is a factor of $\gamma$. Considering the smallest sequence $\eta$ of $(\alpha,\beta)$-letters of $\gamma$ containing $\theta$ as a factor, the sequence obtained by removing the first and the last $(\alpha,\beta)$-letter of $\eta$ has size smaller than $n$ and thus cannot contain $\alpha\beta$ or $\beta\alpha$ as factors, and thus $\eta$ is of the form $\alpha^r$, $\beta^r$, $\alpha^r \beta$, $\beta^r \alpha$, $\beta \alpha^r \beta$, $\alpha \beta^r \alpha$, $\beta \alpha^r$ or $\alpha \beta^r$ for some positive integer $r$. 
		
		In any of these cases, $\eta\in\Sigma(3,|\eta|)$. Indeed, since $(\alpha, \beta) \in \overline{P}$, all of these words are factors of words in $\mathsf{c}(\overline{P}) = P$ (where recall that $\mathsf{c}$ is the concatenation operator $\mathsf{c}(u,v) = uv$). Since, by \Cref{lem:bombieri_old}, the set of factors of words in $P$ coincides with the  set of words $w$ satisfying $w \in \Sigma(3, |w|)$, we obtain that $\eta\in\Sigma(3,|\eta|)$. Therefore, $\theta\in\Sigma(3,n)$.
		
		To complete our proof, we need to show that, for every sufficiently large integer $n$, we have $\Sigma(3 - 6^{-3n}, n)=\Sigma(3, n)$. Indeed, given	$w \in \Sigma(3, n)$, by \Cref{lem:bombieri_old,lem:fundamental_period_bounded}, there exists $\Pi \in P$ containing $w$ such that $|\Pi| \leq 3|w|$. Since $(3+2\sqrt{2})^{3} < 6^3$, if $n$ is sufficiently large then  \Cref{lem:R_S_bound} shows that $\Pi^\infty \in \Sigma(3 - 6^{-3n})$, so $w \in \Sigma(3 - 6^{-3n}, n)$.
	\end{proof}
	
	\section{Improving the estimates} \label{sec:improving}
	
	\medskip

	Bombieri \cite[Lemma 13]{B:markoff_numbers} characterized the words in $\Sigma(3)$ by stating the conditions that the sequence of exponents $(e_i)_{i \in \mathbb{Z}}$ has to satisfy for a Type I or Type II bi-infinite word to belong to $\Sigma(3)$ (where we are using the terminology of \Cref{sec:bombieri}). We begin this section by stating an analog of this fact for words in $\Sigma(3 + e^{-r}, n)$. The proof is essentially applying the renormalization algorithm to a word of the form $w=\alpha^{e_i}\beta\alpha^{e_{i+1}}\beta$ or $w=\beta^{e_i}\alpha\beta^{e_{i+1}}\alpha$, but we need to be careful about the magnitude of $\sizer(w)$.  
	
	\begin{lemma}[Bombieri's characterization]\label{lem:bombieri}
		Let $(\alpha,\beta)\in\overline{P}$. Consider a word $\gamma$ of the form $\gamma=\alpha^{e_0}\beta\alpha^{e_1}\beta\ldots\beta\alpha^{e_\ell}$ or $\gamma=\beta^{e_0}\alpha\beta^{e_1}\alpha\ldots\alpha\beta^{e_\ell}$ with $e_i\geq 1$ for all $1\leq i\leq \ell-1$. Assume that $\gamma \in \Sigma(3+e^{-r},|\gamma|)$ and let $\theta=\alpha$ in the first case and $\theta=\beta$ in the second case. If $\sizer(\theta^{e_i})\leq r-2|\alpha\beta|$, then 
		\begin{itemize}
			\item for $1\leq i\leq\ell-2$, we have $|e_i-e_{i+1}|\leq 1$.
			\item for $i=0$, we have $e_1\geq e_0-1$  when $\theta=\alpha$. When $\theta=\beta$, if moreover $\sizer(\beta^{e_0})\leq r-6|\alpha\beta|$ or $|\alpha|\leq|\beta|$, then $e_1\geq e_0-1$.
			\item for $i=\ell-1$, we have $e_\ell\leq e_{\ell-1}+1$  when $\theta=\beta$. When $\theta=\alpha$, if moreover $\sizer(\alpha^{e_{\ell-1}})\leq r-6|\alpha\beta|$ or $|\beta|\leq|\alpha|$, then $e_\ell\leq e_{\ell-1}+1$.
		\end{itemize}
	\end{lemma}
	
	Before proceeding with the proof, we must comment why we need $\sizer(\theta^{e_i})$ to be smaller at the end of the word in the last two bullet points. Observe that if $\beta=\alpha^sv$ for some $(\alpha,v)\in\overline{P}$, then clearly $e_\ell$ can be much larger than $e_{\ell-1}$, because all powers $\alpha^{e_\ell-e_{\ell-1}}$ could belong to the (potential) next letter $\beta$. Similarly, when $\alpha=u\beta^s$ for some $(u,\beta)\in\overline{P}$, the power $\beta^{e_0-e_1}$ could belong to the (potential) preceding letter $\alpha$. 
	
	\begin{proof}
		Let $\omega$ be a bi-infinite word containing $\gamma$ and such that $\omega \in \Sigma(3 + e^{-r})$. Note that if $\{\theta,\Tilde{\theta}\}=\{\alpha,\beta\}$ then $\sizer(\theta^{e_i}\Tilde{\theta})< r-2|\alpha\beta|+2|\theta|+4\leq r$ by \Cref{lem:rcomparedwithsize}.
		
		Suppose $\gamma=\alpha^{e_0}\beta\alpha^{e_1}\beta\ldots$. Take $k\leq e_{i+1}$ maximal such that $\sizer(\alpha^{e_i}\beta\alpha^k)\leq 2r$. If $e_i\geq k$ then $\sizer(\alpha^k)\leq r-2|\alpha\beta|$ as well, so actually $k=e_{i+1}$ because, otherwise,
		\begin{align*}
			\sizer(\alpha^{e_i}\beta\alpha^{k+1})&\leq \sizer(\alpha^{e_i})+\sizer(\beta\alpha)+\sizer(\alpha^{k})+4 \\
			&\leq r-2|\alpha\beta|+\sizer(\beta\alpha)+r-2|\alpha\beta|+4\leq 2r,
		\end{align*}
		where we used \Cref{lem:rcomparedwithsize} to guarantee that $\sizer(\alpha\beta)\leq 1.8|\beta\alpha|+1.8$. Similarly, we use $\sizer(\beta)\leq 1.8|\beta|+1.8$ (for $\beta=b$  use $\sizer(b)=1$ instead) to get
		\[
		\sizer(\alpha^{e_i}\beta\alpha^{e_{i+1}}\beta)\leq 2r-4|\alpha\beta|+2\sizer(\beta)+6\leq 2r.
		\]
		
		Hence, letting $(\tilde{\alpha},\tilde{\beta})=(\alpha,\alpha^{e_{i+1}}\beta)$ we obtain that
		\[
		\tilde{\gamma}=\alpha^{e_i}\beta\alpha^{e_{i+1}}\beta=\tilde{\alpha}^{e_i-e_{i+1}}\tilde{\beta}\tilde{\beta}
		\]
		is a subword of a word $\omega\in\Sigma(3+e^{-r})$, so if $e_i-e_{i+1}\geq 2$ it will contradict the second part of \Cref{lem:aa_bbr}. If $e_i<k$, then let $(u,v)=(\alpha,\alpha^{e_i-1}\beta)$ and $(\tilde{\alpha},\tilde{\beta})=(u,uv)$ thus by the first case of \Cref{lem:symmetry2}
		\[
		\alpha\beta\alpha^{e_i}\beta\alpha^{e_{i+1}}\beta=\alpha\beta\tilde{\beta}\alpha^{e_{i+1}-e_i}\tilde{\beta}=\beta_a u^b\tilde{\beta}\tilde{\alpha}^{e_{i+1}-e_i}v_au^b=\beta_a(\tilde{\alpha}^{e_{i+1}-e_{i}}\tilde{\beta}\tilde{\beta})^\ast \alpha^b
		\]
		is a subword of $\gamma$ when $i<\ell-1$. If $e_{i+1}-e_i\geq 2$, then we would have that $\tilde{\alpha}\tilde{\alpha}\tilde{\beta}\tilde{\beta}$ is a subword of $\gamma^\ast$ with $\sizer(\tilde{\alpha}\tilde{\alpha}\tilde{\beta}\tilde{\beta})=\sizer(\alpha^{e_i+2}\beta\alpha^{e_i}\beta)\leq 2r$, which contradicts \Cref{lem:aa_bbr}. This finishes the first bullet point for $\theta=\alpha$. In the particular case where $i=\ell-1$, we do not necessarily have $\beta$ after $\alpha^{e_\ell}$. If $|\beta|\leq|\alpha|$, then $(\alpha^b)^\ast(\tilde{\alpha}^{e_{i+1}-e_{i}}\tilde{\beta}\tilde{\beta})(\beta_a)^\ast$ is a subword of $\gamma^\ast$ after removing a $\beta^\ast$ at the beginning, so we still get that $\tilde{\alpha}\tilde{\alpha}\tilde{\beta}\tilde{\beta}$ is a subword of $\gamma^\ast$.
		
		When $(\alpha,\beta)=(u,uv)$ we need to extend the word $\beta\alpha^{e_{\ell-1}}\beta\alpha^{e_\ell}$ by using \Cref{lem:local_extension}. We will extend this word to the left and then to the right. Since $|uv|=|\beta|<r/6$ and $\sizer(\beta)<r-2$ (because $0\leq \sizer(\alpha^{e_{\ell-1}})\leq r-6|\alpha\beta|$), consider the $(u,v)$-semi renormalizable continuation $w\in\Sigma^{(r-2)}(3+e^{-r})$ inside $\omega$ that contains and ends in the leftmost $\beta$ of $\beta\alpha^{e_{\ell-1}}\beta\alpha^{e_\ell}$ (one begins with $\beta$ and then one starts to add the digits of $\omega$ that are to the left of that $\beta$ until one obtains a word $w$ with $\sizer(w)\geq r-2$, which by minimality must be in $\Sigma^{(r-2)}(3+e^{-r})$); in particular it has a $(u,v)$-weakly renormalizable extension $\hat{w}=\hat{w}_1\hat{\gamma}$ where $\hat{\gamma}\in\langle u,v\rangle$ and $\hat{w}_1$ is a suffix of $uv$. We claim that $|u|\leq|\hat{w}_1|$. Otherwise, we use \Cref{lem:rcomparedwithsize} to obtain that
		\begin{equation*}
			\frac{r-2}{1.8}-1\leq|\hat{w}|\leq|\hat{w}_1\beta|\leq|\alpha\beta|-1\leq\frac{r}{6}-1
		\end{equation*}
		which is a contradiction. Hence $|u|\leq|\hat{w}_1|$ and $\hat{w}_1$ ends with $u^b$. Therefore there must be a $u^b\alpha^{f}$ with $f\geq 0$ before the first $\beta$.  
		
		Now we want to extend the word to the right.  Consider now the continuation $w\in\Sigma^{(r-2)}(3+e^{-r})$ that begins at $\alpha\beta\alpha^{e_\ell}$. In particular it has an extension $\hat{w}$ that is $(\alpha,\beta)$-weakly renormalizable by \Cref{lem:local_extension}. Since  
		\begin{align*}
			\sizer(\alpha\beta\alpha^{e_{\ell-1}+1}\alpha\beta)&\leq\sizer(\alpha\beta\alpha)+\sizer(\alpha^{e_{\ell-1}})+\sizer(\alpha\beta)+4 \\
			&\leq 3.6|\alpha\beta|+1.8|\alpha|+r-6|\alpha\beta|+8\leq r-2
		\end{align*} 
		we deduce that $\hat{w}$ contains all $\alpha\beta\alpha^{e_{\ell-1}+2}$ if $e_\ell\geq e_{\ell-1}+2$ and after it must come a $\alpha^g\beta$ or $\alpha^g\hat{w}_2$ where $\hat{w}_2$ starts with $v_a$ and $g\geq 0$. In conclusion
		\[
		u^b\alpha^{f}\beta\alpha^{e_{\ell-1}}\beta\alpha^{e_{\ell-1}+2+g}v_a
		\]
		is a subword of $\omega\in\Sigma(3+e^{-r})$. In this situation the first case of \Cref{lem:symmetry2} yields
		\[
		u^b\alpha^f\beta\alpha^{e_{\ell-1}}\beta\alpha^{e_{\ell-1}+2+g}v_a=(\alpha^{e_{\ell-1}+2+g}\beta\alpha^{e_{\ell-1}}\beta\alpha^f\beta)^\ast 
		\]
		so we still get that $\tilde{\alpha}\tilde{\alpha}\tilde{\beta}\tilde{\beta}=\alpha^{e_{\ell-1}+2}\beta\alpha^{e_{\ell-1}}\beta$ is a subword of $\omega^*\in\Sigma(3+e^{-r})$, contrary to \Cref{lem:aa_bbr} above. This finishes the case $\theta=\alpha$.
		
		Now assume $\gamma=\beta^{e_0}\alpha\beta^{e_1}\alpha\ldots$. Take the maximal integer $k\leq e_{i+1}$ satisfying $\sizer(\beta^{e_i}\alpha\beta^k)\leq 2r$. If $e_i<k$, then letting $(\tilde{\alpha},\tilde{\beta})=(\alpha\beta^{e_i},\beta)$ one gets that $\alpha\beta^{e_i}\alpha\beta^k=\tilde{\alpha}\tilde{\alpha}\tilde{\beta}^{k-e_i}$ is a subword of a word $\omega\in\Sigma(3+e^{-r})$. Observe that
		\begin{align*}
			\sizer(\alpha\beta^{e_i}\alpha\beta^{e_i+2})&\leq \sizer(\alpha)+\sizer(\beta^{e_i})+\sizer(\alpha\beta^2)+\sizer(\beta^{e_i})+6 \\
			&\leq 2r-4|\alpha\beta|+\sizer(\alpha)+\sizer(\alpha\beta^2)+6\leq 2r,
		\end{align*}
		where we used $\sizer(\alpha)+\sizer(\alpha\beta^2)\leq 3.6|\alpha\beta|+3.6$ for $|\alpha\beta|\geq 24$ and the explicit values for $|\alpha\beta|< 24$. Since $k\leq e_{i+1}$, if $e_{i+1}-e_i\geq 2$, we have a contradiction again with \Cref{lem:aa_bbr}. If $e_i\geq k$, then $\sizer(\beta^{k})\leq r-2|\alpha\beta|$ so $k=e_{i+1}$ as before, hence let $(u,v)=(\alpha\beta^{e_{i+1}-1},\beta)$ and $(\tilde{\alpha},\tilde{\beta})=(uv,v)$ so by the second case of \Cref{lem:symmetry2} one gets that
		\[
		\alpha\beta^{e_i}\alpha\beta^{e_{i+1}}\alpha\beta=\tilde{\alpha}\tilde{\beta}^{e_i-e_{i+1}}\tilde{\alpha}\alpha\beta=v_au^b\tilde{\beta}^{e_i-e_{i+1}}\tilde{\alpha}v_a\alpha^b=\beta_a(\tilde{\alpha}\tilde{\alpha}\tilde{\beta}^{e_{i}-e_{i+1}})^\ast\alpha^b
		\]
		is a subword of $\gamma$. If $e_i-e_{i+1}\geq 2$ then one would get that $\tilde{\alpha}\tilde{\alpha}\tilde{\beta}\tilde{\beta}=\alpha\beta^{e_{i+1}}\alpha\beta^{e_{i+1}+2}$ is a subword of $\gamma^\ast$ with $\sizer(\tilde{\alpha}\tilde{\alpha}\tilde{\beta}\tilde{\beta})=\sizer(\alpha\beta^{e_{i+1}}\alpha\beta^{e_{i+1}+2})\leq 2r$, which is impossible. This finishes the first bullet point for $\theta=\beta$. In the particular case where $i=0$, we do not have $\alpha$ before $\beta^{e_0}$. In the case where $|\alpha|\leq|\beta|$ this is no problem because then $(\alpha^b)^\ast\tilde{\alpha}\tilde{\alpha}\tilde{\beta}^{e_{i}-e_{i+1}}(\beta_a)^\ast$ is a subword of $\gamma$ after removing a $\alpha^\ast$ at the end, so we  still get that $\tilde{\alpha}\tilde{\alpha}\tilde{\beta}\tilde{\beta}$ is a subword of $\gamma$. When $(\alpha,\beta)=(uv,v)$, an analogous argument as before shows that the word $\beta^{e_0}\alpha\beta$ has a continuation to the left that is $(\alpha,\beta)$-weakly renormalizable. So before there is either a $\alpha\beta^f$ or a $\hat{w}_1\beta^f$, where $\hat{w}_1$ is a suffix of $\alpha\beta$ that ends with $u^b$. Similarly there is a $\beta^gv_a$ to the right of $\beta^{e_1}\alpha$. In resume the word $u^b\beta^{e_0+f}\alpha\beta^{e_1+g}v_a$ is a subword of $\omega\in\Sigma(3+e^{-r})$. But the second case of \Cref{lem:symmetry2} implies that 
		\[
		u^b\beta^{e_0+f}\alpha\beta^{e_1}\alpha \beta^{g}v_a=(\alpha\alpha\beta^{g}\alpha\beta^{e_1}\alpha\beta^{e_0+f})^*
		\]
		So again if $e_0-e_1\geq 2$, then $\tilde{\alpha}\tilde{\alpha}\tilde{\beta}\tilde{\beta}=\alpha\beta^{e_1}\alpha\beta^{e_1+2}$ is a subword of $\omega^*\in\Sigma(3+e^{-r})$, which contradicts once more \Cref{lem:aa_bbr}.
		
	\end{proof}
	
	\subsection{Constructing renormalizable extensions}
	We start by considering  local extensions. More precisely, if we have a word $w$ that starts (ends) with $\alpha\beta$ where $(\alpha,\beta)\in\overline{P}$, then the alphabet is uniquely determined, and the beginning (end) of $w$ should be $(\alpha,\beta)$-weakly renormalizable. This is a consequence of \Cref{lem:local_extension}.

	Now we want to consider extensions of renormalizable words. The next lemma says that if we have a power $(uv)^s$ with $(u,v)\in\overline{P}$, then it will have a large extension $w$ that is ``almost'' $(u,v)$-weakly renormalizable consisting mostly of powers $(uv)^{s_i}$. This is explained since exponents can only decrease linearly, in fact, they may only decrease by 1 when below some threshold by \Cref{lem:bombieri}. We say that $w=\gamma w_2$ is almost $(u,v)$-weakly renormalizable, because the tail $w_2$ satisfies now the condition $|w_2|<2|uv|$ and $w_2$ is a prefix of a word in $\{u u v, u v v \}$. 
	\begin{lemma}\label{lem:uvextension}
		Let $w \in \Sigma^{(Tr-2)}(3 + e^{-r})$ be a finite word starting with $\theta^s$, where $\theta = uv$ and $(u, v) \in \overline{P}$, $\sizer(\theta^s)\leq r-4|\theta|$, $|uv|<r/9$ and also
		\[
		Tr\leq s^2|\theta|\log((3+\sqrt{5})/2)/2.
		\]
		Then, $w = \gamma w_2$, where $\gamma = (uv)^{s_1} \theta_1 (uv)^{s_2} \theta_2 \ldots (uv)^{s_\ell}$, $s_j\geq 0$, each $\theta_j$ belongs to $\{u u v, u v v \}$ and $w_2$ is a prefix of a word in $\{u u v, u v v \}$. Moreover, $\ell \leq 2.1Tr/|\theta^s|+1$.
	\end{lemma}
	\begin{proof}
		Let $\gamma$ be the largest prefix of the word $w$ than can be written in the form $\gamma = (uv)^{s_1} \theta_1 (uv)^{s_2} \theta_2 \ldots (uv)^{s_\ell}$ where each $\theta_j$ belongs to $\{u u v, u v v \}$ and $s_j\geq 0$. Then, we claim that:
		\begin{itemize}
			\item If $\sizer((uv)^{s_{j+1}})\leq r-6|uv|$, then $\theta_j=\theta_{j+1}$, for $1\leq j\leq \ell-2$.
			\item If $\sizer((uv)^{s_j})\leq r-10|uv|$, then $|s_j-s_{j+1}|\leq 1$ for $1\leq j\leq \ell-2$.
			\item $s_j\geq s-j$, for all $1\leq j\leq \ell-1$.
		\end{itemize} 
		
		Indeed, to see the first claim, if $\theta_j\neq\theta_{j+1}$ note that
		\begin{align*}
			\sizer(\theta_j(uv)^{s_{j+1}}\theta_{j+1})&\leq \sizer(\theta_j)+\sizer((uv)^{s_{j+1}})+\sizer(\theta_{j+1})+4\\
			&\leq 1.8|uuv|+r-6|uv|+1.8|uvv|+7.6\leq r-2,
		\end{align*}
		which is clear for $|uv|\geq 16$, while for $|uv|<16$ we computed the explicit values of $\sizer(uuv)+\sizer(uvv)$ to check that the inequality $\sizer(\theta_j(uv)^{s_{j+1}}\theta_{j+1})\leq r-2$ still holds. If $\theta_j=uuv$, let $(\alpha,\beta)=(u,uv)$ and, if $\theta_{j}=uvv$, let $(\alpha,\beta)=(uv,v)$. So $\alpha\beta\tau'=\theta_j(uv)^{s_{j+1}}\theta_{j+1}$ is a subword of a word in $\Sigma^{(r-2)}(3+e^{-r})$. Then \Cref{lem:local_extension} says that this word is contained in a word that can be written in the alphabet $\{\alpha,\beta\}$. But, if $\theta_{j+1}\neq\theta_j$, this would be impossible.
		
		To prove the second claim, consider the subword $\theta_{j-1}(uv)^{s_j}\theta_{j}(uv)^{s_{j+1}}$. Since $\theta_{j-1}=\theta_{j}$, using the appropriate pair $(\alpha,\beta)\in\overline{P}$ this whole word can be written in that alphabet. An application of Bombieri's characterization (\Cref{lem:bombieri}) yields $\theta_j=\theta_{j+1}$ and also the second claim. Indeed, if $(\alpha,\beta)=(u,uv)$ then $\theta_{j-1}(uv)^{s_j}\theta_{j}(uv)^{s_{j+1}}u=\alpha\beta^{s_j+1}\alpha\beta^{s_{j+1}+1}\alpha$ is a subword of a word in $\Sigma(3+e^{-r})$ so the first bullet point of \Cref{lem:bombieri} gives $|s_j-s_{j+1}|\leq 1$. If $(\alpha,\beta)=(uv,v)$ then $\theta_{j-1}(uv)^{s_j}\theta_{j}(uv)^{s_{j+1}}=\alpha\beta\alpha^{s_j+1}\beta\alpha^{s_{j+1}}$, so the third bullet point of \Cref{lem:bombieri} gives that $s_{j+1}\leq s_j+2$. Therefore,
		\begin{align*}
			\sizer((uv)^{s_{j+1}})&\leq \sizer((uv)^{s_{j}})+\sizer((uv)^2)+2\\
			&\leq r-10|uv|+1.8|uv|+3.8\leq r-6|uv|,
		\end{align*}
		which is clear for $|uv|\geq 10$, while for $|uv|<10$ we checked the inequality $\sizer((uv)^2)+2\leq 4|uv|$ directly. The previous item gives $\theta_j=\theta_{j+1}$, hence $\theta_{j-1}(uv)^{s_j}\theta_{j}(uv)^{s_{j+1}}\theta_{j+1}=\alpha\beta\alpha^{s_j+1}\beta\alpha^{s_{j+1}+1}\beta$ and finally we use the first bullet point of \Cref{lem:bombieri}.
		
		We now prove the third claim. Observe that $s_1\geq s-1$ by construction of $\gamma$. If $s_j\geq s-2$ then we are done. Note that, if $s_j\leq s-3$, then 
		\begin{equation*}
			\sizer((uv)^{s_j})\leq r-4|uv|-\sizer((uv)^3)+1\leq r-6|uv|
		\end{equation*}
		where we used $\sizer((uv)^3)\geq 2.8|uv|-3$ for $|uv|>4$ and for $(u,v)=(a,b)$ we used $\sizer((uv)^3)=16$. Hence the first claim gives that $s_j\leq s-3$ implies $\theta_j=\theta_{j-1}$. In particular $(uv)^{s_{j-1}}\theta_{j-1}(uv)^{s_j}\theta_{j}$ can be written in the appropriate alphabet $(\alpha,\beta)\in\overline{P}$ as $\alpha^{s_{j-1}+1}\beta\alpha^{s_j+1}\beta$ or $\beta^{s_{j-1}}\alpha\beta^{s_j+1}\alpha\beta$. Since $s_j\leq s-3$ implies $\sizer((uv)^{s_j})\leq r-2|\alpha\beta|$ by hypothesis, then the \Cref{lem:bombieri} gives $s_j\geq s_{j-1}-1$. This proves the three claims.
		
		Now, note that $s\geq 2$, since, if $s=1$, we get a contradiction because
		\[
		\sizer(\theta)\leq \sizer(w)\leq Tr\leq |\theta|/2
		\]
		gives $|\theta|\leq 8$, but $\sizer(ab)=5$ and $\sizer(ab)\leq \sizer(\theta)\leq |\theta|/2\leq 4$.

		On the other hand using \Cref{lem:rcomparedwithsize} one has that
		\begin{align*}
			\left(\log\left(\frac{3+\sqrt{5}}{2}\right)\right)^{-1}Tr+4&\geq |w| \geq |\gamma|=(s_1+\dotsb+s_\ell)|uv|+|\theta_1|+\dotsb+|\theta_\ell| \\
			&\geq (s_1+\dotsb+s_\ell+\ell)|\theta| 
		\end{align*}
		If $\ell\geq s$ then one gets that
		\begin{align*}
			\left(\log\left(\frac{3+\sqrt{5}}{2}\right)\right)^{-1}Tr+4&> (s+(s-1)+\dots+1)|\theta| \\
			&\geq \frac{s^2|\theta|}{2}+\frac{s|\theta|}{2} \geq \frac{s^2|\theta|}{2}+4,
		\end{align*}
		
		contradicting the hypothesis (we used irrationality). Thus $\ell<s$. Hence 
		\[
		\left(\log\left(\frac{3+\sqrt{5}}{2}\right)\right)^{-1}Tr+4\geq (\ell-1)(s-\ell/2)|\theta|
		\]
		This shows that $\ell\leq 2.1Tr/|\theta^s|+1$ (for $\ell=2$ we use instead $\sizer(\theta^s)\leq Tr$).

		Now, write $w=\gamma w_2$. We have two cases to consider. When $s_\ell\leq 1$, we choose $\hat{w}\in\Sigma^{(r-2)}(3+e^{-r})$ starting with $\theta_{\ell-1}$. On the other hand, when $s_\ell\geq 2$, we choose $\hat{w}\in\Sigma^{(r-2)}(3+e^{-r})$ starting at the last occurrence of $uvuv$. Since $|uv|<r/6$, \Cref{lem:local_extension} gives that it is $(u,v)$-semi renormalizable. We know that $\hat{w}$ is $(\hat{\alpha},\hat{\beta})$-semi renormalizable for some $(\hat{\alpha},\hat{\beta})\in\overline{P}$ with $|\hat{\alpha}\hat{\beta}|\geq r/6$, because of \Cref{cor:r/6}. Since $|uv|< r/6$ and $\hat{w}$ has $uv$, the word $\hat{w}$ is $(\alpha,\beta)$-semi renormalizable for some $(\alpha,\beta)\in\{(uv,v),(u,uv)\}$. Possibly adding one digit to the right, write $\hat{w}=\hat{\gamma}\hat{w_2}$ with $\hat{\gamma}\in\langle u,v\rangle$ and $\hat{w_2}$ a prefix of $uv$. Observe that $5|uv|<r/1.8<|\hat{w}|$, thus $|\hat{\gamma}|=|\hat{w}|-|\hat{w_2}|>4|uv|$.
		
		We will find a continuation of $\gamma\gamma'$ of $\gamma$ of the form
		\[
		\gamma'\in\{uv,uuv,uvv\}\cap\{\alpha,\beta,\alpha\beta\}.
		\]
		By the maximality of $\gamma$, we have that $w_2$ is contained in $\gamma'$ which will finish the proof of the lemma. When $s_\ell\geq 2$, then if $uvuv=\alpha\alpha$ we have that $\alpha\alpha\beta=uv(uvv)$ extends and $\alpha\alpha\alpha=(uv)^3$ extends $(uv)^{s_\ell}$ to $(uv)^{s_\ell+1}$. If $uvuv=\beta\beta$ then $\beta\beta\beta=(uv)^{3}$ extends $(uv)^{s_\ell}$ to $(uv)^{s_\ell+1}$ while if there is $\alpha\alpha$ after $\beta\beta$, then it must come $\alpha^kv_a$ (if there is no $\beta$ after $\alpha\alpha$ then $\hat{w_2}$ starts with $v_a$). But note that
		\[
		\beta\beta\alpha^k v_a=v_au^b\beta\alpha^kv_a=v_a(\alpha^k\beta\beta)^\ast
		\]
		so $\alpha\alpha\beta\beta$ is a subword of $w^\ast$, a contradiction with \Cref{lem:aa_bbr}. In the situation where $s_\ell= 1$, we have that $uvvuv=\alpha\beta\alpha$, so if there is an $\alpha$ afterwards, then $\alpha\beta\alpha\alpha$ extends $(uv)^{s_\ell}$ while $\alpha\beta\alpha\beta$ gives a contradiction because before $\theta_\ell$ there is $s_{\ell-1}\geq 2$ and we have $\alpha\alpha\alpha\beta\alpha\beta=\alpha\alpha\tilde{\beta}\tilde{\beta}$. When $\theta_{\ell-1}uv=uuvuv=\alpha\beta\beta$, if there is $\beta$ afterwards then $\alpha\beta\beta\beta=\theta_{\ell-1}(uv)^{2}$ extends while $\alpha\beta\beta\alpha\beta=\theta_{\ell-1}(uv)\theta_{\ell}$ also extends, but if there is $\alpha\alpha$ after $\beta$ we have the same contradiction with \Cref{lem:aa_bbr} in the transpose word. In the situation $s_\ell=0$, we have $uvuv\theta_{\ell-1}\in\{\alpha\alpha\alpha\beta,\beta\beta\alpha\beta\}$ and we arrive at the same contradictions or extensions as before. In summary, $w_2$ is a subword of $\gamma'$, so it is a subword of a word in $\{u v, u u v, u v v \}$.
		
	\end{proof}

	In the case where we have a power $a^e$ or $b^e$, then its extensions are not necessarily $(a,b)$-weakly renormalizable, because we could have odd powers of a digit $\{1,2\}$ appearing after. Neverthless, it takes a long time for these powers to decay, which gives us the next extension lemma.
	
	\begin{lemma}\label{lem:digitextension}
		Let $w \in \Sigma^{(Tr-2)}(3 + e^{-r})$ be a finite word starting with $c^s$, where $c \in \{1,2\}$ and $\sizer(c^s)\leq r-2$. Suppose that $Tr\leq s^2(\log x)/4$ where $x=(3+\sqrt{5})/2$ and $x=3+2\sqrt{2}$ for $c=1$ and $c=2$ respectively. Then $w = c^{s_1} \theta' c^{s_2} \theta' \ldots \theta' c^{s_\ell}$, $\theta' = c'c'$,  $c' \in \{1, 2\}$, $c \neq c'$ and $Tr\geq \ell(s-\ell+1)\log x$ (in particular $\ell<2Tr/(s\log x)$). Moreover if $\sizer(c^{s_j})\leq r-8$ then $s_j$ is even and $|s_{j+1}-s_j|\in \{0,2\}$.
	\end{lemma}
	\begin{proof}
		Observe that we only have to prove that in $w$ there is no $c' c' c' c'$, since the sequence $cc' c' c' c$ is forbidden by \Cref{lem:1}. Define $\gamma=c^{s_1} \theta' c^{s_2} \theta' \ldots \theta' c^{s_\ell}$ to be the longest sequence inside $w$ starting with $c^s$ that does not contain $(\theta')^2$. The fact that $\sizer(c^{s_j})\leq r-8$ implies $s_j$ even and $|s_{j+1}-s_j|\in \{0,2\}$, is because of \Cref{lem:1,lem:bombieri}. The inequality \eqref{eq:intervalswith2s} for $c=2$ and the inequality \eqref{eq:intervalswith1s} for $c=1$ gives us $Tr\geq \ell(s-\ell+1)\log x$. Since $\ell\leq s/2-1$ (by hypothesis), we can bound the length  $\ell<2Tr/(s\log x)$. If $w\neq\gamma$, then $\gamma$ must end with $\theta' cc\theta'$, but this would imply that $\ell\geq s/2$.
		
	\end{proof}
	
	The next situation is where we have an $(\alpha,\beta)$-renormalizable word $w$ that does not contain big powers of $\alpha$ or $\beta$. In this case there is an extension that can be written in the same alphabet $(\alpha,\beta)$. In fact, this extension is almost $(\alpha,\beta)$-weakly renormalizable, in the sense that its tail is small; it is a prefix of a word in $\langle\alpha,\beta\rangle$ but is not necessarily a prefix of $\alpha\beta$.
	\begin{lemma}\label{lem:alphabetaextension}
		Let $w\in\Sigma^{(r-2)}(3+e^{-r})$ be an $(\alpha,\beta)$-weakly renormalizable word with $|\alpha\beta|<r/40$. Suppose that for every factor of the form $\alpha^s$ or $\beta^s$ of $w$, we have $|\alpha^s|<\delta|w|$ and $|\beta^s|<\delta|w|$ where $0<\delta<(1/2)\log((3+\sqrt{5})/2)$ is a constant. Assume further that $w$ contains an $\alpha\beta$. If $\overline{w}\in\Sigma^{(Tr)}(3+e^{-r})$ is an extension of $w$ with $T+2\leq\delta' r/(16|\alpha\beta|)$ where $\delta'=1-2\delta/\log((3+\sqrt{5})/2))$, then $w\overline{w}= w_1\gamma\tau'$ where $\gamma\in\langle\alpha,\beta\rangle$, $|\tau'|<|\alpha\beta|$ is a prefix of some word in $\langle\alpha,\beta\rangle$ and $w_1$ is a suffix of $\alpha\beta$.
	\end{lemma}
	
	\begin{proof}
		Write $w=w_1\gamma w_2$ as in the definition of $(\alpha,\beta)$-weakly renormalizable. Take $\gamma'$ to be the largest word inside $w\overline{w}$ starting with $\gamma$ which can be written in the alphabet $\langle\alpha,\beta\rangle$. Write $w\overline{w}=w_1\gamma'\tau'$. We will show that $\tau'$ is a prefix of some word in $\langle\alpha,\beta\rangle$ and $|\tau'|<|\alpha\beta|$. Take the last factor of $\alpha\beta$ in $\gamma'$, say $\gamma'=\eta\alpha\beta\eta'$. In particular by \Cref{lem:determine_alphabet} we obtain $\eta'=\theta^s$ with $\theta\in\{\alpha,\beta\}$. Consider the factor $\hat{w}\in\Sigma^{(r-2)}(3+e^{-r})$ (that possibly extends $w\overline{w}$) starting at this last occurrence of $\alpha\beta$. By \Cref{lem:local_extension}, after possibly adding one digit to the right, it can be written as $\hat{w}=\hat{\gamma}\hat{w_2}$, $\hat{\gamma}=\alpha\beta\ldots$ with $\hat{\gamma}\in\langle\alpha,\beta\rangle$ and $\hat{w_2}$ a prefix of $\alpha\beta$. If $\alpha\beta\theta^s$ is strictly contained on $\hat{\gamma}$, then if $|\tau'|\geq|\alpha\beta|$ we will get that there is at least one more letter $\{\alpha,\beta\}$ of $\hat{\gamma}$ inside $w\overline{w}$ after $\gamma'$, a contradiction with the maximality of $\gamma'$. Hence $|\tau'|<|\alpha\beta|$ and $\alpha\beta\theta^s\tau'$ is a subword of $\hat{\gamma}\alpha\beta$, so  $\tau'$ is the prefix of a word in $\langle\alpha,\beta\rangle$ as claimed. Now suppose that $\alpha\beta\theta^s$ contains $\hat{\gamma}$, so $\sizer(\theta^s)\geq r-4|\alpha\beta|-6$. We will get a contradiction by considering an extension in the transpose word $(w\overline{w})^\ast$.
		
		We need the identity $((uv)^k)^\ast=u^b(uv)^{k-1}v_a$. For $k=1$ this is consequence of \Cref{lem:structure}. To prove the identity for $k\geq 2$, let $(\alpha,\beta)=(uv,v)$, so $(uv)^k=v_au^b(uv)^{k-2}v_au^b=v_au^b\alpha^{k-2}\beta_au^b$. Now using that $u^b,v_a,\alpha^k\beta_a$ are all palindromic, one gets that $((uv)^k)^\ast=u^b\alpha^{k-2}\beta_au^bv_a=u^b(uv)^{k-1}v_a$.
		
		In resume there is a $\theta^{s-1}$ inside $(w\overline{w})^\ast$. Let $\theta^{s'}$ be the maximal suffix of $\theta^{s-1}$ that satisfies $\sizer(\theta^{s'})\leq r-4|\theta|$. Since $\sizer(\theta^{s'})\geq r-6|\alpha\beta|-10$, by \Cref{lem:rcomparedwithsize} we have that $s'|\theta|=|\theta^{s'}|>r/2$. Observe that $\sizer((w\overline{w})^\ast)\leq (T+2)r-4$ by \Cref{prop:rtranspose}, which gives 
		\[
		(s')^2|\theta|/4 \geq (r/8)\cdot\frac{r}{2|\theta|} \geq(T+2)r
		\]
		since $T+2\leq \frac{r}{16|\theta|}$ by hypothesis. Note that $(T+2)r\leq (s')^2/8$ also holds.

		If $\theta=uv$ with $(u,v)\in\overline{P}$ we use \Cref{lem:uvextension} to find a $\tilde{\gamma}=\theta^{s_1}\theta_1\theta^{s_2}\ldots \theta^{s_\ell}$ starting with this $\theta^{s'}=(uv)^{s'}$, where each $\theta_i\in\{uuv,uvv\}$ and such that $(w\overline{w})^\ast=\Tilde{\eta}\Tilde{\gamma}\Tilde{w}_2$ with $\Tilde{w}_2$ a prefix of some word in $\{uuv,uvv\}$. When $|\theta|=2$, we write $\theta=cc$ with $c\in\{1,2\}$, $\theta'\in\{a,b\}$ and use \Cref{lem:digitextension} to find an extension $\tilde{\gamma}=c^{e_1}\theta' c^{e_2}\theta'\ldots c^{e_\ell}$ starting with $\theta^{s^\prime}$ such that $(w\overline{w})^\ast=\Tilde{\eta}\Tilde{\gamma}$.

		Since $(w\overline{w})^\ast=\overline{w}^\ast w^\ast$, when $\theta=uv$ we have that $\theta^{s_\ell}$ (or at least some factor of it) is inside of $w^\ast$. Similarly, when $|\theta|=2$, we have that $c^{e_{\ell}}$ (or at least some factor of it) is inside $w^\ast$, so $\theta^{s_\ell}$ is inside $w$ where $s_\ell=2\lfloor e_\ell/2\rfloor$.
		
		In any case $\theta^{s_\ell-1}$ is inside of $w$ (because of the identity proved above). In any case we will also get that $s_\ell\geq s'-1-\ell$ and also $\ell\leq 2.1Tr/|\theta^{s'}|+1$. But note that using $|\theta^{s'}|>r/2$ and \Cref{lem:rcomparedwithsize}
		\begin{align*}
			|\theta^{s_\ell-1}|&\geq|\theta^{s'-1-\ell}| = |\theta^{s'}|-|\theta^{\ell+1}| \geq r/2-(\ell+1)|\theta| \\
			&\geq r/2-\left(2.1Tr/|\theta^{s'}|+2\right)|\theta| \geq r/2 - 2(2.1T+1)  \\
			&\geq r/2 - \frac{|\theta|}{|\alpha\beta|}\cdot\delta' r+6|\theta|\geq (1/2)(1-\delta')r+6|\theta| \geq \delta|w|
		\end{align*}
		which is a contradiction with the existence of those factors inside $w$. In conclusion  $w\overline{w}= w_1\gamma\tau'$ where $|\tau'|<|\alpha\beta|$ is a prefix of some word in $\langle\alpha,\beta\rangle$ and $w_1$ is a suffix of $\alpha\beta$.
	\end{proof}

	\subsection{Proof of \texorpdfstring{\Cref{thm:mainstrong}}{Theorem 1.2}}
	This section will be devoted to the proof of \Cref{thm:mainstrong}. The main idea is to find a subcovering of the natural covering of $K_t$. Indeed, recall from the introduction that 
	\[
	d(t)=\min\{1,2\cdot\dimH(K_t)\}=\min\{1,2\cdot\dim_{\mathrm{B}}(K_t)\}.
	\]
	In order to prove that the Hausdorff dimension of $K_t$ is at most $d$ we will start with a covering of $K_t$ by a finite union of intervals and then replace each of these intervals by a suitable union of smaller subintervals such that the sum of the $d$-th powers of the sizes of the subintervals is smaller than the $d$-th power of the size of the initial interval. 
	
	The proof is quite long, so it is divided into several subsections. Moreover, we will need the following combinatorial lemma.

	\begin{lemma}\label{lem_Comb}
		Let $U$ be a positive integer, and let $m$ be a positive real number. If $U\le m$ then the number of solutions $(\ell,x_1,x_2,\dots, x_{\ell})$ of 
		$x_1+x_2+\dots+x_{\ell}\le (U-\ell)m$ with each $x_i \in \NN^*$ is at most
		\[
		U\left(\frac{em\varepsilon_m}{1-\varepsilon_m}\right)^{(1-\varepsilon_m)(U+1)}=Ue^{(1-\varepsilon_m)(U+1)/\varepsilon_m}=Ue^{W(m)(U+1)},
		\]
		where $\varepsilon_m$ is the solution in $(0,1)$ of the equation $\log\left(\frac{em\varepsilon}{1-\varepsilon}\right)=\frac1{\varepsilon}$.
		
		In general (including the case when $U>m$), this number of solutions is at most $Ue^{Um/e^{W(m-1)}}$. For $m>1$, this upper bound is equal to $Ue^{U\cdot\frac{m}{m-1}W(m-1)}$, and for $m\ge 5$, this is at most 
		$Ue^{U\cdot \frac{\log m}{m}}e^{U\cdot W(m-1)}<Ue^{U\cdot\frac{\log m}{m}}e^{U\cdot W(m)}$.
		
		In particular, if $U=\ao\left(\frac{\log m}m\right)$, then this number is at most $e^{(\log m-\log\log m+\ao(1))(U+1)}$.
		
	\end{lemma}
	
	\begin{proof}
		We should have $1\le \ell \le U$. Given such $\ell$, the number of solutions of this inequality is the number of natural solutions of $x_0+x_1+\dots+x_{\ell}=\lfloor (U-\ell)m \rfloor$, where $x_0$ is included to transform the inequality into an equality. This is equal to $\binom{\lfloor (U-\ell)m \rfloor+\ell}{\ell}$, and using the inequalities $\binom{n}{k}\le \frac{n^k}{k!} \le \left(\frac{en}k\right)^k$, which hold for $1\le k\le n$, this number of solutions is at most
		\[
		\left(\frac{e(\lfloor (U-\ell)m \rfloor+\ell)}{\ell}\right)^{\ell}\le \left(\frac{e((U-\ell)m+\ell)}{\ell}\right)^{\ell}.
		\]
		If $U\le m$, then $\ell\le U\le m$ and $(U-\ell)m+\ell\le (U+1-\ell)m$, so the previous upper estimate is at most 
		\[
		\left(\frac{e((\tilde U-\ell)m)}{\ell}\right)^{\ell},
		\]
		where $\tilde U\colonequals U+1$. Let $\varepsilon\in (0,1)$ such that $\ell=(1-\varepsilon)\tilde U$, so $\frac{\tilde U-\ell}{\ell}=\frac{\varepsilon}{1-\varepsilon}$. The derivative of $g(\ell)=\log \left(\frac{e((\tilde U-\ell)m)}{\ell}\right)^{\ell}=\ell \log \left(\frac{e((\tilde U-\ell)m)}{\ell}\right)$ is
		\[
		\log\left(\frac{e((\tilde U-\ell)m)}{\ell}\right)-\frac{\tilde U}{\tilde U-\ell}=\log\left(\frac{em\varepsilon}{1-\varepsilon}\right)-\frac1{\varepsilon},
		\]
		and so $g(\ell)$ is maximized for $\ell=(1-\varepsilon_m)\tilde U$. Moreover, since there are $U$ possible values of $\ell$, the number of solutions we are estimating is at most
		\[
		U\cdot e^{g((1-\varepsilon_m)\tilde U)}=Ue^{(1-\varepsilon_m)(U+1)/\varepsilon_m},
		\]
		since, by definition, $\frac{em\varepsilon_m}{1-\varepsilon_m}=e^{1/\varepsilon_m}$.
		
		Notice that, since $\varepsilon_m$ is the solution in $(0,1)$ of the equation $\log\left(\frac{em\varepsilon}{1-\varepsilon}\right)=\frac1{\varepsilon}$, writing $x_m=\frac{1-\varepsilon_m}{\varepsilon_m}$, we have $\frac1{\varepsilon_m}=x_m+1$, and thus $\log\left(\frac{em}{x_m}\right)=x_m+1$, and $\log\left(\frac{m}{x_m}\right)=x_m$. It follows that $x_m e^{x_m}=m$, and thus $x_m=W(m)$, and $Ue^{(1-\varepsilon_m)(U+1)/\varepsilon_m}=Ue^{(U+1)x_m}=Ue^{W(m)(U+1)}$.
		
		In the general case, let us estimate
		$$h(\ell)=\left(\frac{e((U-\ell)m+\ell)}{\ell}\right)^{\ell}.$$
		The derivative of $\log h(\ell)=\log \left(\frac{e((U-\ell)m+\ell)}{\ell}\right)^{\ell}=\ell \log \left(\frac{e((U-\ell)m+\ell)}{\ell}\right)$ is
		\[
		\log \left(\frac{(U-\ell)m+\ell}{\ell}\right)-\frac{\ell(m-1)}{(U-\ell)m+\ell}=\log z-\frac{m-1}z,
		\]
		where $z=\frac{(U-\ell)m+\ell}{\ell}$ is a decreasing function of $\ell$ and so $h(\ell)$ is maximized when $\log z=\frac{m-1}z$, which is equivalent to $z\log z=m-1$ and to $\log z=W(m-1)$. In this case, we have $h(\ell)=(ez)^{\ell}=e^{\ell\log(ez)}$. Since $\ell=\frac{Um}{z+m-1}$ and $\log(ez)=1+\log z=1+\frac{m-1}z=\frac{z+m-1}z$, we have $\ell\log(ez)=\frac{Um}z=\frac{Um}{e^{W(m-1)}}$, which gives our upper estimate for $h(\ell)$:
		$$e^{Um/e^{W(m-1)}},$$
		and as before, since there are $U$ possible values of $\ell$, the number of solutions we are estimating is at most
		$$Ue^{Um/e^{W(m-1)}}.$$
		For $m>1$, we have $e^{W(m-1)}=\frac{m-1}{W(m-1)}$, so our estimate becomes 
		$$Ue^{UmW(m-1)/(m-1)}=Ue^{U\cdot\frac{m}{m-1}W(m-1)}.$$
		We have $\frac{m}{m-1}W(m-1)=W(m-1)+\frac{W(m-1)}{m-1}$, so 
		$$Ue^{U\cdot\frac{m}{m-1}W(m-1)}=Ue^{U\frac{W(m-1)}{m-1}}e^{UW(m-1)},$$
		and for $m\ge 5$, we have $\frac{W(m-1)}{m-1}<\frac{\log m}m$ (indeed, this is equivalent to
		\[
		\frac{(m-1)\log m}{m}m^{(m-1)/m}=\frac{(m-1)\log m}{m}e^{\frac{(m-1)\log m}{m}}>m-1,
		\]
		which is equivalent to $\log m>m^{1/m}$, and thus holds for every $m\ge 5$). Thus, in this case, our upper estimate becomes 
		$$Ue^{U\cdot \frac{\log m}{m}}e^{U\cdot W(m-1)}<Ue^{U\cdot\frac{\log m}{m}}e^{U\cdot W(m)}.$$
		
		If $U=\ao(\frac{m}{\log m})$ then $U=\ao(e^{\log m-\log\log m})$, and thus, since $\frac{m}{m-1}W(m-1)=W(m-1)+\ao(1)=W(m)+\ao(1)=\log m-\log\log m+\ao(1)$, we have 
		\begin{align*}
			e^{(U+1)\cdot\frac{m}{m-1}W(m-1)}&=e^{\log m-\log\log m+\ao(1)}e^{U\cdot\frac{m}{m-1}W(m-1)} \\
			&=(1+\ao(1))\frac{m}{\log m}e^{U\cdot\frac{m}{m-1}W(m-1)}>Ue^{U\cdot\frac{m}{m-1}W(m-1)},
		\end{align*}
		which was our previous upper estimate, and we have 
		\[
		e^{(U+1)\cdot\frac{m}{m-1}W(m-1)}=e^{(\log m-\log\log m+\ao(1))(U+1)},
		\] 
		which concludes the proof.
		
	\end{proof} 
	
	\begin{remark}
		In the proof of \Cref{thm:mainstrong} we only use the case $U\leq m$ of \Cref{lem_Comb}. We use the general case of \Cref{lem_Comb} only on \Cref{sec:error_optimal}. 
	\end{remark}

	\begin{proof}[\emph{\textbf{Proof of \Cref{thm:mainstrong}}}]
		
		We start by recalling some notation. We denote by $\mathsf{m}(\omega) = \sup_{n \in \ZZ} \lambda(\sigma^n(\omega))$ the Markov value of $\omega$, and we have
		\[
		Q_r= \{\alpha= c_1 c_2 \ldots c_n \ \mid\ \sizer(\alpha) \ge r, \sizer(c_1 c_2 \ldots c_{n-1})<r\},
		\]
		that is, $\alpha$ belongs to $Q_r$ if and only if $\sizes(\alpha) < e^{-r-1}$ and $\sizes(\alpha') \geq e^{-r-1}$, where $\alpha'$ is the word obtained by removing the last letter from $w$.
		
		We now recall how the covering of $K_t$ is constructed. We define the sets of words
		\begin{align*}
			C(t,r)&=\{\alpha= c_1 \ldots c_n \in Q_r\ \mid\ K_t\cap I(\alpha)\neq\emptyset\} \\
			&=\{\alpha\in Q_r\ \mid\ \alpha~\text{subword of a word $\omega\in(\NN^\ast)^\ZZ$ with $\mathsf{m}(\omega)\leq t$}\}.
		\end{align*}
		
		Here, $K_t=\{[0;\gamma]\ \mid\ \gamma\in\pi_{+}(\Sigma(t))\}$ where $\pi_{+}\colon \Sigma\to\Sigma^+$ is the projection associated with the decomposition $\Sigma=\Sigma^{-}\times\Sigma^{+}=(\NN^\ast)^{\ZZ_{-}}\times(\NN^\ast)^{\NN}$. That is,
		\[
		\pi_+(\ldots c_{-2} c_{-1} c_0 c_1 c_2 \ldots)=c_0 c_1 c_2 \ldots.
		\]
		Moreover, $\Sigma(t)=\{\omega\in(\NN^\ast)^\ZZ\ \mid\ \mathsf{m}(\omega)\leq t\}$. It is clear that
		\[
		\calM\cap(-\infty,t)\subseteq (\NN^\ast\cap[1,\lfloor t\rfloor])+K_t+K_t.
		\]
		Observe that $K_t$ is covered by all $I(\alpha)$ where $\alpha\in C(t,r)$ for any fixed $r$. 
		
		If $r\leq s$, then the set $C(t,r)$ covers the set $C(t,s)$, in the sense for any interval $I(\alpha)$ with $\alpha=c_1 \ldots c_n \in C(t,s)$ there is $m\leq n$ such that $\tilde{\alpha}= c_1 \ldots c_m \in C(t,r)$ and $I(\tilde{\alpha})\subseteq I(\alpha)$.
		
		Given $d$ depending on $r$, we can prove that the Hausdorff dimension of $K_t$ is at most $d$, by we replacing an interval $I$ (corresponding to a word in $C(t,r)$) with several intervals $I_j$ contained in it, but smaller and of different sizes, each one corresponding a word in $C(t,Tr)$, where $T\in\{10,\lfloor \log^2 r\rfloor, \lfloor r/5\rfloor\}$, whose union still contains the intersection of $K_t$ with $I$ and that satisfy $\sum_j |I_j|^d<|I|^d$. Since this process can be iterated, this shows that the $d$-dimensional Hausdorff measure of $K_t$ is finite for large enough $r$.

		By \Cref{cor:r/6}, if $w\in \Sigma^{(r+2)}(3+e^{-r-4})$, then there is a sequence of alphabets $(\alpha_j,\beta_j)$ such that, for all $0\le j\le m$, $w$ is $(\alpha_j,\beta_j)$-semi renormalizable, with
		\[
		(\alpha_0,\beta_0)=(a,b) \quad \text{ and } \quad (\alpha_{j+1},\beta_{j+1})\in\{(\alpha_j\beta_j,\beta_j), (\alpha_j,\alpha_j\beta_j)\},
		\]
		for each $0 \leq j<m$ and $|\alpha_m\beta_m|\ge r/6$. 
		
		\noindent We consider such a renormalization $(\alpha_t,\beta_t)$ with
		\[
		r/\sqrt{\log r}\le |\alpha_t| + |\beta_t|<2r/\sqrt{\log r}.
		\]
		We will consider words $\tilde w\in\Sigma^{(r)}(3+e^{-r-4})$ such that $w\tilde w\in \Sigma(3+e^{-r-4},|w\tilde w|)$. Then, depending on $w\tilde w$, we will consider continuations $\overline w \in \Sigma^{(Tr)}(3+e^{-r-4})$ for some $T\in\{10,\lfloor \log^2 r\rfloor, \lfloor r/5\rfloor\}$ such that $w\tilde w\overline w\in \Sigma(3+e^{-r-4},|w\tilde w\overline w|)$. 
		
		The strategy is the following: if $w$ (or, more generally, $w\tilde w$) contains a factor $\alpha_t\beta_t$, we may consider the factor $\hat w\in Q_{r+2}$ of $w\tilde w$ starting at this factor $\alpha_t\beta_t$; it should be $(\alpha_t,\beta_t)$-renormalizable. We will attempt to use this argument several times in order to cover the whole word $w\tilde w$ by $(\alpha_t,\beta_t)$-renormalizable words. To determine $\tilde w$, we only need to estimate the number of words in $(\alpha_t,\beta_t)$ after the last factor equal to $\alpha_t\beta_t$ in $w$. For this sake, we consider several cases according to the size of $\alpha_{t}^s$, $\beta_{t}^s$ as a factor of $w$, for some integer $s$.
		
		In Case 1 and Case 2 below, we choose $T = 10$, while in Case 3 we initially choose $T = \lfloor \log^2 r\rfloor$. Furthermore, in all of the following cases, except for Case 3.2.2, we take $d = \frac{\log r - \log \log r}{r}$. In Case 3.2.2, corresponding to when $\alpha_t = 11$ and $\tilde{w}w$ contains a relatively long factor $\alpha_t^s$, we initially choose the estimate $d = \frac{\log r - \log \log r + c_0 + \ao(1)}{r}$, where $c_0 = -\log\log\left(\frac{3+\sqrt{5}}{2}\right) > 0$. This is already enough to obtain the upper bound
		\[
		d(3+e^{-r})\le 2\cdot\frac{\log r-\log\log r+c_0+\ao(1)}r.
		\]
		The only case that produces a ``bad'' estimate is then Case 3.2.2. This estimate can be actually improved by a refined analysis using $T = \lfloor T/5\rfloor$, giving rise to Case 3.2.3. Our final upper bound in \Cref{thm:mainstrong} is derived in this way.
		
		\medbreak
		
		\noindent\textbf{Case 1:} Suppose first that, for every factor of the form $\alpha_t^s$ or $\beta_t^s$ of $w$, we have $|\alpha_t^s|<|w|/3$ and $|\beta_t^s|<|w|/3$. In this case we take $T = 10$.
		
		Then there is a factor $\alpha_t\beta_t$ in the first half of $w$, and until the next appearence of $\alpha_t\beta_t$ (which happens before the end of $w$), we have a factor with total size smaller than $|w|/2$ of the type $\alpha_t\beta_t\alpha_t^{j}\beta_t$ or $\alpha_t\beta_t^{j}\alpha_t\beta_t$ for some positive integer $j$. Suppose we are in the first case, without loss of generality.

		Then, given a continuation $\tilde w \overline w$ of $w$ with $\tilde w\in\Sigma^{(r)}(3+e^{-r-4})$ and $\overline w \in \Sigma^{(10r)}(3+e^{-r-4})$, \Cref{lem:alphabetaextension} gives that $w\tilde{w}\overline{w}=\tau\gamma\tau'$ with $\gamma\in\langle\alpha_t,\beta_t\rangle$, $\tau$ a suffix of $\alpha_t\beta_t$ and $\tau'$ a prefix of some word in $\langle\alpha_t,\beta_t\rangle$ with $|\tau'|<|\alpha_t\beta_t|$. Thus, the continuation of the first factor of the form $\alpha_t\beta_t$ of $w$ in $w \tilde w \overline w$ is a concatenation of factors of the form $\alpha_t^{\tilde \jmath}\beta_t$ or $\alpha_t\beta_t^{\tilde \jmath}$. The number of such factors is at most $|w \tilde w \overline w|/|\alpha_t\beta_t|\leq(13r+10)/|\alpha_t\beta_t|\leq 25\sqrt{\log r}$. Moreover, if we have to consecutive such factors $\alpha_t^{\tilde \jmath_1}\beta_t$ and $\alpha_t^{\tilde \jmath_2}\beta_t$ (or $\alpha_t\beta_t^{\tilde \jmath_1}$ and $\alpha_t\beta_t^{\tilde \jmath_2}$), then $|\tilde \jmath_1-\tilde \jmath_2|\le 1$, and if we have two consecutive factors $\beta_t\alpha_t^{\tilde \jmath_1}\beta_t$ and $\alpha_t\beta_t^{\tilde \jmath_2}\alpha_t$ then $2\le |j_1|+|j_2| \le 3$. This implies that each of these factors of the form $\alpha_t^{\tilde \jmath}\beta_t$ or $\alpha_t\beta_t^{\tilde \jmath}$ has at most $3$ continuations of this form, and so the number of such continuations $\tilde w \overline w$ of $w$ is at most $3^{25\sqrt{\log r}}<r$. Since the number of possible $w\in \Sigma^{(r)}(3+e^{-r-4})$ is $\aO(r^3)$ then the number of possible continuations $w\tilde w \overline w$ is this case is $\aO(r^4)$.
		
		\medbreak
		
		\noindent\textbf{Case 2:} Suppose now that $w$ has a factor $\alpha_t^s$ with $|\alpha_t^s|\geq|w|/3$  (the case of $w$ having a factor $\beta_t^s$ with $|\beta_t^s|\ge |w|/3$ will be analogous). \\
		
		Let us check that we can apply \Cref{lem:uvextension,lem:digitextension} to this factor.  Observe that $r/6-1/3\leq|w|/3\leq s|\alpha_t|$ whence $s> \sqrt{\log r}/12$, so the condition $(T+2)r\leq s^2|\alpha_t|/8$ holds for $T=10$ and $r$ large. In particular $\ell\leq 2.1Tr/|\alpha_t^s|+1< 150$. Hence one has that $(s-\ell/2)|\alpha_t|> r/7$ for sufficiently large $r$. Going back to the inequality $1.1Tr+4\geq (\ell-1)(s-\ell/2)|\alpha_t|$ we get a stronger bound $\ell< (1.1Tr+4)/((s-\ell/2)|\alpha_t|)+1\leq 80$. We consider two subcases depending on the length of $\alpha_t$.
		\smallbreak
		
		\noindent\textbf{Case 2.1:} Suppose that $|\alpha_t|>r^{15/16}$ and  $\alpha_t=uv$ with $(u,v)\in\overline{P}$. 
		
		For big $r$ we can assume further that $\sizer(\alpha_t^s)\leq r-4|\alpha_t|$ holds, since $9/10|w|\leq r+3$. Then, we claim that after $\alpha_t^s=(uv)^s$ the word to be renormalizable with alphabet $\{u,v\}$, and the first appearance of $uu$ or $vv$ determines the new alphabet ($\{u,uv\}$ or $\{uv,v\}$). To prove that, let $\hat{w}\in\Sigma^{((T+2)r)}(3+e^{-r-4})$  be the factor of $w\tilde{w}\overline{w}$ starting at that factor $\alpha_t^s$.  Therefore \Cref{lem:uvextension} gives that $\hat{w} = \hat{\gamma} \hat{w_2}$, where
		\[
		\hat{\gamma} = (uv)^{s_1} \theta_1 (uv)^{s_2} \theta_2 \ldots (uv)^{s_\ell},
		\]
		and  each $\theta_j$ belongs to $\{u u v, u v v \}$ and moreover $\ell\leq 80$.
		
		In particular given a continuation $\tilde w \overline w$ of $w$ with $\tilde w\in\Sigma^{(r)}(3+e^{-r-4})$ and $\overline w \in \Sigma^{(10r)}(3+e^{-r-4})$,  from the first such factor $(uv)^{s_1}$, the sequence should be \[
		(uv)^{s_1}\theta_1(uv)^{s_2}\theta_2\dots(uv)^{s_\ell},
		\]
		with $\theta_j\in\{uuv, uvv\}$, $\ell\le 80$ and $s_1+\dots+s_\ell\le 20 r^{1/16}$, so we have in total at most $2^{80}$ choices for the $\theta_j$, and, given $\ell\le 80$, the number of choices for the $s_j$ is at most the number of natural solutions of $x_1+x_2+\dots+x_{\ell+1}=\lfloor 20 r^{1/16}\rfloor\equalscolon M$, which is $\binom{M + \ell}{\ell}<(21 r^{1/16})^{80}=21^{80}r^{5}$, and so the total number of such words $w\tilde w\overline w$ is $\aO(r^3\cdot 2^{80}\cdot 80\cdot 21^{80}r^{5})=\aO(r^{8})$.
		
		\smallbreak
		
		\noindent\textbf{Case 2.2:} Suppose $|\alpha_t|\leq r^{15/16}$ and that the largest factor $\alpha_t^{s}$ of $w\tilde{w}$ satisfies $\sizer(\alpha_t^{s})\leq r-170|\alpha_t|$.
		
		We claim there is $\beta$ such that $(\alpha_t,\beta)\in\overline{P}$,  $\beta(\alpha_t^{s_1})\beta$ is a subword of $w\tilde{w}$, and the continuation $w\tilde{w}\overline{w}$ has the form $\beta(\alpha_t)^{s_1}\beta(\alpha_t)^{s_2}\ldots(\alpha_t)^{s_\ell}$ with $|s_j-s_{j+1}|\leq 1$ and for all $1\leq j\leq\ell\leq 80$. By hypothesis there is $\beta_t(\alpha_t)^s\beta_t$ inside $w\Tilde{w}$. If $(\alpha_t,\beta_t)=(uv,v)$ for some $(u,v)\in\overline{P}$ or $(\alpha_t,\beta)=(a,b)$ then we set $\beta=\beta_t$, while if $\beta_t=\alpha_t^k\tilde{\beta}$ with $(\alpha_t,\Tilde{\beta})\in\overline{P}$ and $|\Tilde{\beta}|\leq|\alpha_t|$, then we set $\beta=\Tilde{\beta}$.
		
		If $(\alpha_t,\beta)=(uv,v)$ with $(u,v)\in\overline{P}$ then we use \Cref{lem:uvextension} to obtain a continuation $\hat{w}$ with  $\hat{\gamma}=(uv)^{s_1}\theta_1(uv)^{s_2}\theta_2\ldots$ with $\theta_i\in\{uuv,uvv\}$ and $\hat{w}=\hat{\gamma}w_2$. Moreover $\ell\leq 80$. Since $\sizer(\alpha_t^{s_1})\leq r-170|\alpha_t|$, by induction $\sizer(\alpha_t^{s_j})\leq r-(172-2j)|\alpha_t|$ and since $\ell\leq80$ by the proof of that \Cref{lem:uvextension} we get that all $\theta_j$ are equal to $\theta_j=\alpha_t\beta=uvv$ (since $\alpha_t^s\beta_t=(uv)^fv$ for some $f\geq 1$ and the fact that $uuv$ and $uv$ start with $v_a$)  and that $|s_j-s_{j+1}|\leq 1$ for all $j\geq 1$. 
		
		When $\alpha_t=a$, we use \Cref{lem:digitextension} to find a continuation $\hat{w}$ such that $\hat{w}=2^{e_1}b 2^{e_2}b\ldots$. But observe that $2^{e_1}=\alpha_t^{e_1/2}$ is inside $w$, so by hypothesis $\sizer(2^{e_1})<r-170|\alpha_t|$ and there is $b$ before $2^{e_1}$, so $e_2$ is even and $|e_1-e_2|\in\{0,2\}$. By induction we obtain $\sizer(2^{e_j})\leq r-(172-2j)|\alpha_t|$, which forces all $e_j$ to be even and $|e_{j+1}-e_j|\in\{0,2\}$. In this case $\beta=b$, so we get $\hat{w}=2^{e_1}b 2^{e_2}b\ldots=\alpha_t^{s_1}\beta\alpha_t^{s_2}\beta\ldots$ with $|s_{j+1}-s_j|\leq 1$ for all $j\geq 1$.
		
		Therefore, given a continuation $\tilde w \overline w$ of $w$ with $\tilde w\in\Sigma^{(r)}(3+e^{-r-4})$ and $\overline w \in \Sigma^{(10r)}(3+e^{-r-4})$, there is $\beta$ with $(\alpha_t,\beta)\in \overline P$ such that $w\tilde w$ has a factor $\beta(\alpha_t)^{s_1}\beta$, after which the continuation of $w\tilde w \overline w$ is a concatenation of at most $79$ sequences of the type $(\alpha_t)^{s_j}\beta$, $2\le j\le 80$ with $|s_{j+1}-s_j|\le 1$ for every $j\ge 1$. This gives at most $3^{80}$ continuations of $\beta(\alpha_t)^{s_1}\beta$, and so, since we have at most $\aO((r^3)^2)=\aO(r^6)$ choices for $w\tilde w$, we have in total, $\aO(2^{80}\cdot r^6)=\aO(r^6)$ such words $w\tilde w \overline w$.
		
		\smallbreak
		
		In all the previous cases, Case 1, Case 2.1 and Case 2.2, if $d=\frac{\log r-\log\log r}r$, then
		\[
		(e^{-10r})^d=e^{-10(\log r-\log\log r)}=\left(\frac{r}{\log r}\right)^{-10}.
		\]
		Moreover, in these cases we have $\aO(r^8)$ possible such words $w\tilde w \overline w$. Notice that
		\[
		r^8(e^{-10r})^d=r^8\cdot \left(\frac{r}{\log r}\right)^{-10}=\frac{\log^{10} r}{r^2}<\frac1{r}\ll 1.
		\]
		\medbreak
		
		Our third case is derived from Case 2.2, but it is more delicate.
		
		\noindent\textbf{Case 3:}  On the same conditions of Case 2,  suppose that $|\alpha_t|\leq r^{15/16}$ and that $w\tilde{w}$ has a factor $\alpha_t^{s_1}$ satisfying $\sizer(\alpha_t^{s_1})\geq r-170|\alpha_t|$. \\
		We will consider in this case continuations $\overline w \in \Sigma^{(Tr)}(3+e^{-r-4})$ for $T=\lfloor \log^2 r\rfloor$ such that $w\tilde w\overline w\in \Sigma(3+e^{-r-4},|w\tilde w\overline w|)$. Again, consider two subcases depending on the length of $\alpha_t$.
		
		\smallbreak
		
		\noindent\textbf{Case 3.1} Suppose that that $|\alpha_t|>2$.
		
		So,  $\alpha_t=uv$ with $(u,v)\in\overline{P}$. Now let $T=\lfloor(\log r)^2\rfloor$. The condition on $\alpha_t^{s_1}$ implies that $s_1\geq (r-170|\alpha_t|)/(2|\alpha_t|)\geq (1/2)r^{1/16}-85$, so we have that $s_1^2|\alpha_t|/2\geq s_1(r-170r^{15/16})/4\geq(\log r)^2r\geq Tr$.  Let $\hat{w}\in\Sigma^{((T+2)r)}(3+e^{-r-4})$ be the factor of $w\tilde{w}\overline{w}$ starting at that factor $\alpha_t^{s_1}$. \Cref{lem:uvextension} guarantees that
		\[
		\hat{w}=(uv)^{s_1}\theta_1(uv)^{s_2}\ldots(uv)^{s_\ell}\hat{w}_2
		\]
		where each $\theta_i\in\{uuv,uvv\}$ and $\ell\leq 2.1Tr/|\alpha_t^{s_1}|+1\leq 5T+2$ for $r$ big enough. 
		
		Therefore, from this factor $\alpha_t^{s_1}$, the continuation of $w\tilde w\overline w$ is an initial factor of a word of the form $(uv)^{s_1}\theta_1(uv)^{s_2}\theta_2\ldots\theta_{\ell-1}(uv)^{s_\ell}$, and 
		\[
		Tr\le \sizer((uv)^{s_1}\theta_1(uv)^{s_2}\theta_2\dots(uv)^{s_{\ell}})<(T+3)r,
		\]
		with $\theta_j\in\{uuv, uvv\}$, $\ell\le 5T+2$ (and such that $(uv)^{s_1}\theta_1(uv)^{s_2}\theta_2\ldots\theta_{\ell-1}$ is an initial factor of this word beginning in this factor $\alpha_t^{s_1}$ and going till the end of $w\tilde w\overline w$) with $\sizer(\alpha_t^{s_j})\ge r-(173+T)|\alpha_t|>r-\lfloor 6\log^2 r \cdot r^{15/16}\rfloor\equalscolon M$ (notice that if $\sizer((uv)^{s_j})<r-10|\alpha_t|$ then $|s_{j+1}-s_j|\le 1$). Let $s_0$ be the smallest integer satisfying $\sizer(\alpha_t^{s_0})\ge M$. Then, $s_j=s_0+\tilde s_j$ with $\tilde s_j\ge 0$ for each $1 \leq j\le \ell$. 
		
		Since $q_{2s|\alpha_t|}(\alpha_t^{s})\ge q_{2|\alpha_t|}(\alpha_t)^{s}$ and $q_{2|\alpha_t|}(\alpha_t)\ge q_4(1122)=12$, we have $\sizer(\alpha_t^{s})\ge \lfloor \log(q_{2s.|\alpha_t|}(\alpha_t^{s})^2)\rfloor\ge \lfloor\log((12^{s})^2)\rfloor=\lfloor s \log(144)\rfloor>4s$. Hence
		\begin{align*}
			(T+3)r&>\sizer((uv)^{s_1}\theta_1(uv)^{s_2}\theta_2\dots(uv)^{s_{\ell}})\ge \ell \sizer(\alpha_t^{s_0})+(\tilde s_1+\tilde s_2+\dots+\tilde s_{\ell})\sizer(\alpha_t)\\
			&\ge \ell M+4(\tilde s_1+\tilde s_2+\dots+\tilde s_{\ell}).
		\end{align*}
		In particular $\ell\leq T+3$ for $r$ large. Since $(T+4)M=(T+4)(r-\lfloor 6\log^2 r \cdot r^{15/16}\rfloor)>(T+3)r$, it follows that, given $1\le \ell\le T+3$, the number of choices of the $s_j, j\le \ell$ is at most the number of natural solutions of 
		\[
		\tilde s_1+\tilde s_2+\dots+\tilde s_{\ell}\le (T+4-\ell)M/4,
		\]
		which is
		$\binom{\lfloor (T+4-\ell)M/4\rfloor+\ell}{\ell}\le \left(\frac{e((T+4-\ell)M/4+\ell)}{\ell}\right)^{\ell}<\left(\frac{e((\tilde T-\ell)M/4)}{\ell}\right)^{\ell}$, where $\tilde T\colonequals T+5$
		(here we used the inequalities $\binom{n}{k}\le \frac{n^k}{k!} \le \left(\frac{en}k\right)^k$, which hold for $1\le k\le n$). We have at most $2^{\ell}$ choices for the $\theta_j, j\le \ell$; let us estimate the maximum of $f(\ell)=2^{\ell}\left(\frac{e((\tilde T-\ell)M/4)}{\ell}\right)^{\ell}=\left(\frac{e((\tilde T-\ell)M/2)}{\ell}\right)^{\ell}$ for $1\le \ell\le \tilde T-1$. The derivative of $\log f(\ell)$ is $\log\left(\frac{e((\tilde T-\ell)M/2)}{\ell}\right)-\frac{\tilde T}{\tilde T-\ell}$. Since, in this range of $\ell$,
		\[
		\log\left(\frac{e((\tilde T-\ell)M/2)}{\ell}\right)=(1+\ao(1))\log M=(1+\ao(1))\log r
		\]
		and $\tilde T=T+5=\log^2 r+\aO(1)$, we have the maximum attained for $\ell=\tilde T \left(1-\frac{1+\ao(1)}{\log r}\right)=\log^2 r-(1+\ao(1))\log r<T$, and, for such value of $\ell$,
		\begin{align*}
			f(\ell)=\left(\frac{e((\tilde T-\ell)M/2)}{\ell}\right)^{\ell}&=\left(\frac{e((1+\ao(1))M/2)}{\log r}\right)^{\ell}\\
			&<\left(\frac{3M}{2\log r}\right)^{T-(1+\ao(1))\log r}.
		\end{align*}
		We have at most $\tilde T=\log^2 r+\aO(1)$ choices for $\ell$, and we have at most $\aO(r^6)$ choices for $w\tilde w$, so we have at most
		\[
		\aO\left(r^6 \log^2 r \left(\frac{3M}{2\log r}\right)^{T-(1+\ao(1))\log r}\right)<\left(\frac{2r}{\log r}\right)^{T-(1+\ao(1))\log r}
		\]
		such words $w\tilde w\overline w$. 
		
		Notice that, for $d=\frac{\log r-\log\log r}r$, we have $(e^{-Tr})^d=e^{-T(\log r-\log\log r)}$, so \[
		\left(\frac{2r}{\log r}\right)^T (e^{-Tr})^d=\left(\frac{2r}{\log r} e^{-\log r+\log\log r}\right)^T=2^T,
		\]
		and 
		\begin{align*}
			\left(\frac{2r}{\log r}\right)^{T-(1+\ao(1))\log r}(e^{-Tr})^d &\le  \left(\frac{2r}{\log r}\right)^{-(1+\ao(1))\log r}\cdot 2^{\log^2 r} \\ 
			&=e^{-(1+\ao(1))\log^2 r}\cdot e^{\log 2\log^2 r}\\
			&=e^{-(1-\log 2+\ao(1))\log^2 r}\\
			&<e^{-\frac14\log^2 r}\ll 1.
		\end{align*}
		
		\smallbreak
		
		\noindent\textbf{Case 3.2:} Suppose that $|\alpha_t|=2$.\\
		Then, for some $c\in \{1, 2\}$, $w\tilde w$ has a factor $c^{s_1}$ satisfying $\sizer(c^{s_1})\ge r-170$. Let $c'=3-c\in\{1, 2\}$ and $\theta=c'c'$. Observe that $(s+1)\log x\geq \sizer(c^{s})\geq r-170$, so  $(T+2)r\leq s^2(\log x)/4$ holds. Using \Cref{lem:digitextension}, from this factor $c^{s_1}$ the continuation of $w\tilde w\overline w$ has the form $c^{s_1}\theta c^{s_2}\theta\ldots \theta c^{s_\ell}$ with $\ell<2Tr/((s-4)\log x)<2Tr/(r-170)<3T$ for large $r$. Using this information in the inequality $\ell(s-\ell-3)\log x<Tr$, gives that $\ell\leq T+1$ for large $r$. 
		Therefore, from this factor $c^{s_1}$, the continuation of $w\tilde w\overline w$ is an initial factor of a word of the form $c^{s_1}\theta c^{s_2}\theta\ldots \theta c^{s_\ell}$, $Tr\le \sizer(c^{s_1}\theta c^{s_2}\theta\ldots c^{s_{\ell}})\le (T+1)r$, $\ell\le T+1$ (and such that $c^{s_1}\theta c^{s_2}\theta\ldots c^{s_{\ell-1}}\theta$ is an initial factor of this word beginning in this factor $c^{s_1}$ and going till the end of $w\tilde w\overline w$) with $\sizer(c^{s_j})\ge r-(171+2T)>r-\lfloor 3\log^2 r\rfloor\equalscolon N$ (notice that if $\sizer(c^{s_j})<r-7$ then $s_j$ is even and $|s_{j+1}-s_j|\in \{0,2\}$). Let $s_0$ minimum such that $\sizer(c^{s_0})\ge N$. Then $s_j=s_0+\tilde s_j$ with $\tilde s_j\ge 0$ for each $1 \leq j\le \ell$. Notice that, given $c$, $\overline w$ is determined by the choice of $(\ell,s_1,s_2,\dots,s_\ell)$. 
		
		To estimate the number of the corresponding possibilities, we will make use of \Cref{lem_Comb}.  We will consider two last subcases depending on the value of $c$.
		
		\smallbreak
		
		\noindent\textbf{Case 3.2.1:} Assume that $c=2$:
		
		Since $q_{s}(2^{s})\ge 2^{s}$, we have
		\[
		\sizer(2^{s})\ge \lfloor \log(q_{s}(2^{s})^2)\rfloor\ge \lfloor\log((2^{s})^2)\rfloor=\lfloor s \log(4)\rfloor\ge 4s/3-1.
		\]
		We have
		\begin{align*}
			(T+1)r>\sizer(2^{s_1}112^{s_2}11\ldots 2^{s_{\ell}})&\ge \ell \sizer(2^{s_0})+4(\tilde s_1+\tilde s_2+\dots+\tilde s_{\ell})/3-\ell\\
			&\ge \ell N+4(\tilde s_1+\tilde s_2+\dots+\tilde s_{\ell})/3-\ell.
		\end{align*}
		Since
		\[
		(T+2)N=(T+2)(r-\lfloor 3\log^2 r\rfloor)>(T+1)r+(T+1)\ge (T+1)r+\ell,
		\]
		it follows that, given $1\le \ell\le T+1$, the number of choices of the $s_j, j\le \ell$ is at most the number of natural solutions of $\tilde s_1+\tilde s_2+\dots+\tilde s_{\ell}\le 3(T+2-\ell)N/4$. By \Cref{lem_Comb}, it is at most
		\[
		e^{(\log (3N/4-\log\log (3N/4)+\ao(1))(T+3)}=e^{(\log N-\log\log N-\log(4/3)+\ao(1))T}.
		\]
		Since $\log N=\ao(T)$. We have at most $\aO(r^6)$ choices for $w\tilde w$, so we have at most
		\[
		\aO(r^6 e^{(\log N-\log\log N-\log(4/3)+\ao(1))T}=\aO(e^{(\log N-\log\log N-\log(4/3)+\ao(1))T})
		\]
		such words $w\tilde w\overline w$. 
		
		Notice that, for $d=\frac{\log r-\log\log r}r$, we have
		$(e^{-Tr})^d=e^{-T(\log r-\log\log r)}$, and so, since $\log N=\log r+\ao(1)$, 
		\begin{align*}
			e^{(\log N-\log\log N-\log(4/3)+\ao(1))T} (e^{-Tr})^d&=e^{T(\log r-\log\log r-\log(4/3)+\ao(1))}(e^{-Tr})^d \\
			& = e^{T(\ao(1)-\log(4/3))}<e^{-\frac{\log^2 r}{4}}\ll 1.
		\end{align*}
		
		\smallbreak
		
		\noindent\textbf{Case 3.2.2:} Assume that $c=1$.
		
		Observe that
		\[
		(T+1)r>\sizer(1^{s_1}221^{s_2}22\ldots 221^{s_{\ell}})\ge \ell N+(\tilde s_1+\tilde s_2+\dots+\tilde s_{\ell})\log\left(\frac{3+\sqrt5}2\right)
		\]
		by \eqref{eq:intervalswith1s} and \eqref{eq:rofpower1}. Since
		\[
		(T+2)N=(T+2)(r-\lfloor 3\log^2 r\rfloor)>(T+1)r,
		\]
		it follows that, given $1\le \ell\le T+1$, the number of choices of the $s_j, j\le \ell$ is at most the number of natural solutions of
		\[
		\tilde s_1+\tilde s_2+\dots+\tilde s_{\ell}\le (T+2-\ell)N/\log\left(\frac{3+\sqrt5}2\right).
		\]
		By \Cref{lem_Comb}, it is at most 
		\begin{align*}
			&e^{\left(\log \left(N/\log\left(\frac{3+\sqrt5}2\right)\right)-\log\log \left(N/\log\left(\frac{3+\sqrt5}2\right)\right)+\ao(1)\right)(T+3)} \\
			{}={}&e^{\left(\log \left(N/\log\left(\frac{3+\sqrt5}2\right)\right)-\log\log \left(N/\log\left(\frac{3+\sqrt5}2\right)\right)+\ao(1)\right)T},
		\end{align*}
		since $\log N=\ao(T)$. We have at most $\aO(r^6)$ choices for $w\tilde w$, so we have at most
		\begin{align*}
			&\aO\left(r^6 e^{\left(\log \left(N/\log\left(\frac{3+\sqrt5}2\right)\right)-\log\log \left(N/\log\left(\frac{3+\sqrt5}2\right)\right)+\ao(1)\right)T}\right) \\
			{}={}&\aO\left(e^{\left(\log \left(N/\log\left(\frac{3+\sqrt5}2\right)\right)-\log\log \left(N/\log\left(\frac{3+\sqrt5}2\right)\right)+\ao(1)\right)T}\right)
		\end{align*}
		such words $w\tilde w\overline w$. 
		
		Notice that, if $\delta>0$, for $d=\frac{\log r-\log\log r-\log\log\left(\frac{3+\sqrt5}2\right)+\delta}r$, we have
		\[
		(e^{-Tr})^d=e^{-dTr}=e^{-T\left(\log r-\log\log r-\log\log\left(\frac{3+\sqrt5}2\right)+\delta\right)},
		\]
		and so, since $\log N=\log r+\ao(1)$,
		\begin{align*}
			&e^{\left(\log (N/\log\left(\frac{3+\sqrt5}2\right)\right)-\log\log \left(N/\log\left(\frac{3+\sqrt5}2\right)\right)+\ao(1))T} (e^{-Tr})^d \\
			{}={} & e^{T\left(\log r-\log\log r-\log\log\left(\frac{3+\sqrt5}2\right)+\ao(1)-\log r+\log\log r+\log\log\left(\frac{3+\sqrt5}2\right)-\delta\right)}\\
			{}={} & e^{T(\ao(1)-\delta)}<e^{-\frac{\delta\log^2 r}{2}}\ll 1.
		\end{align*}
		Since $c_0\colonequals -\log\log\left(\frac{3+\sqrt5}2\right)=0.03830054\ldots>0$, it follows that
		\[
		d(3+e^{-r})\le 2\cdot\frac{\log r-\log\log r+c_0+\ao(1)}r.
		\]
		
		\bigskip
		
		Up to this point of the proof, we have shown the upper bound
		\[
		d(3+t)\le 2\cdot\frac{\log |\log t|-\log\log |\log t|+c_0+\ao(1)}{|\log t|},
		\]
		which gives us a different proof of the upper bound on the easier bounds stated in the introduction. In fact, the only case that gives the worst bound is the last one with $c=1$ (that is, Case 3.2.2).
		
		\smallbreak
		
		We can actually obtain a more precise upper estimate by choosing $T=\lfloor r/5 \rfloor$, which is what we will do now. For the sake of exposition, we will consider this improved estimate to be a separate case.
		
		\noindent\textbf{Case 3.2.3:} We will derive a more precise estimate for the case $c = 1$.
		
		Observe that it is possible to chose $T=\lfloor r/5 \rfloor$ in \Cref{lem:digitextension}, because $\sizer(1^s)\geq r-170$ gives us that $s\log x\geq r-170$ so one has that
		\[
		s^2(\log x)/4\geq (r-170)^2/(4\log x)\geq r^2/5
		\]
		for large $r$. Moreover
		\[
		\ell<2Tr/(s\log x)\leq 2.1Tr/(r-160)\leq 5/2T\leq r/2
		\]
		for large $r$. Putting this again in the inequality $\ell<Tr/((s-\ell+1)\log x)$ gives further that $\ell<2T+1$, so $\ell\leq 2T$ for large $r$.

		Let $T=\lfloor r/5 \rfloor$. We would have a worst lower estimate for $\sizer(\alpha_t^{s_i})$: for $i\ge 1$, we have $\sizer(\alpha_t^{s_i})\ge r-2(173+i)\ge r/3$. Indeed,
		\begin{align*}
			r^2/5+r&\ge \sizer(1^{s_1}221^{s_2}22\ldots 221^{s_{\ell}})\\
			&\ge \sum_{i=1}^{\min\{\ell,r/2\}}(r-2(173+i))\\
			&=\min\{\ell,r/2\}(r-347-\min\{\ell,r/2\}),
		\end{align*}
		which implies $\ell<3r/10$, and thus $\sizer(\alpha_t^{s_i})\ge r-2(173+i)>r/3$. We will introduce a parameter $j$ equal to the number of values of $i$ for which $\sizer(\alpha_t^{s_i})<r-3$ in $1^{s_1}221^{s_2}22\ldots 221^{s_{\ell}}$, for which we should have $s_{i+1}\in\{s_i,s_i-2,s_i+2\}$ (for the other $\ell-j$ values of $1\le i\le \ell$ we have $\sizer(\alpha_t^{s_i})\ge r-3$); if we consider these $j$ values $i_1<i_2<\dots<i_j$ of $i$, we have $s_{i_1}\ge s_0-100$, so $s_{i_t}>s_0-100-2t, 1\le t\le j$, and $\sum_{1\le i\le j}s_{i_t}>j\cdot (s_0-100-j)$. 
		
		Let $\hat \ell=\ell-j$ and $\{s_i, i\in I=\{1,2,\dots,\ell\}\setminus \{i_t, 1\le t\le j\}=\{\hat s_1,\hat s_2,\dots,\hat s_{\hat\ell}\}$. We have 
		$\ell<3r/10<2T$. Given $\hat\ell$ and $j$ there are at most
		\[
		\binom{\ell}{j}=\binom{\hat{\ell}+j}{j}<\left(\frac{e\ell}j\right)^j<\left(\frac{2eT}j\right)^j
		\]
		choices for the set $\{s_{i_t}, 1\le t\le j\}$. Since for $i\in \{i_t, 1\le t\le j\}$ we have at most $3$ choices for $s_{i+1}$, and the total number of these choices is at most $3^j$. Together with the number of choices for the set $\{s_{i_t}, 1\le t\le j\}$, this gives an estimate of $\left(\frac{6eT}j\right)^j$ for these choices.
		
		Let $\hat s_0$ be the smallest integer that satisfies $\sizer(1^{\hat s_0})\ge r-3$. Then, we have $\hat s_0>(r-5)/\log\left(\frac{3+\sqrt5}2\right)$. The number of solutions of the above inequality is at most the number of natural solutions of
		\begin{align*}
			\hat s_1+\hat s_2+\dots+\hat s_{\hat\ell}&\le (T+2-\hat\ell)(r-5)/\log\left(\frac{3+\sqrt5}2\right)-j\cdot (s_0-100-j)\\
			&<(T+2-\hat\ell-j\cdot (r-104-j)/r)(r-5)/\log\left(\frac{3+\sqrt5}2\right)\\
			&<(T+2-\hat\ell-j/2)(r-5)/\log\left(\frac{3+\sqrt5}2\right).
		\end{align*}
		
		By \Cref{lem_Comb}, the number of solutions of
		\[
		\hat s_1+\hat s_2+\dots+\hat s_{\hat\ell}\le (T+2-\hat\ell-j\cdot (r-104-j)/r)(r-5)/\log\left(\frac{3+\sqrt5}2\right)
		\]
		is at most
		\[
		(T+2)e^{(1-\varepsilon_m)(T+3-j\cdot (r-104-j)/r)/\varepsilon_m},
		\]
		where $\varepsilon_m$ is the solution in $(0,1)$ of the equation
		\[
		\log\left(\frac{em\varepsilon}{1-\varepsilon}\right)=\frac1{\varepsilon}, \qquad \text{ with }  m=(r-5)/\log\left(\frac{3+\sqrt5}2\right).
		\]
		
		Since $\frac{em\varepsilon_m}{1-\varepsilon_m}>\frac{r}{2\log r}$, the factor $\left(\frac{em\varepsilon_m}{1-\varepsilon_m}\right)^{-j\cdot (r-104-j)/r}$ is such that
		\[
		\left(\frac{6eT}j\right)^j\left(\frac{em\varepsilon_m}{1-\varepsilon_m}\right)^{-j\cdot (r-104-j)/r}<\left(\frac{6eT}j\left(\frac{r}{2\log r}\right)^{-(r-104-j)/r}\right)^j.
		\]
		This is smaller than $\left(\frac{6eT}j\left(\frac{r}{2\log r}\right)^{-1/2}\right)^j$, and for $j\ge r^{3/4}$ this is $\ao(1)$ (using $T\le r/4$). For $10\log r\le j<r^{3/4}$, the estimate
		\[
		\left(\frac{6eT}j\left(\frac{r}{2\log r}\right)^{-(r-104-j)/r}\right)^j
		\]
		will be $\ao(1)$ since $-(r-104-j)/r<-1+r^{-1/5}$ and $\left(\frac{r}{2\log r}\right)^{r^{-1/5}}=1+\ao(1)$, so the estimate becomes $\left((3+\ao(1))e/10\right)^{10\log r}=\ao(1)$. On the other hand, for $0\le j<10\log r$, the estimate $\left(\frac{6eT}j\left(\frac{r}{2\log r}\right)^{-(r-104-j)/r}\right)^j$ becomes $\left(\frac{(3+\ao(1))e\log r}j\right)^j<\left(\frac{9\log r}j\right)^j$. The maximum of the function $v(j)=\left(\frac{9\log r}j\right)^j$ is attained at $j=9\log r/e$, and is equal to $e^{9\log r/e}<r^4$. So, using again the fact that we have $\aO(r^6)$ choices for $w\tilde w$, in any case we get an upper estimate for the total number of words $w\tilde w\overline w$ which is
		\[
		\aO(r^6)\cdot r^4\cdot (T+2)\cdot \left(\frac{em\varepsilon_m}{1-\varepsilon_m}\right)^{(1-\varepsilon_m)(T+3)}=\aO(r^{14})\cdot \left(\frac{em\varepsilon_m}{1-\varepsilon_m}\right)^{(1-\varepsilon_m)T}.
		\]
		As before, this gives an upper estimate for the dimension which is 
		\[
		\frac{(1-\varepsilon_m)\log\left(\frac{em\varepsilon_m}{1-\varepsilon_m}\right)+\aO(\log r/T)}r=\frac{(1-\varepsilon_m)\log\left(\frac{em\varepsilon_m}{1-\varepsilon_m}\right)+\aO(\log r/r)}r
		\]
		
		Since $\varepsilon_m$ is the solution in $(0,1)$ of the equation $\log\left(\frac{em\varepsilon_m}{1-\varepsilon_m}\right)=\frac1{\varepsilon_m}$, with $m=(r-5)/\log\left(\frac{3+\sqrt5}2\right)$, $(1-\varepsilon_m)\log\left(\frac{em\varepsilon_m}{1-\varepsilon_m}\right)=\frac{1-\varepsilon_m}{\varepsilon_m}$. Writing $z=\frac{1-\varepsilon_m}{\varepsilon_m}$, the equality $\log\left(\frac{em\varepsilon_m}{1-\varepsilon_m}\right)=\frac1{\varepsilon_m}$ can be written as $\log\left(\frac{em}{z}\right)=z+1$, so $z+\log z=\log m$ and $ze^z=m$, so $z=W(m)=W\left((r-5)/\log\left(\frac{3+\sqrt5}2\right)\right)$, where $W$ is Lambert's function. Since $W'(x)<1/x$,
		\[
		W\left((r-5)/\log\left(\frac{3+\sqrt5}2\right)\right)=W\left(r/\log\left(\frac{3+\sqrt5}2\right)\right)+\aO(1/r),
		\]
		and our upper estimate for the dimension is
		\[
		z/r+\aO\left(\frac{\log r}{r^2}\right)=W\left(r/\log\left(\frac{3+\sqrt5}2\right)\right)/r+\aO\left(\frac{\log r}{r^2}\right).
		\]

	\end{proof}

	\section{The lower bound} \label{sec:lower_bound}
	
	The statements and definitions below are taken from the third author's work \cite{M:geometric_properties_Markov_Lagrange}. 
	
	\begin{definition} Given $B=\{\beta_1,\dots,\beta_\ell\}$, $\ell\geq 2$, a finite alphabet of finite words $\beta_j\in(\mathbb{N}^*)^{r_j}$, which is primitive (in the sense that $\beta_i$ does not begin by $\beta_j$ for all $i\neq j$) then the Gauss-Cantor set $K(B)\subseteq [0,1]$ associated with $B$ is defined as 
		\[
		K(B)\colonequals \{[0;\gamma_1, \gamma_2, \dots] \ \mid\ \gamma_i\in B\}.
		\]
	\end{definition}
	
	The set $K(B)$ is a dynamically defined Cantor set. We will now exhibit its Markov partition and the expanding map which defines it.
	
	For each word $\beta_j\in(\mathbb{N}^*)^{r_j}$, let $I_j=I(\beta_j)$ be the convex hull of the set $\{[0;\beta_j, \gamma_1, \gamma_2, \dots]\ \mid\ \gamma_i\in B\}$ and $\psi|_{I_j}\colonequals G^{r_j}|_{I_j}$ where
	\[
	G(x)=\{1/x\}=1/x-\lfloor 1/x\rfloor
	\]
	is the Gauss map. This defines an expanding map $\psi\colon I(\beta_1)\cup\dotsb\cup I(\beta_\ell)\to I$. Let $I=[\min K(B), \max K(B)]$. Then $I$ is the convex hull of $I_1\cup\dotsb\cup I_\ell$ and $\psi(I_j)=I$ for every $j\le \ell$.

	Let us describe how to estimate $\dimH(K(B))$. 
	\bigskip
	
	According to Palis--Takens \cite[Chapter 4]{PalisTakens}, let 
	\[
	\lambda_j=\inf |\psi'|_{I_j}|, \qquad \Lambda_j=\sup |\psi'|_{I_j}|
	\]
	and $\alpha,\beta\geq 0$ be such that
	\[
	\sum_{i=1}^{\ell}\lambda_j^{-\alpha}=1, \qquad \sum_{i=1}^{\ell}\Lambda_j^{-\beta}=1.
	\]
	Then,
	\begin{equation}\label{Inq-Palis-Takens}
		\beta\leq \dimH(K(B))\leq\alpha.
	\end{equation}
	\bigskip
	Let us discuss how to find estimates for $\alpha$ and $\beta$.
	
	\noindent The iterates of the Gauss map are given explicitly by
	\[
	\psi|_{I_j}(x)=\dfrac{q^{(j)}_{r_{j}}x - p^{(j)}_{r_{j}}}{p^{(j)}_{r_j-1} - q^{(j)}_{r_j-1}x}
	\]
	where $\dfrac{p^{(j)}_k}{q^{(j)}_k} = [0; b^{(j)}_1, \dots, b^{(j)}_k]$ and $\beta_j = (b^{(j)}_1, \dots, b^{(j)}_{r_j})$. 
	
	\noindent Hence 
	\[
	(\psi|_{I_j})'(x)=\frac{(-1)^{r_j-1}}{(p^{(j)}_{r_j-1} - q^{(j)}_{r_j-1}x)^2}.   
	\]
	
	\begin{lemma}\label{lemma_aprox}
		Let $x=[c_0,c_1,c_2,\dots]$ and $\frac{p_n}{q_n}=[c_0,c_1,\dots, c_n]$. Then
		\[
		\frac{1}{2q_nq_{n+1}}<\frac{1}{q_n(q_n+q_{n+1})}<\left|x-\frac{p_{n}}{q_{n}}\right|<\frac{1}{q_nq_{n+1}},
		\]
		and therefore 
		\[
		\frac{1}{2q_{n+1}}<|q_n x-p_n|<\frac{1}{q_{n+1}}.
		\]
	\end{lemma}
	
	\noindent Therefore, \Cref{lemma_aprox} implies that 
	\[
	(q_{r_j}^{(j)})^2<|(\psi|_{I_j})'(x)|=\frac{1}{(p^{(j)}_{r_j-1} - q^{(j)}_{r_j-1}x)^2}<(2q_{r_j}^{(j)})^2.
	\]
	Thus
	\[
	(q_{r_j}^{(j)})^2\leq\lambda_j=\inf|\psi'|_{I_j}|\leq\Lambda_j=\sup|\psi'|_{I_j}|\leq (2q_{r_j}^{(j)})^2.
	\]
	\medskip
	
	\noindent Let $a=22$, $s$ the smallest natural number such that $\sizer(1^s)\ge r$, $k=2r$, $\beta_1=1^k$  and, for $2\le j\le k+1$, $\beta_j=1^{k+1-j} a\,1^s=1^{k+1-j}22\;1^s$. Then, $B=\{\beta_1, \beta_2,\dots,\beta_{k+1}\}$ is primitive.
	
	\medskip
	
	\noindent The alphabet $B=\{\beta_1,\beta_2\,\dots,\beta_{k+1}\}$ as above induces a subshift
	\[
	\Sigma(B)=\{(\gamma_i)_{i\in\mathbb{Z}} \ \mid\ \gamma_i\in B\}.
	\]
	
	\Cref{{lem:calc_s}} implies that, for any $\underline{\theta}\in\Sigma(B)$ and every $n\in\mathbb Z$, 
	\[
	\lambda(\sigma^n(\underline{\theta}))<3+e^{-r}.
	\]
	
	Recall that if $\alpha=c_1c_2\cdots c_m$ and $\beta=\beta_1\beta_2\cdots \beta_n$ are finite words, then
	\[
	q_m(\alpha)q_n(\beta)<q_{m+n}(\alpha\beta)<2q_m(\alpha)q_n(\beta).
	\]
	
	The above estimates give $\Lambda_1=\sup|\psi'|_{I(\beta_1)}|\leq 4\left(\frac{1+\sqrt5}2\right)^{2k}$ and, for $2\le j\le k+1$, 
	\begin{align*}
		\Lambda_j&=\sup|\psi'|_{I(\beta_j)}|\leq 8\cdot\left(\frac{1+\sqrt5}2\right)^{2k-2(j-2)}\cdot (10^2\cdot e^{r+1})\\
		&\leq \left(\frac{1+\sqrt5}2\right)^{2k-2(j-2)}\cdot e^{r+8}
	\end{align*}
	
	Thus, from the above lemma and the third author's work \cite{M:geometric_properties_Markov_Lagrange}, we conclude that 
	\[
	d(3+e^{-r})\geq \dimH(\mathsf{m}(\Sigma(B)))= \min\{1, 2\cdot \dimH(K(B))\}\ge 2\tilde d,
	\]
	where $\mathsf{m}(\omega) =\sup_{n \in \ZZ} \lambda(\sigma^n(\omega))$ denotes the Markov value of $\omega \in \Sigma(B)$, and $\tilde d$ is the solution of
	\[
	\left(4\left(\frac{1+\sqrt5}2\right)^{4r}\right)^{-\tilde d}+\sum_{t=0}^{k-1}\left(\left(\frac{1+\sqrt5}2\right)^{2t}\right)^{-\tilde d}\cdot e^{-(r+8)\tilde d}=1.
	\]
	Since $d(3+e^{-r})=\aO\left(\frac{\log r}r\right)$, we also have $\tilde d=\aO\left(\frac{\log r}r\right)=\ao(1)$. The rest of this section is devoted to finding a lower bound for $\tilde{d}$.
	
	Since $\left(\frac{1+\sqrt5}2\right)^4>e^{3/2}$, $\left(4\left(\frac{1+\sqrt5}2\right)^{2k}\right)^{-\tilde d}\le \left(\frac{1+\sqrt5}2\right)^{-4r\tilde d}\le e^{-\frac32 r\tilde d}$, and we have 
	\begin{equation} \label{eq:ast}
		e^{-(r+8)\tilde d}\cdot\frac{1-\left(\frac{1+\sqrt5}2\right)^{-2k\tilde d}}{1-\left(\frac{1+\sqrt5}2\right)^{-2\tilde d}}=1-\aO(e^{-\frac32 r\tilde d}).
	\end{equation}
	In particular, 
	\begin{align*}
		1&\ge e^{-(r+8)\tilde d}\cdot\frac{1-\left(\frac{1+\sqrt5}2\right)^{-4\tilde d}}{1-\left(\frac{1+\sqrt5}2\right)^{-2\tilde d}}=e^{-(r+8)\tilde d}\cdot\left(1+\left(\frac{1+\sqrt5}2\right)^{-2\tilde d}\right)\\
		&\ge 2e^{-(r+8)\tilde d}\cdot\left(\frac{1+\sqrt5}2\right)^{-2\tilde d}\ge 2e^{-(r+9)\tilde d},
	\end{align*}
	and so $\tilde d\ge\frac{\log 2}{r+9}\ge \frac1{2r}$. So we have
	\[
	\left(\frac{1+\sqrt5}2\right)^{-2k\tilde d}=\left(\frac{1+\sqrt5}2\right)^{-4r\tilde d}\le e^{-\frac32 r\tilde d}\le e^{-3/4}<1/2
	\]
	and thus 
	\[
	1\ge e^{-(r+8)\tilde d}\cdot\frac{1-\left(\frac{1+\sqrt5}2\right)^{-2k\tilde d}}{1-\left(\frac{1+\sqrt5}2\right)^{-2\tilde d}}\ge \frac{e^{-(r+8)\tilde d}}{2\left(1-\left(\frac{1+\sqrt5}2\right)^{-2\tilde d}\right)}.
	\]
	Since $\tilde d=\ao(1)$, writing $c_1=\log\frac{3+\sqrt5}2=0.9624\ldots$, we have 
	\[
	\left(\frac{1+\sqrt5}2\right)^{-2\tilde d}=e^{-c_1 \tilde d}=1-c_1 \tilde d+\aO(\tilde d^2),
	\]
	and therefore $1-\left(\frac{1+\sqrt5}2\right)^{-2\tilde d}=c_1\tilde d+\aO(\tilde d^2)=(1+\aO(\tilde d))c_1 \tilde d$. It follows that 
	\[
	1\ge \frac{e^{-(r+8)\tilde d}}{2\left(1-\left(\frac{1+\sqrt5}2\right)^{-2\tilde d}\right)}=\frac{e^{-(r+8)\tilde d}}{(2+\aO(\tilde d))c_1 \tilde d}\ge \frac{e^{-(r+8)\tilde d}}{2\tilde d}
	\]
	and thus $0\ge -(r+8)\tilde d-\log 2-\log\tilde d$. It follows that $-r\tilde d\le \log\tilde d+\aO(1)$, and thus $\left(\frac{1+\sqrt5}2\right)^{-4r\tilde d}\le e^{-\frac32 r\tilde d}=\aO(\tilde d^{3/2})$. From \eqref{eq:ast}, we get 
	\begin{align*}
		1-\aO(\tilde d^{3/2})&=e^{-(r+8)\tilde d}\cdot\frac{1-\left(\frac{1+\sqrt5}2\right)^{-2k\tilde d}}{1-\left(\frac{1+\sqrt5}2\right)^{-2\tilde d}}\\
		&=e^{-(r+8)\tilde d}\cdot\frac{1-\aO(\tilde d^{3/2})}{(1+\aO(\tilde d))c_1 d}=(1+\aO(\tilde d))\frac{e^{-r\tilde d}}{c_1 d},
	\end{align*}
	and thus $\aO(\tilde d^{3/2})=-r\tilde d+\aO(\tilde d)+c_0-\log\tilde d$ and therefore 
	\begin{equation}\label{eq:doble_ast}
		r\tilde d=-\log\tilde d+c_0+\aO(\tilde d)=|\log\tilde d|+c_0+\aO(\tilde d),
	\end{equation}
	where $c_0=-\log c_1=0.03830054\ldots$. 
	
	In particular, $r\tilde d=(1+\aO(1/|\log\tilde d|))|\log \tilde d|=(1+\ao(1))|\log\tilde d|$, and thus $\log\tilde d+\log r=\log|\log\tilde d|+\ao(1)$ and
	\[
	\log r=-\log\tilde d+\log|\log\tilde d|+\ao(1)=(1-\ao(1))|\log\tilde d|.
	\]
	It follows that $|\log\tilde d|=(1+\ao(1))\log r$ and $\log|\log\tilde d|=\log\log r+\ao(1)$, and so
	\[
	\log\tilde d+\log r=\log|\log\tilde d|+\ao(1)=\log\log r+\ao(1)
	\]
	and $|\log\tilde d|=-\log\tilde d=\log r-\log\log r+\ao(1)=\log r (1-(1+\ao(1))\log\log r/\log r)$, which implies $\log|\log\tilde d|=\log \log r-(1+\ao(1))\log\log r/\log r$. 
	
	From $r\tilde d=(1+\aO(1/|\log\tilde d|))|\log\tilde d|$ it follows that
	\begin{align*}
		\log\tilde d+\log r&=\log|\log\tilde d|+\aO\left(\frac{1}{|\log\tilde d|}\right)\\
		&=\log|\log\tilde d|+\aO\left(\frac{1}{|\log r|}\right)\\
		&=\log \log r-\frac{(1+\ao(1))\log\log r}{\log r},
	\end{align*}
	so $|\log\tilde d|=-\log\tilde d=\log r-\log\log r+(1+\ao(1))\log\log r/\log r$ and, from $r\tilde d=|\log\tilde d|+c_0+\aO(\tilde d)=|\log\tilde d|+c_0+\aO(\log r/r)$, we get 
	\begin{align*}
		\tilde d&=\frac{|\log\tilde d|+c_0+\aO(\log r/r)}r \\
		&=\frac{\log r-\log\log r+c_0+(1+\ao(1))\log\log r/\log r}{r}\\
		&>\frac{\log r-\log\log r+c_0}r,
	\end{align*}
	and thus
	\[
	d(3+e^{-r})>2\cdot\frac{\log r-\log\log r+c_0}r.
	\]
	
	We can give a more precise asymptotic expression for $\tilde d$ (and thus for $d(3+e^{-r})$), using the Lambert function $W\colon [e^{-1},+\infty)\to [-1,+\infty)$, which is the inverse function of $f\colon[-1, +\infty)\to [e^{-1},+\infty), f(x)=x e^x$ (which is increasing in the domain $[-1,+\infty)$): let $g\colon (0,+\infty)\to\mathbb R$ given by $g(x)=rx+\log x$. We have $g(\tilde d)=r\tilde d+\log \tilde d=c_0+\aO(\tilde d)$. Let $d_0\in (0,+\infty)$ be the solution of $g(d_0)=c_0$. Since $g'(x)=r+1/x>r$ for every $x\in (0,+\infty)$, and there exists $t$ between $d_0$ and $\tilde d$ such that $|g(\tilde d)-c_0|=|g(\tilde d)-g(d_0)|=|g'(t)(\tilde d-d_0)|\ge r|\tilde d-d_0|$, it follows that
	\[
	|\tilde d-d_0|\le \frac1r |g(\tilde d)-c_0|=\aO(\tilde d/r)=\aO(\log r/r^2)
	\]
	and $\tilde d=d_0+\aO(\log r/r^2)=(1+\aO(1/r))d_0$. On the other hand, since $rd_0+\log d_0=g(d_0)=c_0$, we have $d_0 e^{rd_0}=e^{c_0}$, and so $f(rd_0)=rd_0e^{rd_0}=re^{c_0}$ and thus $rd_0=W(re^{c_0})$, which gives a closed expression for $d_0$: $d_0=\frac1r W(re^{c_0})$, from which we get 
	\[
	\tilde d=\frac{W(re^{c_0})}r+\aO\left(\frac{\log r}{r^2}\right)=\frac{1+\aO(1/r)}r\cdot W(re^{c_0}).
	\]
	(for a detailed discussion on the function $W$, including its asymptotic expansion, we refer the reader to the work of Corless et al.\ \cite{lambertfunction}).
	
	The improved estimates of the previous section (using $T=\lfloor r/5 \rfloor$ in the case of $1^{s_1}221^{s_2}\ldots$) give the same asymptotic expression for $\frac12 d(3+e^{-r})$, so the proof of \Cref{thm:mainstrong} is complete.
	
	\section{The error term is optimal} \label{sec:error_optimal}
	
	In the case $c=1$, the Markov values larger than $3$ are due to two types of ``contradictions'' that we analyze as two separate subcases:
	
	\smallbreak
	
	\noindent\textbf{Case 1:} Words of the form $1^{s_1}221^{2k+1}221^{2k+j}221^{s_2}$, where $2k+1$ is of the order of $\hat{s}_0$, and $s_1, s_2$ are at least $\hat{s}_0-4$. In this case the Markov value associated with the cut $1^{s_1}221^{2k+1}|221^{2k+j}221^{s_2}$ is $3+x$, where 
	\begin{align*}
		x&=[0;1^{2k+1}221^{s_1}\dots]-[0;1^{2k+2+j}221^{s_2}\dots]\\
		&=(1+\ao(1))\frac{2(3\varphi-4)}{3\varphi^4}\left(\frac1{\varphi^{4k}}+\frac{(-1)^j}{\varphi^{4k+2+2j}}\right),
	\end{align*}
	where $\varphi=\frac{1+\sqrt5}2$, and so $x$ belongs to an interval of the type 
	\[
	\frac{2(3\varphi-4)}{3\varphi^{4k+4}}\left[(1+\ao(1))\left(1-\frac1{\varphi^4}\right),(1+\ao(1))\left(1+\frac1{\varphi^2}\right)\right].
	\]
	
	Indeed, we have
	\[
	[0;1^{2k+1}221^{s_1}\dots]=[0;1^{2k+1}22\overline{1}]+\aO(\varphi^{-8k})
	\]
	and
	\[
	[0;1^{2k+2+j}221^{s_2}\dots]=[0;1^{2k+2+j}22\overline{1}]+\aO(\varphi^{-8k}).
	\]
	Moreover, we have 
	\begin{align*}
		[0;1^{n}22\overline{1}]&=\left[0;1^{n},2+\frac1{2+\varphi^{-1}}\right]\\
		&=[0;1^{n},4-\varphi]=\frac{(4-\varphi)\mathrm{F}_n+\mathrm{F}_{n-1}}{(4-\varphi)\mathrm{F}_{n+1}+\mathrm{F}_n}\\
		&=\frac{\mathrm{F}_{n-1}/\mathrm{F}_n+(4-\varphi)}{(4-\varphi)\mathrm{F}_{n-1}/\mathrm{F}_n+5-\varphi}.
	\end{align*}
	On the other hand, the identity $\frac{au+b}{cu+d}-\frac{av+b}{cv+d}=\frac{(ad-bc)(u-v)}{(cu+d)(cv+d)}$ applied for $a=1, b=4-\varphi, c=4-\varphi, d=5-\varphi, u=\mathrm{F}_{2k}/\mathrm{F}_{2k+1}$ and $v=\mathrm{F}_{2k+1+j}/\mathrm{F}_{2k+2+j}$ together with $(cu+d)(cv+d)=(1+\ao(1))(c\varphi^{-1}+d)^2=(1+\ao(1))(3\varphi)^2$ gives 
	\[
	x=(1+\ao(1))\frac{12-6\varphi}{(3\varphi)^2}(v-u)=(1+\ao(1))\frac2{3\varphi^4}(v-u).
	\]
	In order to estimate $v-u$, let us estimate $\mathrm{F}_n/\mathrm{F}_{n+1}-\varphi^{-1}$: we have 
	\begin{align*}
		\frac{\mathrm{F}_n}{\mathrm{F}_{n+1}}-\frac1{\varphi}&=\frac{\varphi^n-(-\varphi^{-1})^n}{\varphi^{n+1}-(-\varphi^{-1})^{n+1}}-\frac1{\varphi}\\
		&=(1+\ao(1))\frac{(-1)^{n+1}(\varphi+\varphi^{-1})\varphi^{-n}}{\varphi^{n+2}}\\
		&=\frac{(-1)^{n+1}(3\varphi-4+\ao(1))}{\varphi^{2n}}.
	\end{align*}
	Using this for $n=2k+1+j, n=2k$ and subtracting, we get the above estimate for $x$.
	
	\medbreak
	
	\noindent\textbf{Case 2:} Words of the form $1^{s_1}221^{2k}221^{2k+3+j}221^{s_2}$, where $2k$ is of the order of $\hat{s}_0$, and $s_1, s_2$ are at least $\hat{s}_0-4$. In this case the Markov value associated with the cut $1^{s_1}221^{2k}22|1^{2k+3+j}221^{s_2}$ is $3+y$, where 
	\begin{align*}
		y&=[0;1^{2k+3+j}221^{s_2}\dots]-[0;1^{2k+2}221^{s_1}\dots]\\
		&=(1+\ao(1))\frac{2(3\varphi-4)}{3\varphi^4}\left(\frac1{\varphi^{4k+2}}+\frac{(-1)^j}{\varphi^{4k+4+2j}}\right).
	\end{align*}
	The proof of this estimate is analogous to the previous one, applying the above estimate of $\mathrm{F}_n/\mathrm{F}_{n+1}-\varphi^{-1}$ for $n=2k+1, n=2k+2+j$ and subtracting. 
	
	Hence, $y$ belongs to an interval of the type 
	\[
	\frac{2(3\varphi-4)}{3\varphi^{4k+6}}\left[(1+\ao(1))\left(1-\frac1{\varphi^4}\right),(1+\ao(1))\left(1+\frac1{\varphi^2}\right)\right].
	\]

	Since
	\[
	1-\frac1{\varphi^4}>0.854>0.528>\left(1+\frac1{\varphi^2}\right)\cdot\frac1{\varphi^2}
	\]
	and
	\[
	1+\frac1{\varphi^2}<1.382<2.236<\left(1-\frac1{\varphi^4}\right)\varphi^2,
	\]
	it follows that, for large $k$, none of these Markov values belong to the interval 
	\[
	[3+x_k,3+y_k]=3+\frac{2(3\varphi-4)}{3\varphi^{4k+4}}[1.382,2.236],
	\]
	whose size is comparable to the value of its endpoints, and so there are no sequences of the type $\ldots 1^{s_1}221^{s_2}221^{s_3}22 \ldots$ with $s_j>3k/2$ for all $j$ whose Markov values belong to $[3+x_k,3+y_k]$. Indeed, we have the same characterization of sequences of this type whose Markov values are smaller than $y_k$ and whose values are smaller than $x_k$:
	for $s=2k$, if $s_j<s$ then $s_j$ is even and $s_{j-1}-s_j, s_{j+1}-s_j\in \{-2,0,2\}$ (and there are no other restrictions).

	Let again $s=2k$ and $T=\lfloor r\log r\rfloor$, where $r=\lfloor|\log y_k|\rfloor$. For each $\tilde T$ with $T/2<\tilde T\le T$, let $M(\tilde T)$ be the number of elements of the set $B(\tilde T)$ of the sequences $1^{s_1}221^{s_2}22\ldots 221^{s_t}22$ with
	\[
	r\cdot (\tilde T-1)<\sizer(1^{s_1}221^{s_2}22\ldots 221^{s_t}22)\le r\cdot\tilde T,
	\]
	$s_j>3s/4$ for every $j\le t$, $s_1, s_t\ge s$ and such that, for each $j\le t$ with $s_j<s$, $s_j$ is even and $s_{j-1}-s_j, s_{j+1}-s_j\in \{-2,0,2\}$. Let $\tilde d=\max\left\{\frac{\log M(\tilde T)}{r\tilde T}\right\}$. Then $d(3+x_k)\ge 2\tilde d$. Indeed, $\mathsf{m}(\Sigma(B(\tilde T))) \subseteq \calM \cap (-\infty,3+x_k)$, where $\mathsf{m}(\omega) = \sup_{n \in \ZZ} \lambda(\sigma^n(\omega))$ denotes the Markov value of $\omega \in \Sigma(B(\tilde T))$.
	
	Let us now give upper estimates: suppose that $w\tilde w$ does not have a factor $1^{s_1}$ satisfying $\sizer(1^{s_1})\ge r-170$, where $r=\lfloor|\log y_k|\rfloor$ and consider an infinite continuation $\theta$ of it contained in $\Sigma(3+e^{-r})\supseteq\Sigma(3+y_k)$. Then the previous discussion provides a continuation $\overline w \in \Sigma^{(Tr)}(3+e^{-r})$ for some $T\in\{10,\lfloor \log^2 r\rfloor\}$ depending on $\tilde w$ such that $w\tilde w\overline w$ is the continuation of $w\tilde w$ in $\theta$, $w\tilde w\overline w\in \Sigma(3+e^{-r},|w\tilde w\overline w|)$, and the number $K$ of these words $w\tilde w\overline w$ satisfies $K\cdot e^{-Trd}<1/r$ for $d=\frac{\log r-\log\log r}r$. 
	
	Suppose now that $w\tilde w$ has a factor $1^{s_1}$ satisfying $\sizer(1^{s_1})\ge r-170$, where
	\[
	r=\lfloor|\log y_k|\rfloor\in \left(s\log\left(\frac{3+\sqrt5}2\right),(s+2)\log\left(\frac{3+\sqrt5}2\right)\right).
	\]
	Let us consider continuations $\overline w \in \Sigma^{(m)}(3+e^{-r})$ for some $r^{3/2}<m\le r\lfloor r\log r\rfloor$ such that $w\tilde w\overline w\in \Sigma(3+e^{-r},|w\tilde w\overline w|)$ and the continuation of $1^{s_1}$ in $w\tilde w\overline w$ is $1^{s_1}221^{s_2}22\ldots 1^{s_t}22$ such that there are at least $r^{15/16}$ values of $i\le t$ with $s_i<s$, and such that $t$ is minimum with this property. Then $s_t>r-3r^{15/16}$. We will introduce a parameter $j$ equal to the number of values of $i\le t$ with $s_i<s$; consider these $j$ values $i_1<i_2<\dots<i_j$ of $i$. We have $j\ge r^{15/16}$. There are at most $\binom{t}{j}<(\frac{et}j)^j$ choices for the set $\{i_t, 1\le t\le j\}$. Since for $i\in \{i_v, 1\le v\le j\}$ we have at most $3$ choices for $s_{i+1}$, and the total number of these choices is at most $3^j$. Together with the number of choices for the set $\{i_t, 1\le t\le j\}$, this gives an estimate of $(\frac{3et}j)^j$ for these choices of the set $\{(i_t,s_{i_t}), 1\le t\le j\}$. Let $\overline t=t-j$. The number of choices of the remaining values of the $s_i$ is at most the number of solutions of $\hat s_1+\hat s_2+\dots+\hat s_{\overline t}\le m/\log\left(\frac{3+\sqrt5}2\right)-s\overline t-j(r-3r^{15/16})\le (U-\overline t)s$, where $U=m/(r-2)-j/2<r\log r$, which is at most $Ue^{U\cdot\frac{\log s}{s}}e^{U\cdot W(s)}$. As before, $e^{W(s)}=(1+\ao(1))s/\log s$, and (since $e^{U\cdot\frac{\log s}{s}}\le e^{\frac{r\log r\log s}{s}}=e^{\aO(\log^2 s)}=e^{\ao(j)}$),
	the total number $\tilde K$ of these sequences is
	\begin{align*}
		&\aO\left(r^4 \left(\frac{6em}{jr}\right)^j (m/r) ((1+\ao(1))s/\log s)^{-j/2}e^{U\cdot\frac{\log s}{s}}
		e^{W(s)m/(r-2)}\right) \\
		{}={}&\aO(s^{-j/4}e^{W(s)m/r}),
	\end{align*}
	and, since $j\ge r^{15/16}$, for $d=\frac{W\left(r/\log\left(\frac{3+\sqrt5}2\right)\right)}r-\frac1{r^{3/2}}$ we have $\tilde K \cdot e^{-md}<e^{-\sqrt{r}}$.
	
	Consider now the remaining case where there are less than $r^{15/16}$ values of $i\le t$ with $s_i<s$ and consider the largest continuation of $1^{s_1}$ in $w\tilde w\overline w\in \Sigma^{(Tr)}(3+e^{-r})$, $T=\lfloor r\log r\rfloor$ of the form $1^{s_1}221^{s_2}22\ldots 1^{s_t}$, $s_j>s-3r^{15/16}$ for each $j$. Taking $j_1$ minimum and $j_2$ maximum with $s_{j_1}, s_{j_2}\ge s$ (notice that $j_1+t-j_2\le r^{15/16}$), the number $\hat N$ of such words is at most $3^{j_1+t-j_2}M<3^{r^{15/16}}M$, where $M$ is the number of elements of $B(\tilde T)$, where $r\cdot(\tilde T-1)<\sizer(1^{s_{j_1+1}}221^{s_2}22\ldots 1^{s_{j_2-1}})\le r\cdot \tilde T$. We have $\tilde T<T-(j_1+t-j_2)/2$ and $M(\tilde T)\le e^{r\tilde d\tilde T}<e^{r\tilde d (T-(j_1+t-j_2)/2)}$, so $\hat N\le e^{r\tilde d T}(3e^{-r\tilde d/2})^{j_1+t-j_2}$. Since, by our lower estimates on $d(3+\varepsilon)$, $\frac{\log \hat N}{Tr}>\frac{\log r-\log\log r+0.03}r$, it follows that $\tilde d\ge \frac{\log M}{r\tilde T}>\frac{\log r-\log\log r}r$, and thus
	\[
	(3e^{-r\tilde d/2})^{j_1+t-j_2}<(3(\log r/r)^{1/2})^{j_1+t-j_2}\le 1
	\]
	and, adding these estimates for all possible choices of $(j_1,t-j_2)$, we get $\hat N\le 2e^{r\tilde d T}$. This, together with the previous estimates, implies that $d(3+y_k)\le 2\tilde d+\aO(1/r^2)$. Indeed, $(e^{-Tr})^{\tilde d+1/r^2}=e^{-T/r}e^{-r\tilde d T}<e^{1-\log r}e^{-r\tilde d T}$, and thus $2e^{r\tilde d T}(e^{-Tr})^{\tilde d+1/r^2}\le 2e^{1-\log r}=2e/r=\ao(1)$.
	
	Finally suppose that $F$ is a twice continuously-differentiable function such that
	\[
	d(3+\varepsilon)=F(\varepsilon)+\ao\left(\frac{\log|\log \varepsilon|}{|\log \varepsilon|^2}\right). 
	\]
	
	By the mean value theorem there is $\xi_k\in(x_k,y_k)$ such that
	\[
	F'(\xi_k)=\frac{F(y_k)-F(x_k)}{y_k-x_k}=\ao\left(\frac{\log|\log y_k|}{y_k|\log y_k|^2}\right).
	\]
	Let $c_1>1$ be a constant we will chose later. By \Cref{thm:mainstrong} we have
	\begin{align*}
		F(c_1y_k)-F(y_k)&=g_1(c_1y_k)-g_1(y_k)+\aO\left(\frac{\log|\log y_k|}{|\log y_k|^2}\right) \\
		&=(2\log(c_1)+\ao(1))\frac{\log|\log y_k|}{|\log y_k|^2}+\aO\left(\frac{\log|\log y_k|}{|\log y_k|^2}\right).
	\end{align*}
	By choosing $c_1>1$ large enough and using the mean value theorem, we obtain $\tilde{\xi_k}\in (y_k,c_1y_k)$ such that
	\[
	F'(\tilde{\xi_k})>C\cdot\frac{\log|\log y_k|}{y_k|\log y_k|^2}
	\]
	Hence for each $k$, we can find a point in $(\xi_k,\tilde{\xi}_k)$ where the second derivative of $F$ is positive and also a point in $(\tilde{\xi}_\ell,\xi_k)$ (for $\ell$ large enough) where the second derivative of $F$ is negative.

	\appendix	
	\section{Basic facts and estimates on continued fractions} 
	
	\label{sec:continued_fractions}
	
	Let $\alpha = c_1 \ldots c_n \in (\NN^*)^n$ be a finite word of length $n > 0$. We define $K(c_1 \ldots c_n)$ to be the \emph{continuant} of $\alpha$, that is, the denominator of the fraction $[0; c_1, \dotsc, c_n]$. The following lemma can be found in the book by Cusick--Flahive \cite[Appendix 2]{Cusick-Flahive}.
	\begin{lemma}[Euler's property of continuants]\label{lem:eulersproperty}
		The continuant $K(c_1 \ldots c_n)$ is equal to a sum of certain products of the integers $c_1, \dotsc, c_n$. Moreover, the products that appear in this sum can be determined in the following way. Start with the product $c_1 \ldots c_n$. Now, include all products obtained by removing pairs of adjacent integers. Continue by including all products obtained by removing two separate pairs of adjacent integers, and follow this procedure until no pair remains. Observe that if $n$ is even, then the empty product, equal to $1$, must be also included.
		
		As a corollary, we obtain that
		\[
		K(c_1 \ldots c_n) {}={} K(c_1 \ldots c_m) K(c_{m+1} \ldots c_n) +  K(c_1 \ldots c_{m-1}) K(c_{m+2} \ldots c_n) 
		\]
		for any $1 \leq m < n$.
	\end{lemma}
	
	In particular, the previous lemma implies that
	\begin{align*}
		K(c_1) &= c_1 \\
		K(c_1c_2) &= c_1 \cdot c_2 + 1 \\
		K(c_1c_2c_3) &= c_1 \cdot c_2 \cdot c_3 + c_1 + c_3 \\
		K(c_1c_2c_3c_4) &= c_1 \cdot c_2 \cdot c_3 \cdot c_4 + c_1 \cdot c_2 + c_1 \cdot c_4 + c_3 \cdot c_4 + 1.
	\end{align*}
	
	Let $\theta=\theta_1 \ldots \theta_n\in(\NN^*)^n$, $a=(2,2)$ and $b=(1,1)$. Using Euler's property of continuants we can find a gap between the size of the intervals of the following words:
	\begin{align}\label{eq:sisesgap}
		\begin{split}
			\sizes(b\theta b)^{-1}&\leq \left(5+\frac{2}{\theta_1}+\frac{2}{\theta_n}\right)^2q_n(\theta)^2, \\
			\sizes(a\theta a)^{-1}&\geq \left(25+\frac{10}{\theta_1+1}+\frac{10}{\theta_n+1}\right)^2q_n(\theta)^2.
		\end{split}
	\end{align}
	Indeed, using the convention $q_{0}=1$ and $q_{-1}=0$ we have
	\begin{align*}
		q_{n+4}(b\theta b)&=
		\begin{multlined}[t]
			4q_n(\theta)+2q_{n-1}(\theta_1 \ldots \theta_{n-1})+2q_{n-1}(\theta_2 \ldots \theta_n) \\
			+q_{n-2}(\theta_2 \ldots \theta_{n-1}) 
		\end{multlined} \\
		&\leq \left(5+\frac{2}{\theta_1}+\frac{2}{\theta_n}\right)q_n(\theta),
	\end{align*}
	\begin{align*}
		q_{n+4}(a\theta a)&=
		\begin{multlined}[t]
			25q_n(\theta)+10q_{n-1}(\theta_1 \ldots \theta_{n-1})+10q_{n-1}(\theta_2 \ldots \theta_n) \\
			+4q_{n-2}(\theta_2 \ldots \theta_{n-1}) 
		\end{multlined} \\
		&\geq\left(25+\frac{10}{\theta_1+1}+\frac{10}{\theta_n+1}\right)q_n(\theta),
	\end{align*}
	and finally we use $q_m(a_1 \ldots a_m)^2\leq\sizes(a_1 \ldots a_m)\leq 2q_m(a_1 \ldots a_m)^2$.

	\begin{lemma}\label{lem:rcomparedwithsize}
		Let $w$ be a nonempty finite word in $1$ and $2$ of length $n \in \NN^*$. We have that
		\[
		(n-3) \log\left( \frac{3 + \sqrt{5}}{2} \right) \leq \sizer(w) \leq (n+1) \log(3 + 2 \sqrt{2}).
		\]
	\end{lemma}
	
	\begin{proof}
		Given $\alpha = c_1 \ldots c_n \in (\NN^*)^n$, we have that
		\[
		\sizes(\alpha) = \frac{1}{q_n(q_n + q_{n-1})},
		\]
		so $\sizes(\alpha)$ is minimized when $q_n$ and $q_{n-1}$ are maximized; and maximized when $q_n$ and $q_{n-1}$ are minimized. This happens, respectively, when $q_n = \mathrm{P}_n$ (where $\mathrm{P}_n$ is the $n$-th Pell number) and where $q_n = \mathrm{F}_n$ (where $\mathrm{F}_n$ is the $n$-th Fibonacci number). Hence,
		\[
		\sizer(1^n) \leq \sizer(w) \leq \sizer(2^n).
		\]
		Moreover, we have that
		\begin{align}
			\sizes(1^n)^{-1} &= \mathrm{F}_{n+1}(\mathrm{F}_{n+1} + \mathrm{F}_{n})\nonumber \\
			&= -\frac{1}{5}(-1)^{n+1} + \frac{\sqrt{5}+1}{10}\left( \frac{3 + \sqrt{5}}{2} \right)^{n+1} - \frac{\sqrt{5}-1}{10}\left(\frac{3 - \sqrt{5}}{2}\right)^{n+1} \nonumber\\
			&\geq -\frac{1}{5} + \frac{\sqrt{5}+1}{10}\left( \frac{3 + \sqrt{5}}{2} \right)^{n+1} - \frac{\sqrt{5}-1}{2} \nonumber\\
			&=\frac{\sqrt{5}+1}{10} \left( \left( \frac{3 + \sqrt{5}}{2} \right)^{n+1} - 1 \right) \geq \left( \frac{3 + \sqrt{5}}{2} \right)^{n-1}, \label{eq:boundfibonacciproduct}
		\end{align}
		and, on the other hand, we have that
		\begin{align*}
			\sizes(2^n)^{-1} &= \mathrm{P}_n(\mathrm{P}_n + \mathrm{P}_{n-1}) = \frac{(3 + 2\sqrt{2})^n - (3 - 2\sqrt{2})^n}{4\sqrt{2}} \\
			&\leq \frac{(3 + 2\sqrt{2})^{n+1}}{4\sqrt{2}} \leq (3 + 2\sqrt{2})^{n+1}.
		\end{align*}
		Thus, we obtain that
		\[
		(n - 1) \log\left(\frac{3 + \sqrt{5}}{2}\right) \leq \log \sizes(w)^{-1} \leq (n+1) \log(3 + 2\sqrt{2}).
		\]
		Finally, since $2\log\left(\frac{3 + \sqrt{5}}{2}\right) > 1$, we get that
		\[
		(n - 3) \log\left(\frac{3 + \sqrt{5}}{2}\right) \leq \sizer(w) = \lfloor\log \sizes(w)^{-1}\rfloor \leq (n+1) \log(3 + 2\sqrt{2}).
		\]

	\end{proof}

	\begin{lemma} \label{lem:bounds_subwords}
		Let $w$ be a finite word and let $v$ be a factor of $w$. Then, $\sizes(w) \leq \sizes(v)$ and $\sizer(w) \geq \sizer(v)$.
	\end{lemma}
	\begin{proof}
		Assume first that $v$ is a prefix of $w$, so $w = v \beta$ for some word $\beta$. Then, $\sizes(w) = \sizes(v \beta) = |I(v \beta)| \leq |I(v)| = \sizes(v)$, since, by definition, $I(v \beta) \subseteq I(v)$.
		
		Assume now that $w = \alpha v \beta$ for some words $\alpha, \beta$, where $\alpha$ is nonempty. Then, $\sizes(w) = \sizes(\alpha v \beta) \leq \sizes(\alpha v) < 2 \sizes(\alpha) \sizes(v)$. Moreover, if $\alpha$ starts with the letter $c$, then we have that $\sizes(\alpha) \leq \sizes(c)$. Since $\sizes(c) = 1/(c^2+c)$, we have that $\sizes(c) \leq 1/2$. We obtain that $\sizes(w) < \sizes(v)$, as desired.
	\end{proof}
	
	A property that is useful to simplify some computations is
	\[
	\sizer(w_1 k_1k_2w_2)\geq \sizer(w_1)+\sizer(w_2)
	\]
	for any positive integers such that $(k_1,k_2)\neq (1,1)$ and any words $w_1,w_2$. Indeed, it follows from
	\[
	\sizes(w_1 k_1k_2w_2)\leq 4\sizes(k_1k_2)\sizes(w_1)\sizes(w_2)\leq \sizes(w_1)\sizes(w_2)/3.
	\]
	For $(k_1,k_2)=(1,1)$ we have that $\sizer(w_1bw_2)\geq \sizer(w_1)+\sizer(w_2)-1$, since $\sizer(b)=1$.

	Nevertheless, we will prove some sharper bounds that we will use to get cleaner statements of the lemmas. 
	
	Let $s_1,\dots,s_\ell$ be nonnegative integers with $\ell\geq 2$. We will show that
	\begin{equation}\label{eq:intervalswith1s}
		\sizer(1^{s_1}221^{s_2}22\ldots 221^{s_\ell})\geq (s_1+\dots+s_\ell+3(\ell-2))\log\left(\frac{3+\sqrt{5}}{2}\right).
	\end{equation}
	and 
	\begin{equation}\label{eq:intervalswith2s}
		\sizer(2^{s_1}112^{s_2}11\ldots 112^{s_\ell})\geq (s_1+\dots+s_\ell+\ell-2)\log(3+2\sqrt{2}).
	\end{equation}
	
	First, we will show inductively that
	\begin{equation}\label{eq:convergents}
		q(1^{s_1}221^{s_2}22\ldots 221^{s_\ell})\geq \mathrm{F}_{s_1+\dots+s_\ell+3(\ell-1)+1},
	\end{equation}
	and
	\begin{equation}\label{eq:convergents2}
		q(2^{s_1}112^{s_2}11\ldots 112^{s_\ell})\geq \mathrm{P}_{s_1+\dots+s_\ell+\ell}.
	\end{equation}
	Using Euler's property of continuants (\Cref{lem:eulersproperty})
	\[
	q(1^{s_1}221^{s_2})=q(1^{s_1})q(221^{s_2})+q(1^{s_1-1})q(21^{s_2}).
	\]
	Since $q(1^s)=\mathrm{F}_{s+1}$ one has
	\[
	q(221^{s})=5q(1^{s})+2q(1^{s-1})=5\mathrm{F}_{s+1}+2\mathrm{F}_{s}=3\mathrm{F}_{s+1}+2\mathrm{F}_{s+2},
	\]
	\[
	q(21^{s})=2q(1^{s})+q(1^{s-1})=2\mathrm{F}_{s+1}+\mathrm{F}_{s}.
	\]
	From the identity
	\[
	\mathrm{F}_n\mathrm{F}_m+\mathrm{F}_{n-1}\mathrm{F}_{m-1}=\mathrm{F}_{n+m-1},
	\]
	we get
	\begin{align}
		\mathrm{F}_{n+1}q(221^m)+\mathrm{F}_nq(21^m)&=\mathrm{F}_{n+1}(2\mathrm{F}_{m+2}+3\mathrm{F}_{m+1})+\mathrm{F}_{n}(2\mathrm{F}_{m+1}+\mathrm{F}_{m})\nonumber \\
		&=2\mathrm{F}_{n+m+2}+\mathrm{F}_{n+m+1}+2\mathrm{F}_{n+1}\mathrm{F}_{m+1} \nonumber \\
		&=\mathrm{F}_{n+m+4}+2\mathrm{F}_{n+1}\mathrm{F}_{m+1}. \label{eq:fibeq}
	\end{align}
	Thus
	\begin{align*}
		q(1^{s_1}221^{s_2})
		&=q(1^{s_1})q(221^{s_2})+q(1^{s_1-1})q(21^{s_2}) \\
		&=\mathrm{F}_{s_1+1}q(221^{s_2})+\mathrm{F}_{s_1}q(21^{s_2}) \\
		&=\mathrm{F}_{s_1+s_2+4}+2\mathrm{F}_{s_1+1}\mathrm{F}_{s_2+1}.
	\end{align*}
	Hence \eqref{eq:convergents} is true for $\ell=2$. Assuming it for $\ell$, we use \eqref{eq:fibeq} with $n=s_1+\dots+s_\ell+3(\ell-1)+1$ and $m=s_{\ell+1}$ to obtain
	\begin{align*}
		q(1^{s_1}221^{s_2}22\ldots 221^{s_{\ell+1}})&=
		\begin{multlined}[t]
			q(1^{s_1}221^{s_2}22\ldots 221^{s_\ell})q(221^{s_{\ell+1}})+ \\
			q(1^{s_1}221^{s_2}22\ldots 221^{s_\ell-1})q(21^{s_{\ell+1}})
		\end{multlined} \\
		&\geq \mathrm{F}_{n}q(221^{s_{\ell+1}})+\mathrm{F}_{n-1}q(21^{s_{\ell+1}}) \\
		&\geq \mathrm{F}_{s_1+\dots+s_{\ell+1}+3\ell+1}
	\end{align*}
	Finally, using \eqref{eq:boundfibonacciproduct}
	\begin{alignat*}{2}
		\sizes(1^{s_1}221^{s_2}22\ldots 221^{s_\ell})^{-1}&{}\geq{} & &\mathrm{F}_{s_1+\dots+s_\ell+3(\ell-1)+1} \\
		&&&(\mathrm{F}_{s_1+\dots+s_\ell+3(\ell-1)+1}+\mathrm{F}_{s_1+\dots+s_\ell+3(\ell-1)}) \\
		&{}\geq{}&&\left( \frac{3 + \sqrt{5}}{2} \right)^{s_1+\dots+s_\ell+3(\ell-1)-1}.
	\end{alignat*}
	
	On the other hand, using that $\mathrm{F}_{n+2}\leq 3\mathrm{F}_n$ we get
	\begin{align*}
		\sizes(1^n)^{-1}=\mathrm{F}_{n+1}\mathrm{F}_{n+2}\leq \frac{3}{4}\mathrm{F}_{2n+2}\leq\left( \frac{3 + \sqrt{5}}{2} \right)^{n}
	\end{align*}
	so
	\begin{equation}\label{eq:rofpower1}
		\sizer(1^n)\leq n\log((3+\sqrt{5})/2).
	\end{equation}
	
	Similarly, one has that $q(2^s)=\mathrm{P}_{s+1}$ and $q(112^s)=\mathrm{P}_{s+2}$. The Pell numbers also satisfy the identity
	\[
	\mathrm{P}_n \mathrm{P}_m + \mathrm{P}_{n-1} \mathrm{P}_{m-1}=\mathrm{P}_{n+m-1}.
	\]
	Hence
	\begin{align*}
		\mathrm{P}_{n+1}q(112^m)+\mathrm{P}_nq(12^m)&=\mathrm{P}_{n+1}\mathrm{P}_{m+2}+\mathrm{P}_n(\mathrm{P}_{m}+\mathrm{P}_{m+1}) \\
		&=\mathrm{P}_{n+m+2}+\mathrm{P}_n\mathrm{P}_m
	\end{align*}
	Therefore by induction
	\begin{align*}
		q(2^{s_1}112^{s_2}11\ldots 112^{s_{\ell+1}})&\geq \mathrm{P}_{s_1+\dots+s_\ell+\ell}q(112^{s_{\ell+1}})+\mathrm{P}_{s_1+\dots+s_\ell+\ell-1}q(12^{s_{\ell+1}}) \\
		&\geq \mathrm{P}_{s_1+\dots+s_\ell+s_{\ell+1}+(\ell+1)}
	\end{align*}
	
	Finally to show \eqref{eq:intervalswith2s} we use that
	\begin{align*}
		\sizes(2^{s_1}112^{s_2}11\ldots 112^{s_{\ell}})^{-1}\geq \mathrm{P}_n^2 \geq 4(3+2\sqrt{2})^{n-2}
	\end{align*}
	where $n=s_1+\dots+s_\ell+\ell$.

	\begin{lemma}\label{prop:rtranspose}
		If $\alpha=c_1 \ldots c_n \in(\NN^\ast)^n$ then
		\[
		\frac{[1;c_n+1]}{[1;c_1]}\leq \frac{\sizes(\alpha^\ast)}{\sizes(\alpha)}\leq \frac{[1;c_n]}{[1;c_1+1]}
		\]
		and, hence,
		\[
		-\log\left(1+\frac{1}{c_n+1}\right)-1\leq \sizer(\alpha)-\sizer(\alpha^\ast)\leq\log\left(1+\frac{1}{c_n}\right)+1.
		\]
	\end{lemma}
	\begin{proof}
		By Euler's property of continuants (\Cref{lem:eulersproperty}) we have $q_n(c_1\ldots c_n)=q_n(c_n\ldots c_1)$ thus
		\begin{align*}
			\sizes(\alpha^\ast)^{-1}&=q_n(c_n\ldots c_1)(q_n(c_n\ldots c_1)+q_{n-1}(c_n\ldots c_2)) \\
			&\geq \left(1+\frac{1}{c_1+1}\right)q_n(c_n\ldots c_1)^2 \\
			&= \left(1+\frac{1}{c_1+1}\right)q_n(c_1\ldots c_n)^2,
		\end{align*}
		while 
		\begin{align*}
			\sizes(\alpha)^{-1}&=q_n(c_1\ldots c_n)(q_n(c_1\ldots c_n)+q_{n-1}(c_1\ldots c_{n-1})) \\
			&\leq\left(1+\frac{1}{c_n}\right)q_n(c_1\ldots c_n)^2.
		\end{align*}
		Hence
		\begin{align*}
			\frac{\sizes(\alpha^\ast)}{\sizes(\alpha)}\leq\frac{[1;c_n]}{[1;c_1+1]}
		\end{align*}
		By symmetry we obtain the lower bound.
	\end{proof}
	
	\sloppy\printbibliography
	
\end{document}